\documentclass[12pt,oneside,reqno]{amsart}
\usepackage{}
\usepackage{amssymb}
\usepackage[colorlinks=true, linkcolor=blue, citecolor=blue]{hyperref}
\usepackage{bbm}
\usepackage{cases}
\usepackage{amsmath}
\usepackage{graphicx}
\usepackage{mathrsfs}
\usepackage{stmaryrd}
\usepackage{color}
\usepackage{soul}
\usepackage[dvipsnames]{xcolor}
\usepackage{amsfonts}
\usepackage{cite}
\usepackage{enumerate,amsmath,amssymb,amsthm}

\numberwithin{equation}{section}

\newcommand{\be}{\begin{eqnarray}}
\newcommand{\mE}{\end{eqnarray}}
\newcommand{\ce}{\begin{eqnarray*}}
\newcommand{\de}{\end{eqnarray*}}
\newtheorem{theorem}{Theorem}[section]
\newtheorem{lemma}[theorem]{Lemma}
\newtheorem{remark}[theorem]{Remark}
\newtheorem{definition}[theorem]{Definition}
\newtheorem{proposition}[theorem]{Proposition}
\newtheorem{example}[theorem]{Example}
\newtheorem{corollary}[theorem]{Corollary}

\def\e{{\mathrm{e}}}
\def\eps{\varepsilon}

\def\a{\alpha}

\def\b{\beta}
\def\p{\partial}
\def\d{\delta}

\def\l{\lambda}

\def\[{{\Big[}}
\def\]{{\Big]}}
\def\<{{\langle}}
\def\>{{\rangle}}
\def\({{\Big(}}
\def\){{\Big)}}

\def\bx{{\mathbf{x}}}

\def\dif{{\mathord{{\rm d}}}}

\def\no{\nonumber}
\def\={&\!\!=\!\!&}
\def\bt{\begin{theorem}}
\def\et{\end{theorem}}
\def\bl{\begin{lemma}}
\def\el{\end{lemma}}
\def\br{\begin{remark}}
\def\er{\end{remark}}

\def\bd{\begin{definition}}
\def\ed{\end{definition}}
\def\bp{\begin{proposition}}
\def\ep{\end{proposition}}
\def\bc{\begin{corollary}}
\def\ec{\end{corollary}}
\def\bx{\begin{example}}
\def\ex{\end{example}}

\def\cG{{\mathcal G}}

\def\cI{{\mathcal I}}

\def\cL{{\mathcal L}}

\def\cP{{\mathcal P}}

\def\cU{{\mathcal U}}
\def\cV{{\mathcal V}}
\def\cW{{\mathcal W}}

\def\mC{{\mathbb C}}

\def\mE{{\mathbb E}}
\def\mF{{\mathbb F}}

\def\mI{{\mathbb I}}

\def\mN{{\mathbb N}}

\def\mP{{\mathbb P}}

\def\mR{{\mathbb R}}

\def\sC{{\mathscr C}}
\def\sD{{\mathscr D}}

\def\sF{{\mathscr F}}

\def\sJ{{\mathscr J}}
\def\sK{{\mathscr K}}
\def\sL{{\mathscr L}}

\def\sP{{\mathscr P}}

\def\geq{\geqslant}
\def\leq{\leqslant}

\def\T{\mathord{{\rm Tr}}}

\allowdisplaybreaks

\begin{document}

\title{Asymptotic limit of fully coupled multi-scale non-linear stochastic system: the non-autonomous approximation method}

\date{}

\author{Yuewen Hou, Yun Li  and Longjie Xie}

\address{Yuewen Hou: School of Mathematics and Statistics, Fujian Normal University,
Fuzhou, Fujian 350007, P.R. China\\ Email: hyw.fjnu@foxmail.com
}

\address{Yun Li:
       School of Mathematics and Statistics, Jiangsu Normal University,
	Xuzhou, Jiangsu 221000, P.R. China\\
       Email: liyun@jsnu.edu.cn
}

\address{Longjie Xie:
	School of Mathematics and Statistics, Jiangsu Normal University,
	Xuzhou, Jiangsu 221000, P.R. China\\
	Email: longjiexie@jsnu.edu.cn
}

\thanks{
	This work is supported by the National Key R$\&$D program of China (No. 2023YFA1010103), NNSF of China (No.  12471140, 12090011, 12401580) and NSF of Jiangsu Province (No. BK20241047).
}


\begin{abstract}
In this paper, we develop a novel argument, the non-autonomous approximation method, to seek the asymptotic limits of the fully coupled multi-scale McKean-Vlasov stochastic systems with irregular coefficients, which, as summarized in \cite[Section 7]{BS2},  remains an open problem in the field.  We provide an explicit characterization for the averaged limit of the non-linear stochastic system, where both the choice of the frozen equation and the    definition of the averaged  coefficients  are more or less unexpected  since  new integral terms with respect to the measure variable appear. More importantly,
in contrast with the classical theory of multi-scale systems which focuses on the averaged limit of the slow process, we propose a new perspective that  the asymptotic behavior of the entire system is actually governed by the limit of the fast motion. By studying the long-time estimates of the solution of the Kolmogorov equation in Wasserstein space, we identify  the limiting distribution of the  fast motion of the non-linear system, which, to the best of our knowledge,  is   new even for the classical multi-scale It\^o SDEs. Furthermore,
rates of convergence are  also obtained, which are rather sharp and depend only on the regularity of the coefficients with respect to the slow variable.
The  innovation of our argument is to transform the non-linear system into a sequence of linear but non-autonomous systems, which is rather simple insofar as it avoids to involve the mean-field type PDEs associated with non-linear stochastic  system, and at the same time, it turns out to be quite effective as it enables us to show that the strong convergence in the averaging principle of the non-linear stochastic system follows directly from the weak convergence, which significantly simplified the proof.

	\bigskip
	
	
\end{abstract}

\maketitle

\tableofcontents

\section{Introduction  and main results}

In this paper, we  consider the following fully coupled multi-scale McKean-Vlasov stochastic system in $\mR^{d_1}\times\mR^{d_2}$:
\begin{equation} \label{sde1}
\left\{ \begin{aligned}
&\dif X^{\eps}_t =b(X^{\eps}_t,\cL_{X_t^\eps},Y^{\eps}_t,\cL_{Y_t^\eps})\dif t
+\sigma(X^{\eps}_t,\cL_{X_t^\eps},\cL_{Y_t^\eps})\dif W^1_t,\qquad\qquad\qquad\, X^{\eps}_0=\xi,\\
&\dif Y^{\eps}_t =\frac{1}{\eps}F(X^{\eps}_t,\cL_{X_t^\eps},Y^{\eps}_t,\cL_{Y_t^\eps})\dif t+\frac{1}{\sqrt{\eps}}G(X^{\eps}_t,\cL_{X_t^\eps},Y^{\eps}_t,\cL_{Y_t^\eps})\dif W_t^2,\qquad Y^{\eps}_0=\eta,
\end{aligned} \right.
\end{equation}
where $d_1, d_2\geq 1$, $\xi$, $\eta$ are two random variables, $b,\sigma,F,G$ are measurable functions, $W^1_t$, $W^2_t$ are $d_1$, $d_2$-dimensional independent standard Brownian motions both defined on some probability space $(\Omega,\sF,\mP)$, respectively, and the small parameter $0<\eps\ll 1$ represents the separation of time scales between the slow component $X_t^\eps$  and the fast motion $Y_t^\eps$. Here and throughout, we denote by  $\cL_X$ the distribution of  a random
variable $X$, and
$\sP_2(\mR^{d})$ ($d\geq 1$) the space of all square integrable probability measures over $\mR^d$ equipped with the Wasserstein metric, i.e.,
$$
\cW_2(\mu_1,\mu_2):=\inf_{\pi\in\cP(\mu_1,\mu_2)}\left(\int_{\mR^d\times\mR^d}|x-y|^2\pi(\dif x,\dif y)\right)^{\frac{1}{2}},\quad \forall \mu_1,\mu_2\in \sP_2(\mR^{d}),
$$
where $\cP(\mu_1,\mu_2)$ is the class of measures on $\mR^d\times\mR^d$ with marginals $\mu_1$ and $\mu_2$.

\vspace{2mm}
The McKean-Vlasov  stochastic differential equation (SDE for short), also known as  the mean-field  SDE or the distribution dependent SDE, describes the limiting behavior of an individual particle involving within a system of particles interacting through their empirical measure, as the size of the
population grows to infinity (the so-called
propagation of chaos, see e.g. \cite{S0}), and the  solution is often called the non-linear diffusion. The pioneer work on such system was indicated by Kac \cite{K} in kinetic theory and  McKean \cite{M} in the study of non-linear partial differential equations (PDEs for short). So far, the McKean-Vlasov SDEs have been investigated in various aspects such as well-posedness, ergodicity, large deviation and connection with non-linear Fokker-Planck equations as well as porous media and granular flows, etc, we refer the readers to \cite{BR,CDLL,CD,CGPS,CGM,CM,EGZ,GS,HSS,MV,O,RST,TT,W1} and the references therein.
Meanwhile, the presence of multiple scales arises naturally in many applications  ranging from climate modeling to chemical physics, and has been the central
topic of study in science and engineering, see the monograph \cite{PS}, and \cite{C2,HL,HS,K1,KK,KY} among others. In particular,  multiple scales can leads to hysteresis loops in the bifurcation diagram and induce phase transitions of certain McKean-Vlasov equations as studied in \cite{CGPS,DGPS,GP}, and the asymptotic limit of the system (\ref{sde1}) as $\eps\to0$ is closely related to  the limit theorem for solutions of PDEs  with singularly perturbed terms in Wasserstein space, which has its own interest, see e.g. \cite{HP} and \cite[Chapter IV]{Fr}.
Averaging results for multi-scale McKean-Vlasov SDEs can be found in, see e.g. \cite{BS1, BS2, HLL2,RSX,Z} and many others. However, in
all the previous works the coefficients of the multi-scale systems are not allowed to depend on the distribution of the fast motion.
A system of weakly interacting diffusions in a two-scale potential relying on the faster empirical measure was considered in \cite{DGP}, the combined mean field and diffusive limits were investigated. Recently, the authors in \cite{LWX} considered the diffusion approximation for the multi-scale McKean-Vlasov SDEs by using a non-linear PDE as the corrector, where the coefficients can depend on the distributions of both the slow component and the fast motion, but the coefficients in the fast motion are  not allowed to depend on the slow component itself. We shall explain in subsection 1.1 that the  fully cross interactions between the slow and fast modes as well as their distributions will made the system (\ref{sde1}) totally non-linear and more difficult to deal with. So far in the literature and as summarized in  \cite[Section 7]{BS2}, the averaging principle for the fully coupled McKean-Vlasov stochastic system (\ref{sde1}) {\bf remains to be an open problem}.

\vspace{2mm}
We also point out that two methods are commonly used in the literature to study  the asymptotic limit of multi-scale systems as $\eps\to0$: the Khasminskii's time discretisation  argument and the technique of Poisson equation. But both seem to be {\bf not feasible for the fully coupled McKean-Vlasov system (\ref{sde1})}. On the one hand, the Khasminskii's time discretisation  argument essentially requires the Markov property of the corresponding frozen equation, which does not hold for the non-linear system (\ref{sde1}) anymore since its frozen equation will be a McKean-Vlasov type equation (the solution does  not define a flow) due to the dependence of the distribution of the fast motion. On the other hand, for the fully coupled stochastic system, even if the Poisson equation in the Wasserstein space associated with the non-linear system was studied in \cite{LWX} and the regularities of the corresponding solution are obtained therein, there is still an essential problem in using the Poisson equation to prove the averaging principle as explained in \cite[Remark 3.2]{LWX}.

\vspace{2mm}
The aim in this paper is to develop a novel and robust method (see Section 2 for a brief introduction of the
main idea) to investigate  the asymptotic behavior of the entire McKean-Vlasov system (\ref{sde1})   as $\eps\to0$.
More precisely, we shall identify the averaged limit of the non-linear stochastic system (\ref{sde1}), and establish not only the strong and weak convergence in the averaging principle for the slow process $X^{\eps}_t$, but also characterize the limiting distribution of the fast motion $Y^{\eps}_t$ (which seems to be totally new even for the classical multi-scale It\^o SDEs). The main results are presented in {\bf Theorem \ref{main1}} below.
The  innovation of our argument is to transform the non-linear system (\ref{sde1}) into a sequence of linear but non-autonomous systems, which is rather simple insofar as it avoids to involve the  mean-field type PDEs associated with non-linear system (\ref{sde1}) (the mean-field type backward and forward  Kolmogorov equations as well as the Poisson equation,  whose optimal regularities of the solutions still seems to be unknown), and the corresponding frozen equations of the non-autonomous   systems turn out to be the autonomous approximations of the frozen equation of the original non-linear system.
We believe that our approach goes beyond the scope of the results established here and may be of interest for applications to homogenization of non-linear equations (see e.g. \cite{HP,PV1,PV2}).
{\bf Besides the novelty of the  method, the main highlights of our work  can be summarized as follows} (see also Remark \ref{r1} below):

\vspace{2mm}

{\bf (i) Identification of the averaged limit for the non-linear system (\ref{sde1}).}
As concluded  in \cite[Section 7]{BS2}, the averaging principle for the fully coupled McKean-Vlasov SDE (\ref{sde1}) remains to be an open problem in the field. We will explain in subsection 1.1 what we  would expect the averaged system for (\ref{sde1}) to be by adopting the  conventional intuitive derivation used in the previous literature, but such formal derivation will lead to a wrong limit. We shall identify the averaged equation for the fully coupled non-linear  system (\ref{sde1}), and prove both the strong convergence and the convergence of the distribution of the slow process to its averaged limit. Compared with the existing results  (see e.g. \cite{BS1, BS2, HLL2,LWX,RSX,Z}), it turns out that both the choice of the frozen equation and the obtained averaged  limit for the  non-linear system (\ref{sde1}) are more or less unexpected  since  new integral terms with respect to measure variable appear, which are exactly due to the dependence on the distribution of the fast motion in system (\ref{sde1}), see also {\it (ii)} of Remark \ref{br} below for partial  explanation why  the limit we derived is reasonable. We also point out that the convergence of the distribution   of the slow process we obtained is even more general than the classical weak convergence of  multi-scale SDEs, see Remark \ref{r1} {\it (ii)} below.
Besides, counter example is known which shows that the strong convergence in the averaging principle of classical multi-scale It\^o SDEs does not hold when the diffusion coefficient $\sigma$ in the slow process relies on the fast motion $Y_t^\eps$, see e.g. \cite{L} (this is the reason for our choice to focus on system (1.1) with $\sigma$ being independent of the  $y$-variable). However, our study demonstrate that    despite the coefficient $\sigma$ depends on the distribution $\sL_{Y_t^\eps}$ of the fast motion, the strong convergence still holds.

\vspace{2mm}

{\bf (ii) Characterization of the limit of the fast motion.}
The classical theory of averaging principle for   multi-scale systems   focuses on seeking the limit  of the slow component, which can be thought of as the mathematical model for a phenomenon appearing at the natural time scale, and the fast motion  is referred to as the  random environment taking place at a faster time scale. Obviously, characterization of the limit of the fast motion should be more difficult. But we shall provide a new perspective  that identifying the limit of the fast motion is more important, and  the asymptotic behavior of the whole multi-scale system is in fact governed by the limit of the distribution of the  fast motion, whereas the averaged limit of the slow process follows directly as a byproduct, see subsection 1.2  for more detailed explanation. We give explicit characterization for the limit of the distribution of the  fast motion $Y^\eps_t$ in the McKean-Vlasov system (\ref{sde1}).
To the best of our knowledge, this is new even for the classical multi-scale It\^{o} SDEs. The proof of the convergence of the distribution of the fast motion relies on a new tool:  the long time estimates of the solution of Kolmogorov equation on the product measure space $\sP_2(\mR^{d_1}\times\mR^{d_2})$, see equation (\ref{cp0}) in subsection 3.2. In a very particular case where the fast motion in the non-linear system (\ref{sde1}) does not depend on the slow process and its distribution, our arguments also provide  an autonomous approximation method to prove the  exponential ergodic in the weighted total variation distance for the McKean-Vlasov SDEs, see Remark \ref{r1} {\it (iii)}, which should be of independent interest.

\vspace{2mm}
{\bf (iii) Derivation of strong convergence from weak convergence.}
Usually, the strong convergence in the averaging principle of classical multi-scale SDEs implies the weak convergence (but  the weak convergence require weaker assumptions on the coefficients, and the weak convergence rate is faster than the strong convergence). However, for the McKean-Vlasov stochastic system (\ref{sde1}), we find that it is enough to prove  the weak convergence (i.e., the convergence of the distributions of  the slow process and the fast motion), and we shall show that the strong convergence in the averaging principle follows directly from the weak convergence  (which is a significant distinction from the classical theory of multi-scale It\^o SDEs) by reviewing the non-linear system (\ref{sde1}) as a  linear but non-autonomous system, and the optimal strong convergence rate can be obtained simultaneously. This newfound perspective allows for a more straightforward treatment of the non-linear stochastic systems with irregular coefficients and aviods to involve the Zvonkin's transformation, making the proof  significantly simplified.

\vspace{2mm}

{\bf (iv) Treatment of irregular coefficients and rates of convergence.} All the existing results concerning the asymptotic behavior for multi-scale McKean-Vlasov systems require very strong regularity assumptions on  the coefficients, even if the system does not involve the distribution of the fast motion, see e.g. \cite{BS1, BS2, HLL2,LWX,RSX,Z}.
We assume only  H\"{o}lder continuity of the coefficients in both the space and the measure variables, where the H\"{o}lder continuous with respect to the measure component being for the Wasserstein distance (which is not Lions differentiable, see Remark \ref{r1} {\it (i)} below). This reflects the regularization of noises on the multi-scale non-linear   system (\ref{sde1}). Besides, we obtain  the strong and weak convergence rates in the averaging principle (which are rather sharp and coincide with the cases of classical It\^o SDEs)  as well as the rate of convergence for the distribution of the fast motion. These rates  depend only on the regularities of the coefficients with respect to the slow variable, and do not rely on their regularities with respect to the fast component.
For these, we need to study the optimal regularities for the solutions of two kinds of Kolmogorov equation on Wasserstein space, and introduce an mollifying argument for functions on Wasserstein space  with explicit approximating rate and bounds on the Lions derivatives of the approximation sequence, which might be of independent interest, see e.g. \cite[Section 3]{MZ}.

\subsection{Formal derivation leads to wrong limit}

Let us briefly explain  what we would expect to arise from (\ref{sde1}) as $\eps\to0$ by adopting the formal idea used in the previous method. Meanwhile, we point out the key difference between the fully coupled stochastic system (\ref{sde1}) and the existing results in the literature.

\vspace{2mm}
{\it (i) Intuitive derivation of the averaged limit for multi-scale SDEs.} For simplicity, let us consider
\begin{equation} \label{sde01}
\left\{ \begin{aligned}
&\dif X^{\eps}_t =b(X^{\eps}_t,Y^{\eps}_t)\dif t
+\dif W^1_t,\qquad\qquad\qquad X^{\eps}_0=x\in\mR^{d_1},\\
&\dif Y^{\eps}_t =\frac{1}{\eps}F(X^{\eps}_t,Y^{\eps}_t)\dif t+\frac{1}{\sqrt{\eps}}\dif W_t^2,\qquad\qquad\!\!   Y^{\eps}_0=y\in\mR^{d_2}.
\end{aligned} \right.
\end{equation}
The intuitive idea for deriving the averaged limit equation of the system (\ref{sde01}) is based on the observation
that during the fast transients, the slow variable
remains ``constant", and by the time its changes become noticeable, the fast variable has almost reached its
``quasi-steady state". More explicitly, let us first look at the fast process $Y_t^\eps$. The natural way is to slow it down by making the   time scaling that $t\mapsto \eps t$. Namely, define  $\tilde Y_t^\eps:=Y_{\eps t}^\eps$, then the process  $\tilde Y_t^\eps$ should satisfy
\begin{align}\label{sde002}
\dif \tilde Y_t^\eps=F(X^{\eps}_{\eps t},\tilde Y^{\eps}_{t})\dif t+\dif \tilde W^2_t,\quad \tilde Y^{\eps}_0=y\in\mR^{d_2},
\end{align}
where $\tilde W_t^2:=\eps^{-1/2}W^2_{\eps t}$ is a Brownian motion.
Since we are interested in the limit that $\eps\to0$, it is natural to consider the auxiliary process $Y_t^{\bar x}$ which is the solution of the following frozen equation:
\begin{align}\label{sde02}
\dif Y_t^{\bar x}=F({\bar x},Y_t^{\bar x})\dif t+\dif W_t,\quad Y_0^{{\bar x},\mu}=y\in\mR^{d_2},
\end{align}
where ${\bar x}\in\mR^{d_1}$ is a parameter and $W_t$ is a new standard Brownian motion.
Under certain dissipative condition, the process $Y_t^{\bar x}$ admits a unique invariant measure $\zeta^{\bar x}(\dif y)$.
Taking this back into the slow equation of the system (\ref{sde01}) and  averaging the coefficient with respect to parameter in the fast variable, we obtain that the slow component $X_t^\eps$ will converge  as $\eps\to0$ to the solution of the following averaged equation:
\begin{align}\label{sde03}
\dif \bar X_t=\bar b(\bar X_t)\dif t+\dif W_t^1,\quad\bar X_0=x\in\mR^{d_1},
\end{align}
where the new   drift  is defined by
$$
\bar b(x):=\int_{\mR^{d_2}}b(x,y)\zeta^{x}(\dif y).
$$

{\it (ii) McKean-Vlasov system without involving the distribution of the fast motion.} The above intuitive derivation is still suitable for the following McKean-Vlasov stochastic system:
\begin{equation} \label{sde04}
\left\{ \begin{aligned}
&\dif X^{\eps}_t =b(X^{\eps}_t,\cL_{X_t^\eps},Y^{\eps}_t)\dif t
+\dif W^1_t,\qquad\qquad\qquad X^{\eps}_0=\xi,\\
&\dif Y^{\eps}_t =\frac{1}{\eps}F(X^{\eps}_t,\cL_{X_t^\eps},Y^{\eps}_t)\dif t+\frac{1}{\sqrt{\eps}}\dif W_t^2,\qquad\qquad\!\!   Y^{\eps}_0=\eta.
\end{aligned} \right.
\end{equation}
Note that the coefficients do not depend on the distribution of the fast motion. Again, with the time scaling  $t\mapsto \eps t$, we have that  the re-scaled fast process  $\tilde Y_t^\eps:=Y_{\eps t}^\eps$ satisfies
$$
\dif \tilde Y_t^\eps=F(X^{\eps}_{\eps t},\cL_{X_{\eps t}^\eps},\tilde Y^{\eps}_{t})\dif t+\dif \tilde W^2_t,\quad \tilde Y^{\eps}_0=\eta.
$$
As $\eps\to0$ and arguing as above, we could freeze the position variable of the  slow process $X_{\eps t}^\eps$ as a parameter ${\bar x}$ and its distribution as a parameter $\mu$. Thus it is natural to seek the frozen equation as
\begin{align}\label{sde05}
\dif Y_t^{{\bar x},\mu}=F({\bar x},\mu,Y_t^{{\bar x},\mu})\dif t+\dif W_t,\quad Y_0^{{\bar x},\mu}=\eta.
\end{align}
We remark that the only difference between (\ref{sde02}) and (\ref{sde05}) is that there exists an additional parameter $\mu$ in (\ref{sde05}). Under exactly the same dissipative condition as before, the process $Y_t^{{\bar x},\mu}$ admits a unique invariant measure $\zeta^{{\bar x},\mu}(\dif y)$ (where $\mu$ is also a parameter).
Then, following the same idea as in case {\it (i)}, the slow component $X_t^\eps$ in system (\ref{sde04}) will converge  as $\eps\to0$ to $\bar X_t$ which satisfies the following averaged equation:
\begin{align}\label{sde003}
\dif \bar X_t=\bar b(\bar X_t,\cL_{\bar X_t})\dif t+\dif W_t^1,\quad\bar X_0=\xi,
\end{align}
where the new   drift  is defined by
$$
\bar b(x,\mu):=\int_{\mR^{d_2}}b(x,\mu,y)\zeta^{x,\mu}(\dif y).
$$
Throughout the whole procedure and in comparison with (\ref{sde01}), the distribution of the slow process in (\ref{sde04}) only appears as a parameter.

\vspace{2mm}
{\it (iii) Fully coupled McKean-Vlasov system.} Now we consider the following McKean-Vlasov system involving the cross interactions of the slow process and the fast motion as well as their distributions (especially the distribution of the fast motion):
\begin{equation} \label{sde06}
\left\{ \begin{aligned}
&\dif X^{\eps}_t =b(X^{\eps}_t,\cL_{X_t^\eps},Y^{\eps}_t,\cL_{Y_t^\eps})\dif t
+\dif W^1_t,\qquad\qquad\qquad X^{\eps}_0=\xi,\\
&\dif Y^{\eps}_t =\frac{1}{\eps}F(X^{\eps}_t,\cL_{X_t^\eps},Y^{\eps}_t,\cL_{Y_t^\eps})\dif t+\frac{1}{\sqrt{\eps}}\dif W_t^2,\qquad\qquad\!\!   Y^{\eps}_0=\eta.
\end{aligned} \right.
\end{equation}
In this case, with the time scaling  $t\mapsto \eps t$, we have that  the re-scaled fast process  $\tilde Y_t^\eps:=Y_{\eps t}^\eps$ satisfies
$$
\dif \tilde Y_t^\eps=F(X^{\eps}_{\eps t},\cL_{X_{\eps t}^\eps},\tilde Y^{\eps}_{t},\cL_{\tilde Y_{t}^\eps})\dif t+\dif \tilde W^2_t,\quad \tilde Y^{\eps}_0=\eta.
$$
Obviously, only the slow process and its distribution should be freezed as $\eps\to0$. Thus one might except as before that  we could choose the frozen equation as
\begin{align}\label{sde07}
\dif Y_t^{\bar x,\mu}=F(\bar x,\mu,Y_t^{\bar x,\mu},\cL_{Y_t^{\bar x,\mu}})\dif t+\dif W_t,\quad Y_0^{\bar x,\mu}=\eta,
\end{align}
where $(\bar x,\mu)$ are parameters. The situation now is quite different with the cases in {\it (i)} and {\it  (ii)}, since the frozen system (\ref{sde07}) is a McKean-Vlasov equation while the systems (\ref{sde02}) and (\ref{sde05}) are classical It\^o SDEs.
As a result, we need to ensure that the non-linear system (\ref{sde07}) admits a unique invariant measure $\zeta^{\bar x,\mu}(\dif y)$.
Then one might except that the slow component $X_t^\eps$ in system (\ref{sde06}) will converge  as $\eps\to0$ to $\bar X_t$ which satisfies the following averaged equation:
\begin{align*}
\dif \bar X_t=\bar b(\bar X_t,\cL_{\bar X_t})\dif t+\dif W_t^1,\quad\bar X_0=\xi,
\end{align*}
where the new   drift  is defined by
\begin{align}\label{bb}
\bar b(x,\mu):=\int_{\mR^{d_2}}b(x,\mu,y,\zeta^{x,\mu})\zeta^{x,\mu}(\dif y).
\end{align}
But the above formally derived averaged equation (\ref{sde03}) for the system (\ref{sde06}) turns out to be not the correct one: {\bf both the choice of the frozen equation (\ref{sde07}) and the definition of the averaged drift (\ref{bb}) are wrong}.

\subsection{New perspective: the fast motion governs the limit of the whole system}

In contrast with the classical theory of averaging
principle for multi-scale systems which focuses on seeking the  limit of the slow process,
let us explain that, the asymptotic behavior of the whole multi-scale system is in fact determined by the limit of the distribution of the fast component, whereas the averaged limit of the slow process follows directly as a byproduct.

\vspace{2mm}
We start with the multi-scale SDE (\ref{sde01}) again. To study  the averaged limit of the slow process $X_t^\eps$ as $\eps\to0$, we may suppose that the limit is denoted by $\bar X_t$, and proceed to seek the equation satisfied by $\bar X_t$. In this way, we naturally have (assume that the coefficient is regular enough) that as $\eps\to0$,
\begin{align*}
&\dif X_t^\eps=b(X_t^\eps,\cdot)\dif t+\dif W_t^1,\\
&\,\,\downarrow\qquad\quad\downarrow\\
&\dif \bar X_t=b(\bar X_t,\cdot)\dif t+\dif W_t^1.
\end{align*}
Thus, the key point  to determine the equation for $\bar X_t$ is to identify the limit of $b(\cdot, Y_t^\eps)$ as $\eps\to0$, i.e.,  the limit  of the distribution of the fast motion. Recall that we have $Y_t^\eps=\tilde Y_{t/\eps}^\eps$, where $\tilde Y_t^\eps$ satisfies the equation (\ref{sde002}), and we have freezed the term $X_{\eps t}^\eps$ in  (\ref{sde002}) as a parameter $\bar x$ to get the frozen equation (\ref{sde02}). Now, intuitively, taking the time as $t/\eps$ in (\ref{sde02}),  replacing the parameter $\bar x$ by $X^\eps_{\eps\cdot t/\eps}=X_t^\eps$ and letting $\eps\to0$, we should have that
\begin{align}\label{yyy}
\mE b(\cdot, Y_t^\eps)=\mE b(\cdot, \tilde Y_{t/\eps}^\eps)\stackrel{\eps\to0}{\longrightarrow }
\mE\left(\int_{\mR^{d_2}}b(\cdot,y)\zeta^{\bar X_t}(\dif y)\right),
\end{align}
where we have used the fact that $X_t^\eps\to\bar X_t$ (the convergence of the parameter term) as pre-assumed.
As a result, we can conclude that the limit equation for the multi-scale SDE (\ref{sde01}) is given by (\ref{sde03}).
Let us point out that the convergence in (\ref{yyy}) can also be seen from the perspective of fluctuation estimate, see e.g. \cite[Lemma 4.2]{RX1} and \cite{RX2}, which says that as $\eps\to0$,
$$
\mE\left(\int_0^tb(\cdot, Y_s^\eps)\dif s\right)\longrightarrow \mE\left(\int_0^t\!\!\int_{\mR^{d_2}}b(\cdot,y)\zeta^{\bar X_s}(\dif y)\dif s\right).
$$

\vspace{1mm}
Using the above agrument, it is quite easy to seek the averaged limit of $X_t^\eps$ in the McKean-Vlasov system (\ref{sde04}).
Namely, if we assume that the limit of $X_t^\eps$ is denoted by $\bar X_t$ (and at the same time, we would have that $\cL_{X_t^\eps}\to\cL_{\bar X_t}$), then we can deduce that as $\eps\to0$,
\begin{align*}
&\dif X_t^\eps=b(X_t^\eps,\cL_{X_t^\eps},\cdot)\dif t+\dif W_t^1,\\
&\,\,\downarrow\qquad\qquad\downarrow\\
&\dif \bar X_t=b(\bar X_t,\cL_{\bar X_t},\cdot)\dif t+\dif W_t^1.
\end{align*}
This is why we said before that the distribution of the slow process in (\ref{sde04}) is only a parameter and does not play an important role, the   equation for $\bar X_t$ will be determined by the limit of $b(\cdot, \cdot, Y_t^\eps)$ as $\eps\to0$. Since the frozen equation (\ref{sde05}) is of the same type as (\ref{sde02}) (i.e., classical It\^o SDE), arguing as in  (\ref{yyy}) (where $\mu$ is only a parameter in the frozen equation (\ref{sde05}), and  as $\eps\to0$, $\mu$ should be replaced by the distribution of $X_{\eps\cdot t/\eps}^\eps=X_t^\eps$ which converges to the distribution of $\bar X_t$ immediately as pre-assumed), we should have that
$$
\mE b(\cdot, \cdot,Y_t^\eps)\longrightarrow \mE\left(\int_{\mR^{d_2}}b(\cdot,\cdot,y)\zeta^{\bar X_t,\cL_{\bar X_t}}(\dif y)\right),
$$
which yields the averaged limit equation (\ref{sde003}).

\vspace{2mm}
Now, for the fully coupled system (\ref{sde06}), assume that the limit of $X_t^\eps$ is denoted by $\bar X_t$, then we have
\begin{align*}
&\dif X_t^\eps=b(X_t^\eps,\cL_{X_t^\eps},\cdot,\cdot)\dif t+\dif W_t^1,\\
&\,\,\downarrow\qquad\qquad\downarrow\\
&\dif \bar X_t=b(\bar X_t,\cL_{\bar X_t},\cdot,\cdot)\dif t+\dif W_t^1.
\end{align*}
As before, the key point is to seek the limits for
\begin{align*}
b(\cdot,\cdot, Y_t^\eps,\cdot)\quad\text{and}\quad b(\cdot,\cdot, \cdot,\cL_{Y_t^\eps}).
\end{align*}
We point out that the former one  is easier  since it involves the distribution of the fast motion only  linearly, and indeed it is a particular case of the later one (allowing non-linear dependence of the distribution).  Unlike (\ref{sde02}) and (\ref{sde05}), the corresponding frozen equation for system (\ref{sde06}) should be  a non-linear one, thus the formal derivation of (\ref{sde07}) and (\ref{bb}) are wrong.  As our result showed  below,  the frozen system of (\ref{sde06}) shall be  given by the following McKean-Vlasov type equation:
\begin{align}\label{ad}
\dif Y_t^{x,\mu}=&F\bigg(x,\mu,Y_t^{x,\mu},\int_{\mR^{d_1}}\cL_{Y_t^{\tilde x,\mu}}\mu(\dif \tilde x)\bigg)\dif t+\dif \tilde W_t,
\end{align}
where $(x,\mu)$ are parameters, and   for test function, we have (see the estimate (\ref{es2}) below)
$$
\psi(\cL_{Y_t^\eps})\stackrel{\eps\to0}{ \longrightarrow }\psi\left(\int_{\mR^{d_1}}\zeta^{x,\cL_{\bar X_t}}\cL_{\bar X_t}(\dif x)\right),
$$
where $\zeta^{x,\mu}$ is the invariant measure of (\ref{ad}). In particular,
$$
\mE \hat \psi(Y_t^\eps)\longrightarrow \int_{\mR^{d_1}}\int_{\mR^{d_2}}\hat\psi(y)\zeta^{x,\cL_{\bar X_t}}(\dif y)\cL_{\bar X_t}(\dif x)=\mE\left(\int_{\mR^{d_2}}\hat\psi(y)\zeta^{\bar X_t,\cL_{\bar X_t}}(\dif y)\right).
$$
Thus, the averaged limit equation of system (\ref{sde06}) should be given by
$$
\dif \bar X_t=\int_{\mR^{d_2}}b\Big(\bar X_t,\cL_{\bar X_t},y,\tilde \mE\big(\zeta^{{\tilde{\bar{X_t}}},\cL_{\bar X_t}}\big)\Big)\zeta^{\bar X_t,\cL_{\bar X_t}}(\dif y)\dif t+\dif W_t^1,
$$
where ${\tilde{\bar{X_t}}}$ is a copy of the limit $\bar X_t$, and the expectation $\tilde \mE$ is taken with respect to ${\tilde{\bar{X_t}}}$.

\subsection{Main results: the asymptotic limit for  the entire system}

We shall show  that as $\eps\to0$, the averaged limit for the fully coupled multi-scale non-linear stochastic  system (\ref{sde1}) is given by following  McKean-Vlasov SDE:
\begin{align}\label{ave}
\dif \bar X_t=\bar b(\bar X_t,\cL_{\bar X_t})\dif t+\bar \sigma(\bar X_t,\cL_{\bar X_t})\dif W^1_t,\qquad  \bar X_0=\xi,
\end{align}
where the averaged coefficients are defined by
\begin{eqnarray}
\begin{split}\label{barb}
&\bar b(x,\mu):=\int_{\mR^{d_2}}b\bigg(x,\mu,y,\int_{\mR^{d_1}}\zeta^{\tilde x,\mu}\mu(\dif \tilde x)\bigg)\zeta^{x,\mu}(\dif y),\\
&\bar \sigma(x,\mu):=\sigma\bigg(x,\mu,\int_{\mR^{d_1}}\zeta^{\tilde x,\mu}\mu(\dif \tilde x)\bigg),
\end{split}
\end{eqnarray}
and $\zeta^{x,\mu}(\dif y)$ is the unique invariant measure of the following McKean-Vlasov type frozen equation:
\begin{align}\label{frozen0}
\dif Y_t^{x,\mu}=&F\bigg(x,\mu,Y_t^{x,\mu},\int_{\mR^{d_1}}\cL_{Y_t^{\tilde x,\mu}}\mu(\dif \tilde x)\bigg)\dif t\no\\
&\qquad+G\bigg(x,\mu,Y_t^{x,\mu},\int_{\mR^{d_1}}\cL_{Y_t^{\tilde x,\mu}}\mu(\dif \tilde x)\bigg)\dif \tilde W_t,\quad Y_0^{x,\mu}=\eta,
\end{align}
where  $(x,\mu)\in\mR^{d_1}\times\sP_2(\mR^{d_1})$ are regarded as parameters, and $\tilde W_t$ is a new standard Brownian motion.
We establish both the strong  convergence   in the averaging principle and the convergence of the distribution (with different rates of convergence) for the slow process $X_t^\eps$.  Moreover, we shall show that the distribution of the fast motion $Y_t^\eps$ will converge to $\cL_{\bar Y_t}$ which is given by
\begin{align}\label{ybar}
\cL_{\bar Y_t}(\dif y):=\int_{\mR^{d_1}}\zeta^{x,\cL_{\bar X_t}}(\dif y)\cL_{\bar X_t}(\dif x)=\mE\Big(\zeta^{\bar X_t,\cL_{\bar X_t}}(\dif y)\Big).
\end{align}
Before stating the main results, let us provide the following comments on the above limits.

\br\label{br}
{\it (i)} The frozen equation (\ref{frozen0}) is indeed a McKean-Vlasov stochastic system. To see this, we define for
every $x\in\mR^{d_1}$, $y\in\mR^{d_2}$, $\mu\in\sP_2(\mR^{d_1})$ and $\nu^x\in\sP_2(\mR^{d_2})$ that
\begin{align*}
&\tilde F\big(x,\mu,y,\nu^x\big):=F\bigg(x,\mu,y,\int_{\mR^{d_1}}\nu^{\tilde x}\mu(\dif \tilde x)\bigg),\\
&\tilde G\big(x,\mu,y,\nu^x\big):=G\bigg(x,\mu,y,\int_{\mR^{d_1}}\nu^{\tilde x}\mu(\dif \tilde x)\bigg).
\end{align*}
Then the  equation (\ref{frozen0}) can be rewritten as
\begin{align}\label{frozen}
\dif Y_t^{x,\mu}=\tilde F\big(x,\mu,Y_t^{x,\mu},\cL_{Y_t^{x,\mu}}\big)\dif t+\tilde G\big(x,\mu,Y_t^{x,\mu},\cL_{Y_t^{ x,\mu}}\big)\dif \tilde W_t,
\end{align}
where $(x,\mu)$ are parameters, and the coefficients depend on the solution $Y_t^{x,\mu}$ as well as its distribution $\cL_{Y_t^{x,\mu}}$. The trick is that the dependence of $\tilde F$ and $\tilde G$ on the parameter $\mu$ come from two parts:
$$
F\bigg(\cdot,\mu,\cdot,\int_{\mR^{d_1}}\nu^{\tilde x}\mu(\dif \tilde x)\bigg)\quad {\text{and}}\quad G\bigg(\cdot,\mu,\cdot,\int_{\mR^{d_1}}\nu^{\tilde x}\mu(\dif \tilde x)\bigg).
$$

{\it (ii)} In contrast with (\ref{bb}) and (\ref{sde07}),   there exists an integral with respect to the measure $\mu$ in the coefficients. This is reasonable in the sense that, if we consider a particular case of the non-linear system (\ref{sde1}) with $\sigma=\mI_d$ and the drift $b$ depends only on the distribution of the fast motion, i.e., for $b: \sP_2(\mR^{d_2})\to\mR^{d_1}$,
$$
\dif X_t^\eps=b(\cL_{Y_t^\eps})\dif t+\dif W_t^1,
$$
then obviously, the limit of the term $b(\cL_{Y_t^\eps})$ should be a deterministic one. The limit equation (\ref{ave}) (or the limit in (\ref{ybar})) implies that this term will converge to
$$
b\left(\int_{\mR^{d_1}}\zeta^{x,\cL_{\bar X_t}}(\cdot)\cL_{\bar X_t}(\dif x)\right)\in\mR^{d_1},
$$
whereas (\ref{bb}) becomes
$$
b\Big(\zeta^{\bar X_t,\cL_{\bar X_t}}(\cdot)\Big),
$$
which is a stochastic process due to the existence of $\bar X_t$.
\er

\vspace{2mm}
To study the asymptotic limit  of the non-linear stochastic system (\ref{sde1}) with irregular coefficients, we  assume the following basic non-degeneracy conditions on the diffusion coefficients:

\vspace{2mm}
\begin{description}
\item[{\bf ($\bf H_{1}$)}] the coefficients $a=\sigma\sigma^*$ and $\cG=GG^*$ are   non-degenerate in the sense that there exist constants $k, \varrho>0$ such that for any $(x,y)\in\mR^{d_1}\times\mR^{d_2}$, $\mu\in\sP_2(\mR^{d_1})$ and $\nu\in\sP_2(\mR^{d_2})$,
\begin{align*}
|\sigma^*(x,\mu,\nu)z_1|^2\geq\varrho|z_1|^2,\qquad \forall z_1\in\mR^{d_1},
\end{align*}
and
\begin{align*}
|G^*(x,\mu,y,\nu)z_2|^2\geq\varrho(1+|y|)^{-k}|z_2|^2,\qquad \forall z_2\in\mR^{d_2}.
\end{align*}
\end{description}

\vspace{1mm}
Given a function $\cV: \mR^{d_2}\to\mR_+$, recall that the weighted total variation distance between two probability measures $\nu_1$ and $\nu_2$ on $\mR^{d_2}$ is defined by
\begin{align}\label{tv}
\rho_\cV(\nu_1,\nu_2)&:=\int_{\mR^{d_2}}\big(1+\cV(y)\big)|\nu_1-\nu_2|(\dif y)\no\\
&=\sup_{\|f\|_{1+\cV}\leq 1}\int_{\mR^{d_2}}f(y)(\nu_1-\nu_2)(\dif y),
\end{align}
where the weighted supremum norm is given by
$$
\|f\|_{1+\cV}:=\sup_{y\in\mR^{d_2}}\frac{|f(y)|}{1+\cV(y)}.
$$
We  make the following dissipative assumptions on the coefficients of the fast motion to ensure the existence of a unique invariant measure for the frozen system (\ref{frozen0}):

\vspace{2mm}
\begin{description}
\item[{\bf ($\bf H_{2}$)}]
For any $q\geq2$, there exist constants $C_1>C_2\geq 0,$ $C_3\geq0$ such that for any $(x,y)\in\mR^{d_1}\times\mR^{d_2}$, $\mu\in\sP_2(\mR^{d_1})$ and $\nu\in\sP_2(\mR^{d_2})$,
\begin{align}\label{dissf}
2\<F(x,\mu,y,\nu),y\>+(q-1)\|G(x,\mu,y,\nu)\|^2\leq -C_1|y|^{2}+C_2\|\nu\|_{2}^{2}+C_3,
\end{align}
and there exists  $\kappa>0$ small  and $p\geq 1$ such that for any $\nu_1, \nu_2\in\sP_2(\mR^{d_2})$,
\begin{align}\label{dissf2}
&|F(x,\mu,y,\nu_1)-F(x,\mu,y,\nu_2)|\no\\
&\quad+\|G(x,\mu,y,\nu_1)-G(x,\mu,y,\nu_2)\|\leq \kappa\,\rho_\cV(\nu_1,\nu_2),
\end{align}
where $\cV(y)=1+|y|^p$.
\end{description}

\vspace{1mm}
We give  the following comments on the above assumption.

\br
The dissipative condition (\ref{dissf}) is mainly used to establish  the existence of invariant measures for the frozen system (\ref{frozen0}), which is much weaker than the one-side Lipschitz assumptions (see e.g. \cite{LWX,W1}). Since (\ref{frozen0}) is a McKean-Vlasov system, it is well known (see e.g. \cite{D,DYY,T,ZS}) that the existence of several invariant measures may occur for non-convex confining potential. This phenomenon is referred to as phase transition. The assumption (\ref{dissf2}) ensures the uniqueness of the invariant measure, and the smallness of  $\kappa$  in (\ref{dissf2}) is essential (see similar assumption in \cite[Theorem 3.1]{W2}) in view of the work of D. A. Dawson: \cite{D} established  the phase transition for the McKean-Vlasov equation with a particular
double-well confinement, which shows that there exists a $\kappa_c$ such that if $\kappa\geq \kappa_c$, then the corresponding system admits three invariant measures. In this case, the characterization of the basin of attractions of these different invariant measures is more difficult and there exists very few results in the literature, see \cite{T} for partial result. Since we can identify the limiting distribution of the fast motion, the condition (\ref{dissf2}) is naturally needed  (otherwise, this will imply the characterization of the basin of attractions  even if in the  very particular case where the fast motion in the non-linear system (\ref{sde1}) does not depend on the slow process and its distribution). Note that we do not need the coefficients of the slow equation satisfy such condition. If $F$ and $G$  admit a linear functional derivative  with
$$
\Big|\frac{\delta F}{\delta\nu}(\cdot,\cdot,\cdot,\nu)(\tilde y)\Big|\leq \kappa (1+|\tilde y|^p)\quad \mathrm{and} \quad \Big|\frac{\delta G}{\delta\nu}(\cdot,\cdot,\cdot,\nu)(\tilde y)\Big|\leq \kappa (1+|\tilde y|^p),
$$
then (\ref{dissf2}) holds. We also remark that the coefficients may not be Lipschitz continuous with respect to the $\cW_2$-Wasserstein distance, see Remark \ref{r1} below.
\er

\vspace{2mm}
Fix $T>0$. Let $(X_t^\eps, Y_t^\eps)$ and $\bar X_t$  satisfy the McKean-Vlasov equations (\ref{sde1}) and (\ref{ave}), respectively. The following is the main result of this paper. For brevity, the spaces of functions mentioned below are introduced in the Notation part at the end of this section.

\bt\label{main1}
Let {\bf ($\bf H_{1}$)}  and {\bf ($\bf H_{2}$)} hold. Assume that $b,F,G\in C_p^{\a,(2,\a),\b,(2,\b)}$ and $\sigma\in C_p^{\a,(2,\a),(2,\b)}$ with $\a,\b>0$. Then  for any  $t\in[0,T]$, we have

\vspace{1mm}
\noindent (i) (convergence of the distribution of the slow process) for every $\varphi\in C_b^{(2,\a)}(\sP_2(\mR^{d_1}))$,
\begin{align}\label{es1}
|\varphi(\cL_{X_t^\eps})-\varphi(\cL_{\bar X_t})|\leq C_T\,\eps^{\frac{\a}{2}\wedge1},
\end{align}
where  $C_T>0$ is a constant  independent of $\eps$;

\vspace{1mm}
\noindent (ii) (limit for the distribution of the fast motion) for every $\psi\in C_p^{(2,\beta)}(\sP_2(\mR^{d_2}))$,
\begin{align}\label{es2}
|\psi(\cL_{Y_t^\eps})-\psi(\tilde\zeta^{\cL_{\bar X_t}})|\leq C_T\,\eps^{\frac{\a}{2}\wedge1}+C_0\,\e^{-\frac{\gamma_0t}{\eps}},
\end{align}
where $C_0>0$ and $\gamma_0>0$ are constants independent of $T$ and $\eps$, $\zeta^{x,\mu}$ is the unique invariant measure for the frozen equation (\ref{frozen0}),  and
\begin{align}\label{zmu}
\tilde \zeta^{\mu}(\cdot):=\int_{\mR^{d_1}}\zeta^{x,\mu}(\cdot)\mu(\dif x);
\end{align}

\vspace{1mm}
\noindent (iii) (strong convergence of the slow process) assume further that $\sigma\in C_b^{1,(2,\a),(2,\b)}$, then
\begin{align}\label{es3}
\mE|X_t^\eps-\bar X_t|^2\leq C_T\,\eps^{\a\wedge1}.
\end{align}

\vspace{1mm}
In particular, we have for every $\hat\varphi\in C_b^\a(\mR^{d_1})$,
\begin{align*}
\sup_{t\in[0,T]}|\mE\hat\varphi(X_t^\eps)-\mE\hat\varphi(\bar X_t)|\leq C_T\,\eps^{\frac{\a}{2}\wedge1},
\end{align*}
and for every $\hat\psi\in C_p^\beta(\mR^{d_2})$,
\begin{align*}
\bigg|\mE \hat\psi(Y_t^\eps)-\mE\bigg[\int_{\mR^{d_2}} \hat\psi(y)\zeta^{\bar X_t,\cL_{\bar X_t}}(\dif y)\bigg]\bigg|\leq C_T\,\eps^{\frac{\a}{2}\wedge1}+C_0\,\e^{-\frac{\gamma_0t}{\eps}}.
\end{align*}
\et

\vspace{2mm}
We provide the following remark for the above result.

\br\label{r1}
(i) Note that all the coefficients are  not differentiable with respect to the measure variables in the sense of Lions. In fact, they are only H\"{o}lder continuous with respect to the measures in the  Wasserstein distance. Let us explain this for the coefficient $b$ with respect to the $\mu$-variable when $0<\a< 1$. By (\ref{lin}) below, we have for  every $\mu_1,\mu_2\in\sP_2(\mR^{d_1})$,
\begin{align*}
&|b(x,\mu_1,y,\nu)-b(x,\mu_2,y,\nu)|\\
&=\left|\int_0^1\int_{\mR^{d_1}}\frac{\d b}{\d\mu}(x,\theta\mu_1+(1-\theta)\mu_2,y,\nu)(z)(\mu_1-\mu_2)(\dif z)\dif\theta\right|\\
&\leq\int_0^1\int_{\mR^{d_1}\times\mR^{d_1}}\Big|\frac{\d b}{\d\mu}(x,\theta\mu_1+(1-\theta)\mu_2,y,\nu)(z_1)\\
&\quad-\frac{\d b}{\d\mu}(x,\theta\mu_1+(1-\theta)\mu_2,y,\nu)(z_2)\Big|\pi(\dif z_1,\dif z_2)\dif\theta\\
&\leq C_0(1+|y|^p+\|\nu\|_2^p)\cdot\int_{\mR^{d_1}}|z_1-z_2|^\a\pi(\dif z_1,\dif z_2),
\end{align*}
where the measure $\pi$ is an arbitrary coupling of $\mu_1$ and $\mu_2$. Due the arbitrariness of $\pi$, we arrive at
\begin{align*}
&|b(x,\mu_1,y,\nu)-b(x,\mu_2,y,\nu)|\\
&\leq C_0(1+|y|^p+\|\nu\|_2^p)\cdot\Big(\inf_{\pi\in\cP(\mu_1,\mu_2)}\int_{\mR^{d_1}}|z_1-z_2|^\a\pi(\dif z_1,\dif z_2)\Big)\\
&\leq C_0(1+|y|^p+\|\nu\|_2^p)\cdot \cW_2(\mu_1,\mu_2)^\a.
\end{align*}
Meanwhile, in view of the estimate (\ref{lip}) below,  the coefficients are Lipschitz in the weighted total variation distance. This seems sharp when one considers using the finite dimensional noise to regularise a function defined on the infinite dimensional space $\sP_2(\mR^d)$. Under our assumptions, the  well-posedness of the system (\ref{sde1}) can be obtained by \cite[Theorem 3.4, Corollary 3.5]{CF2} or \cite{HSS}. The existence of  invariant measures for the McKean-Vlasov SDE (\ref{frozen}) under condition (\ref{dissf}) can be found in  \cite[Theorem 2.2]{ZS}, see also Lemma \ref{inv} below for an alternate proof using our arguments. The uniqueness of invariant measure for system (\ref{frozen}) under (\ref{dissf2}) can be proved similarly as in \cite[Theorem 3.1]{W2}, since we will not need this property in our proof and for the sake of simplicity, we do not
deal with this problem in the present article and postpone it to another work. In addition, we shall show that the averaged coefficient $\bar b, \bar \sigma$ defined in (\ref{barb}) are also H\"{o}lder continuous, i.e., $\bar b\in C_b^{\a,(2,\a)}$ and $\bar\sigma\in C_b^{\a,(2,\a)}$ (and $\bar\sigma\in C_b^{1,(2,\a)}$ in the case of {\it(iii)}) , see Lemma \ref{lem3} below. Thus, the weak and strong well-posedness of   the averaged equation (\ref{ave}) follows by \cite{C,CF2}.

\vspace{2mm}
(ii) The estimate (\ref{es1}) for the convergence of the distribution of the slow process is even more general than the classical weak convergence in the averaging principle of the  multi-scale It\^o SDEs (see also Remark \ref{resde} below for more explanations) since non-linear test functions of the distribution are allowed. Such kind of estimate was also obtained in \cite[Theorem 3.1]{BS2} when  the system (\ref{sde1}) does not involve the distribution of the fast motion, and  our regularity assumptions on the coefficients as well as the test functions are much weaker than the previous results. In addition, both the strong and weak convergence rates in the averaging principle obtained in the estimates (\ref{es3}) and (\ref{es1}) coincide with the cases in the classical It\^o SDEs, see e.g. \cite{RX1,RX2}.

\vspace{2mm}
(iii) The estimate (\ref{es2}) is new and seems to be the first result established for the convergence of the distribution of the fast motion even for classical multi-scale It\^o SDEs.
Estimate (\ref{es2}) implies that for every $t>0$, the distribution $\cL_{Y_t^\eps}$ of the fast motion will converge  to $\mE\zeta^{\bar X_t,\cL_{\bar X_t}}$ as $\eps\to0$, and the rate of convergence is given by $\eps^{\frac{\a}{2}\wedge1}+\e^{-\frac{\gamma_0 t}{\eps}}$, which is independent of the regularity index $\b$ (the regularity of the coefficients with respect to the fast motion and its distribution).
Note that the constant $C_0$ in (\ref{es2}) is independent of the time variable.                                         The exponential decay term in the rate is natural since even if the multi-scale system is not fully coupled, i.e., the fast motion in the McKean-Vlasov system (\ref{sde1}) does not depend on the slow process and its distribution, we would  have that for every $t>0$,
\begin{align*}
|\psi(\cL_{Y_t^\eps})-\psi(\zeta)|\leq C_0\,\e^{-\frac{\gamma_0t}{\eps}},
\end{align*}
where $\zeta$ is the unique invariant measure for the fast McKean-Vlasov SDE. In view of the assumption on the test function $\psi$, this in particular implies the exponential ergodic of the fast motion in the weighted total variation distance, which is of independent interest.
\er

The rest of this paper is structured as follows. In Section 2, we briefly explain the idea used to study the asymptotic behavior of the McKean-Vlasov system (\ref{sde1}). Section 3 is devoted to study the optimal regularities for the solutions of two kinds of Kolmogorov equation in Wasserstein space. In Section 4, we state some results about the Poisson equation with parameters and introduce an mollifying approximation on Wasserstein space. In Section 5, we study the asymptotic behavior of the non-autonomous multi-scale SDEs (\ref{sde2}) by using the results obtained in Section 3.  Finally, we give the proof of Theorem \ref{main1} in Section 6.

\vspace{2mm}
\noindent{\bf Notations.} Let us first briefly recall two kinds of differential calculus on the space of measures $\sP_2(\mR^{d})$, for more  complete and detailed exposition, we refer the readers to \cite[Section 2]{CF2} or \cite{CD,Lions}.
The  first is the linear functional derivative, which is a standard notion of differentiability for functions of measures relying on the convexity of $\sP_2(\mR^{d})$. Given a real-valued function $f$  on  $\sP_2(\mR^{d})$, we say that $f$ admits a linear functional derivative if there exists a real-valued and continuous function $[\d f/\d\mu](\mu)(x)$ defined on $\sP_2(\mR^{d})\times\mR^{d}$   such that for all $\mu_1,\mu_2\in\sP_2(\mR^{d})$,
\begin{align}\label{lin0}
\lim_{\theta\to0}\frac{f(\mu_1+\theta(\mu_2-\mu_1))-f(\mu_1)}{\theta}=\int_{\mR^d}\frac{\d f}{\d\mu}(\mu_1)(x)(\mu_2-\mu_1)(\dif x).
\end{align}
The map $x\mapsto[\d f/\d\mu](\mu)(x)$ being defined up to an additive constant, we will follow the usual normalization convention $\int_{\mR^d}[\d f/\d\mu](\mu)(x)\mu(\dif x)=0$. Note that by definition, we have  for every $\mu_1,\mu_2\in\sP_2(\mR^{d})$,
\begin{align}\label{lin}
f(\mu_1)-f(\mu_2)=\int_0^1\int_{\mR^d}\frac{\d f}{\d\mu}(\theta\mu_1+(1-\theta)\mu_2)(x)(\mu_1-\mu_2)(\dif x)\dif\theta.
\end{align}
In particular,  if  $f$ admits a bounded linear functional derivative, then it is Lipschitz continuous with respect to the total variation distance since by (\ref{lin}) we have
\begin{align}\label{lip}
|f(\mu_1)-f(\mu_2)|\leq \sup_{\mu\in\sP_2(\mR^d)}\Big\|\frac{\delta f}{\delta\mu}(\mu)(\cdot)\Big\|_\infty\cdot\|\mu_1-\mu_2\|_{{\text{TV}}}
\end{align}
The second notion of differentiation with respect to the measure variable we shall used was introduced by Lions. We say that $f$ is $L$-differentiable if its lifting defined by
$\mF: L^2(\Omega)\ni X\rightarrow f(\cL_X)\in\mR$,
is Fr\'echet differentiable. Moreover, there exists a function $\p_\mu f(\mu)(\cdot): \mR^d\rightarrow\mR^d$ such that
$$
D\mF(X)=\p_\mu f(\cL_X)(X).
$$
The function $\p_\mu f(\mu)(x)$ is then called the Lions derivative ($L$-derivative for short) of $f$ at $\mu$. The advantage of the $L$-derivative is that it permits to use the tools of differential calculus on Banach spaces. If $f$ is continuously $L$-differentiable and if the Fr\'echet derivative of its lift  is bounded in $L^2(\Omega)$, then it is Lipschitz continuous with respect to the $\cW_2$-Wasserstein distance since by Cauchy-Schwarz's inequality we have for all $\mu_1,\mu_2\in\sP_2(\mR^{d})$,
\begin{align*}
|f(\mu_1)-f(\mu_2)|\leq \|D\mF\|_{L^2(\Omega)}\cdot\cW_2(\mu_1,\mu_2).
\end{align*}
As underlined in \cite{CD}, the following relation holds between the $L$-derivative and the linear functional derivative:
$$
\p_\mu f(\mu)(x)=\p_x\frac{\d f}{\d\mu}(\mu)(x).
$$
The higher order derivatives of $f$ at $\mu$ can be defined similarly.

\vspace{2mm}

To end this section and for simplicity, we provide the following notations used in this paper.	
Given a function space, the subscript $b$ will stand for boundness, while the subscript $p$ stands for polynomial growth in $y$ and $\nu$. More precisely, for a function $f(t,x,\mu,y,\nu)\in L^\infty_p:=L^{\infty}_p(\mR_+\times\mR^{d_1}\times\sP_2(\mR^{d_1})\times\mR^{d_2}\times\sP_2(\mR^{d_2}))$, we mean there exist constants $C, p>0$ such that for any $t>0, x\in\mR^{d_1}, \mu\in\sP_2(\mR^{d_1}), y\in\mR^{d_2}$ and $\nu\in\sP_2(\mR^{d_2})$,
$$
|f(t,x,\mu,y,\nu)|\leq C(1+|y|^p+\|\nu\|_2^p),
$$
where $\|\nu\|_2$ represents the $2$-order moments of $\nu$, and we let
$$
\|f\|_{L^\infty_p}:=\sup_{t,x,\mu,y,\nu}\frac{|f(t,x,\mu,y,\nu)|}{1+|y|^p+\|\nu\|_2^p}.
$$
We shall use the following spaces of functions:
for $k_1,k_2\in\mN_+$ and $0<\a,\b\leq2$,

\begin{itemize}

\item The space $C^{(k,\a)}_b:=C^{(k,\a)}_b(\sP_2(\mR^d))$. A function $f(\mu)$ is in $C^{(k,\a)}_b$ if $f$ admits $k$-order bounded linear functional derivatives $[\delta^k f/\delta \mu^k](\mu)(x_1,\cdots,x_k)$ such that for every $k\geq 1$, the map $x_k\mapsto[\delta^k f/\delta \mu^k](\mu)(x_1,\cdots,x_k)$ is $\a$-H\"{o}lder continuous uniformly with respect to other variables.

\item The space $C_p^{(k,\b)}:=C_p^{(k,\b)}(\sP_2(\mR^{d}))$. A function $f(\nu)$ is in $C_p^{(k,\b)}$  if $f$ admits $k$-order linear functional derivative $[\delta^{k} f/\delta \nu^{k}](\nu)(y_1,\cdots,y_{k})$ such that the derivative is polynomially grow in $(y_1,\cdots,y_{k})$ uniformly with respect to other variables, and for every $k\geq 1$, the map  $y_{k}\mapsto[\delta^{k} f/\delta \nu^{k}](\nu)(y_1,\cdots,y_{k})$ is local $\b$-H\"{o}lder continuous with polynomial growth in $(y_1,\cdots,y_{k})$.

\item The space $C_p^{(k_1,\a),(k_2,\b)}:=C_p^{(k_1,\a),(k_2,\b)}(\sP_2(\mR^{d_1})\times\sP_2(\mR^{d_2}))$. A function $f(\mu,\nu)$ is in $C_p^{(k_1,\a),(k_2,\b)}$  if $f$ admits $k_1$-order linear derivatives $[\delta^{k_1} f/\delta \mu^{k_1}](\mu,\nu)(x_1,\cdots,x_{k_1})$ and $k_2$-order linear functional derivatives $[\delta^{k_2} f/\delta \nu^{k_2}](\mu,\nu)(y_1,\cdots,y_{k_2})$ such that these derivatives are polynomially grow in $(y_1,\cdots,y_{k})$ uniformly with respect to other variables, and the map  $x_{k_1}\mapsto[\delta^{k_1} f/\delta \mu^{k_1}](\mu,\nu)(x_1,\cdots,x_{k_1})$ is $\a$-H\"{o}lder continuous, and the map $y_{k_2}\mapsto[\delta^{k_2} f/\delta \nu^{k_2}](\mu,\nu)(y_1,\cdots,y_{k_2})$ is local $\b$-H\"{o}lder continuous with polynomial growth in $(y_1,\cdots,y_{k})$. Similarly, we also can define the space $C_b^{(k_1,\a),(k_2,\b)}(\sP_2(\mR^{d_1})\times\sP_2(\mR^{d_2}))$.

  \item The space $C_b^{\a/2,\a,(k_1,\a),\b,(k_2,\b)}:=C_b^{\a/2,\a,(k_1,\a),\b,(k_2,\b)}(\mR^+\times\mR^{d_1}
      \times\sP_2(\mR^{d_1})\times\mR^{d_2}\times\sP_2(\mR^{d_2}))$. A function $f(t,x,\mu,y,\nu)$ is in $C_b^{\a/2,\a,(k_1,\a),\b,(k_2,\b)}$  if for every $\mu\in\sP_2(\mR^{d_1})$ and $\nu\in\sP_2(\mR^{d_2})$,
      $f(\cdot,\cdot,\mu,\cdot,\nu)\in C_b^{\a/2,\a,\b}(\mR^+\times\mR^{d_1}\times\mR^{d_2})$(i.e., the usual H\"{o}lder space), and for every $(t,x,y)\in\mR^+\times\mR^{d_1}\times\mR^{d_2}$, $f(t,x,\cdot,y,\cdot)\in C_b^{(k_1,\a),(k_2,\b)}(\sP_2(\mR^{d_1})\times\sP_2(\mR^{d_2})) $ .

  \item The space $C_p^{\a/2,\a,(k_1,\a),\b,(k_2,\b)}:=C_p^{\a/2,\a,(k_1,\a),\b,(k_2,\b)}(\mR^+\times\mR^{d_1}
      \times\sP_2(\mR^{d_1})\times\mR^{d_2}\times\sP_2(\mR^{d_2}))$. A function $f(t,x,\mu,y,\nu)$ is in $C_p^{\a/2,\a,(k_1,\a),\b,(k_2,\b)}$  if $f\in L_p^\infty$ and for every $\mu\in\sP_2(\mR^{d_1})$ and $\nu\in\sP_2(\mR^{d_2})$,
      $f(\cdot,\cdot,\mu,\cdot,\nu)\in C_p^{\a/2,\a,\b}(\mR^+\times\mR^{d_1}\times\mR^{d_2})$ (where $C_p^{\a/2,\a,\b}$ consists of all functions that are $\a/2$-H\"{o}lder continuous in $t$, $\a$-H\"{o}lder continuous in $x$ and $\b$-H\"{o}lder continuous with polynomial growth in $y$), and for every $(t,x,y)\in\mR^+\times\mR^{d_1}\times\mR^{d_2}$, $f(t,x,\cdot,y,\cdot)\in C_p^{(k_1,\a),(k_2,\b)}(\sP_2(\mR^{d_1})\times\sP_2(\mR^{d_2}))$.
\end{itemize}

\section{Idea of method: the non-autonomous approximation argument}

The aim of this section is to explain our idea of studying the asymptotic behavior of the entire non-linear stochastic system (\ref{sde1}) by using the non-autonomous approximations. At the same time, we shall point out the main difficulties, especially those caused by the low regularity of the coefficients.

\vspace{2mm}
For clarity, we  divide the explanation of the idea for  proving the strong convergence in the averaging principle and the convergence of the distributions of $X_t^\eps$ and $Y_t^\eps$ into the following four steps:

\vspace{2mm}
\noindent{\bf Step 1.} We first consider the following linear but non-autonomous multi-scale stochastic system in $\mR^{d_1}\times\mR^{d_2}$:
\begin{equation} \label{sde2}
\left\{ \begin{aligned}
&\dif \hat X^{\eps}_t =b_\eps(t,\hat X^{\eps}_t,\hat Y^{\eps}_t)\dif t
+\sigma_\eps(t,\hat X^{\eps}_t)\dif W^1_t,\qquad\qquad\qquad\, \hat X^{\eps}_0=\xi,\\
&\dif \hat Y^{\eps}_t =\frac{1}{\eps}F_\eps(t,\hat X^{\eps}_t,\hat Y^{\eps}_t)\dif t+\frac{1}{\sqrt{\eps}}G_\eps(t,\hat X^{\eps}_t,\hat Y^{\eps}_t)\dif W_t^2,\quad\quad   \hat Y^{\eps}_0=\eta.
\end{aligned} \right.
\end{equation}
Note that the coefficients in the system (\ref{sde2}) depend not only on the time variable $t$, but also on the scale parameter $\eps$. Roughly speaking, we shall show that: if the coefficients
$b_\eps(t,x,y)$, $\sigma_\eps(t,x)$, $F_\eps(t,x,y)$ and $G_\eps(t,x,y)$ converge (in certain sense) to some $\hat b(t,x,y)$, $\hat \sigma(t,x)$, $\hat F(t,x,y)$ and $\hat G(t,x,y)$ as $\eps\to 0$, respectively,
then we can identify the averaged limit $\bar{\hat X}_t$ of the system (\ref{sde2}) and prove the strong convergence as well as the convergence of the distribution of the slow process $\hat X^{\eps}_t$. More importantly, we  give explicit characterization for the limit $\cL_{\bar{\hat Y}_t}$ of the distribution of the fast motion $\hat Y_t^\eps$. Meanwhile,  explicit rates of convergence  depending on the convergence of $b_\eps$, $\sigma_\eps$, $F_\eps$ and $G_\eps$ to $\hat b$, $\hat \sigma$, $\hat F$ and $\hat G$ are also obtained.
These  results are presented in Theorem \ref{non-aut}, which  will play an important role in the study of the asymptotic behavior of the non-linear stochastic system (\ref{sde1}).

\vspace{2mm}
Here, we give the following important comments.

\begin{itemize}
  \item The convergence of the distributions  of the slow process and the fast motion is more general than the classical weak convergence of the multi-scale SDEs. More precisely, we obtain the convergence of
\begin{align}\label{77}
 \varphi\big(\cL_{\hat X_t^\eps}\big)\longrightarrow\varphi\big(\cL_{\bar{\hat X}_t}\big)\quad\text{and}\quad \psi\big(\cL_{\hat Y_t^\eps}\big)\longrightarrow\psi\big(\cL_{\bar{\hat Y}_t}\big)
 \end{align}
  as $\eps\rightarrow0$, where $\varphi, \psi$ are test functions defined on space of measures (allowing non-linear test functions). This is essential to study the non-linear system (\ref{sde1}) as the coefficients depend non-linearly on the distributions of the solutions.
Especially, to prove the convergence of the distribution of the fast motion in the above sense, we need to study the long time decay of the solution of the forward Kolmogorov equation on the product Wasserstein space $\sP_2(\mR^{d_1}\times\mR^{d_2})$, see subsection 3.2.

  \item The regularity assumptions on the coefficients and the test functions in (\ref{77}) are low. Especially, $\varphi$ and $\psi$  are even not Lipschitz continuous w.r.t. the $W_2$-distance (and thus not Lions differentiable). This allows us to make low regularity assumptions on the coefficients of the original non-linear system (\ref{sde1}). To overcome this difficulty, we need to derive the optimal regularities of the solutions of the (backward and forward) Kolmogorov equations on the Wasserstein space. The simplicity lies in that,  we only need to handle the  Kolmogorov equations associated with the linear system (\ref{sde2}) but not the mean-field type PDEs associated with the original non-linear system (\ref{sde1}). Besides, we need to seek an mollifying argument on Wasserstein space with explicit approximating rate and bounds on the Lions derivatives of the approximation functions, which is important to derive the optimal rates of convergence for the system (\ref{sde1}).

      \item The explicit dependence on the convergence of $b_\eps$, $\sigma_\eps$, $F_\eps$ and $G_\eps$ to $\hat b$, $\hat \sigma$, $\hat F$ and $\hat G$ in the convergence rates of (\ref{77}) will play an important role in studying the asymptotic behavior of the non-linear system (\ref{sde1}).

\end{itemize}

\vspace{2mm}
\noindent{\bf Step 2.}
For every $t\geq 0$ and $\eps>0$, let us denote
$$
X_t^{0,\eps}\equiv \xi\quad\text{and}\quad Y_t^{0,\eps}\equiv \eta.
$$
For $n\geq 1$,
we introduce the following approximations of the non-linear system (\ref{sde1}):
\begin{equation} \label{sden}
\left\{ \begin{aligned}
&\dif X^{n,\eps}_t =b\big(X^{n,\eps}_t,\cL_{X_t^{n-1,\eps}},Y^{n,\eps}_t,\cL_{Y_t^{n-1,\eps}}\big)\dif t\\
&\qquad\qquad\quad+\sigma\big(X^{n,\eps}_t,\cL_{X_t^{n-1,\eps}},\cL_{Y_t^{n-1,\eps}}\big)\dif W^1_t,\qquad\qquad\qquad X^{n,\eps}_0=\xi,\\
&\dif Y^{n,\eps}_t =\frac{1}{\eps}F\big(X^{n,\eps}_t,\cL_{X_t^{n-1,\eps}},Y^{n,\eps}_t,\cL_{Y_t^{n-1,\eps}}\big)\dif t\\
&\qquad\qquad\quad+\frac{1}{\sqrt{\eps}}G\big(X^{n,\eps}_t,\cL_{X_t^{n-1,\eps}},Y^{n,\eps}_t,\cL_{Y_t^{n-1,\eps}}\big)\dif W_t^2,\qquad
 Y^{n,\eps}_0=\eta.
\end{aligned} \right.
\end{equation}
Note that for every fixed $n\geq 1$, the above system is a linear one (classical It\^o SDE with time-dependent coefficients)
since the distributions appearing in the coefficients are not $\cL_{X_t^{n,\eps}}$ and $\cL_{Y_t^{n,\eps}}$ but rather $\cL_{X_t^{n-1,\eps}}$ and $\cL_{Y_t^{n-1,\eps}}$, i.e., the distributions of
the solutions of the approximations at the previous step. Each of the approximation system (\ref{sden}) can be viewed as a linear  non-autonomous system of the form (\ref{sde2}).
In fact, for $n\geq 1$, define
$$
b_{n,\eps}(t,x,y):=b\big(x,\cL_{X_t^{n-1,\eps}},y,\cL_{Y_t^{n-1,\eps}}\big), \quad \sigma_{n,\eps}(t,x):=\sigma\big(x,\cL_{X_t^{n-1,\eps}},\cL_{Y_t^{n-1,\eps}}\big),
$$
and
$$
F_{n,\eps}(t,x,y):=F\big(x,\cL_{X_t^{n-1,\eps}},y,\cL_{Y_t^{n-1,\eps}}\big),\quad G_{n,\eps}(t,x,y):=G\big(x,\cL_{X_t^{n-1,\eps}},y,\cL_{Y_t^{n-1,\eps}}\big).
$$
Then, the system (\ref{sden}) can be rewritten as
\begin{equation} \label{sde4}
\left\{ \begin{aligned}
&\dif X^{n,\eps}_t =b_{n,\eps}(t,X^{n,\eps}_t,Y^{n,\eps}_t)\dif t
+\sigma_{n,\eps}(t,X^{n,\eps}_t)\dif W^1_t,\qquad\qquad\qquad X^{n,\eps}_0=\xi,\\
&\dif Y^{n,\eps}_t =\frac{1}{\eps}F_{n,\eps}(t,X^{n,\eps}_t,Y^{n,\eps}_t)\dif t+\frac{1}{\sqrt{\eps}}G_{n,\eps}(t,X^{n,\eps}_t,Y^{n,\eps}_t)\dif W_t^2,\quad\,   Y^{n,\eps}_0=\eta.
\end{aligned} \right.
\end{equation}
For every fixed $n$, the convergence of the coefficients $b_{n,\eps}(t,x,y)$, $\sigma_{n,\eps}(t,x)$, $F_{n,\eps}(t,x,y)$ and $G_{n,\eps}(t,x,y)$ as $\eps\to0$ follows from the  convergence of the distributions of $X_t^{n-1,\eps}$ and $Y_t^{n-1,\eps}$ obtained in the previous step (by taking properly test functions relying on the coefficients in the estimate (\ref{77})). Thus, as a direct consequence of Step 1 and the induction argument, we obtain immediately the asymptotic limit $\cL_{\bar X_t^n}$ and $\cL_{\bar Y_t^n}$ for the distributions of the slow process $X_t^{n,\eps}$ and the fast motion $Y_t^{n,\eps}$ in the system (\ref{sde4}) with explicit rates of convergence for every $n\geq 1$, respectively, see Theorem \ref{weakn}. Moreover,
we show that the  convergence rates are uniform with respect to $n$ (this is essentially due to the coefficients in (\ref{sde4}) satisfy conditions uniformly with respect to $n$), which will play a crucial role below for  taking the limit as $n\rightarrow\infty$.

\vspace{2mm}
\noindent{\bf Step 3.} We seek the limit as $n\to\infty$. Suppose that the limits of the distributions of  $X_t^\eps$ and $Y_t^\eps$ of the non-linear stochastic system (\ref{sde1}) are denoted by $\cL_{\bar X_t}$ and $\cL_{\bar Y_t}$, respectively. To identify the equation satisfied by $\bar X_t$ and characterize $\cL_{\bar Y_t}$, we deduce that for test functions,
\begin{align}\label{add}
&\big|\varphi(\cL_{X_t^\eps})-\varphi(\cL_{\bar X_t})\big|+\big|\psi(\cL_{Y_t^\eps})-\psi(\cL_{\bar Y_t})\big|\no\\
&\leq \Big[\big|\varphi(\cL_{X_t^\eps})-\varphi(\cL_{X_t^{n,\eps}})\big|+ \big|\psi(\cL_{Y_t^\eps})-\psi(\cL_{Y_t^{n,\eps}})\big|\Big]\no\\ &\quad+\Big[\big|\varphi(\cL_{X_t^{n,\eps}})-\varphi(\cL_{\bar X^n_t})\big|+\big|\psi(\cL_{Y_t^{n,\eps}})-\psi(\cL_{\bar Y^n_t})\big|\Big]\no\\
&\quad+\Big[\big|\varphi(\cL_{\bar X^n_t})-\varphi(\cL_{\bar X_t})\big|+\big|\psi(\cL_{\bar Y^n_t})-\psi(\cL_{\bar Y_t})\big|\Big],
\end{align}
where $(X_t^{n,\eps},Y_t^{n,\eps})$ satisfy the system (\ref{sden}), and $(\cL_{\bar X_t^n}, \cL_{\bar Y_t^n})$ are their limits obtained in Step 2. For the first term, it is easy to get that
$$
\big|\varphi(\cL_{X_t^\eps})-\varphi(\cL_{X_t^{n,\eps}})\big|+ \big|\psi(\cL_{Y_t^\eps})-\psi(\cL_{Y_t^{n,\eps}})\big|\leq C_\eps\cdot o(n),
$$
where $C_\eps>0$ is a constant depending on $\eps$, and $o(n)$ is a constant satisfying $\lim_{n\to\infty}o(n)=0$, i.e., we have for every fixed $\eps>0$,
$$
\lim_{n\to\infty}\Big[\big|\varphi(\cL_{X_t^\eps})-\varphi(\cL_{X_t^{n,\eps}})\big|+ \big|\psi(\cL_{Y_t^\eps})-\psi(\cL_{Y_t^{n,\eps}})\big|\Big]=0.
$$
By the result of Step 2, we can control the second term by
$$
\big|\varphi(\cL_{X_t^{n,\eps}})-\varphi(\cL_{\bar X^n_t})\big|+\big|\psi(\cL_{Y_t^{n,\eps}})-\psi(\cL_{\bar Y^n_t})\big|\leq C_0\cdot o(\eps),
$$
where $C_0>0$ is a constant independent of $n$, and $o(\eps)$ is  a constant satisfying $\lim_{\eps\to0}o(\eps)=0$, i.e., we have
$$
\lim_{\eps\to0}\Big[\big|\varphi(\cL_{X_t^{n,\eps}})-\varphi(\cL_{\bar X^n_t})\big|+\big|\psi(\cL_{Y_t^{n,\eps}})-\psi(\cL_{\bar Y^n_t})\big|\Big]=0.
$$
Thus, the limits $\cL_{\bar X_t}$ and $\cL_{\bar Y_t}$ are determined by the limits of $\cL_{\bar X_t^n}$ and $\cL_{\bar Y_t^n}$ as $n\to\infty$. Once we  obtain the equation of $\bar X_t$ and identify $\cL_{\bar Y_t}$, and show that
\begin{align}\label{ess}
\lim_{n\to\infty}\Big[\big|\varphi(\cL_{\bar X^n_t})-\varphi(\cL_{\bar X_t})\big|+\big|\psi(\cL_{\bar Y^n_t})-\psi(\cL_{\bar Y_t})\big|\Big]=0,
\end{align}
we can let $n\to\infty$ first and then $\eps\to0$ in (\ref{add}) to obtain the convergence of $\cL_{X_t^\eps}$ and $\cL_{Y_t^\eps}$ to $\cL_{\bar X_t}$ and $\cL_{\bar Y_t}$, respectively. This is done in subsection 6.2.

\vspace{2mm}
To conclude, we study the asymptotic behavior of the entire {\bf non-linear} stochastic system (\ref{sde1}), but each approximation system in (\ref{sden}) we dealt with is a linear one, the essential difficulty  of nonlinearity appears when taking the limits as $n\to\infty$ in (\ref{ess}) since the convergence of $\cL_{\bar X^n_t}$ relies on the convergence of  $\cL_{\bar Y^n_t}$, and vice versa.

\vspace{2mm}
\noindent{\bf Step 4.}
To prove the strong convergence of $X_t^\eps$ to $\bar X_t$, we transform the non-linear stochastic system (\ref{sde1})   into a non-autonomous linear system by freezing the distributions in the coefficients. Namely, since the system is weakly well-posed, we  define for every $t>0$ that
\begin{eqnarray*}
\begin{split}
b_\eps(t,x,y)&:=b(x,\cL_{X_t^\eps},y,\cL_{Y_t^\eps}),\\ \sigma_\eps(t,x)&:=\sigma(x,\cL_{X_t^\eps},\cL_{Y_t^\eps}),\\ F_\eps(t,x,y)&:=F(x,\cL_{X_t^\eps},y,\cL_{Y_t^\eps}),\\
 G_\eps(t,x,y)&:=G(x,\cL_{X_t^\eps},y,\cL_{Y_t^\eps}).
\end{split}
\end{eqnarray*}
Then, the McKean-Vlasov system (\ref{sde1}) can be rewritten as
\begin{equation*}
\left\{ \begin{aligned}
&\dif X^{\eps}_t =b_\eps(t,X^{\eps}_t,Y^{\eps}_t)\dif t
+\sigma_\eps(t,X^{\eps}_t)\dif W^1_t,\qquad\qquad\quad\, X^{\eps}_0=\xi,\\
&\dif Y^{\eps}_t =\frac{1}{\eps}F_\eps(t,X^{\eps}_t,Y^{\eps}_t)\dif t+\frac{1}{\sqrt{\eps}}G_\eps(t,X^{\eps}_t,Y^{\eps}_t)\dif W_t^2,\quad   Y^{\eps}_0=\eta.
\end{aligned} \right.
\end{equation*}
This system is exactly the form of (\ref{sde2}), and the convergence of the distributions of $X_t^\eps$ and $Y_t^\eps$ obtained in the previous step implies the convergence of the coefficients $b_\eps,\sigma_\eps, F_\eps$ and $G_\eps$. Thus, the strong convergence of $X_t^\eps$ to $\bar X_t$ can be obtained directly  by the strong convergence of the non-autonomous system (\ref{sde2}) established in Step 1. Such argument avoids to prove the strong convergence of the approximation system (\ref{sden}) to the original system (\ref{sde1}) (which will need the Zvonkin's transformation for the McKean-Vlasov system), see subsection 6.3 for more details.

\section{Optimal regularity for  Kolmogorov PDEs in Wasserstein space}

This section is devoted to study the optimal regularities for the solutions of two kinds of Kolmogorov equations on Wasserstein space. The first one is the backward  Kolmogorov equation on $[0,T]\times\sP_2(\mR^{d_1})$ with a fixed finite terminal time $T>0$, which will be used to prove the convergence of distribution of the slow process in (\ref{77}) (with non-linear test function $\varphi$ on $\sP_2(\mR^{d_1})$) for the non-autonomous system (\ref{sde2}). The second one is the forward Kolmogorov equation on the entire timeline and the product measure space $\mR_+\times\sP_2(\mR^{d_1}\times\mR^{d_2})$, which is essential for the proof of the convergence of the distribution of the fast process in (\ref{77}).  It is important to note that these two kinds of Kolmogorov equations we  handled are associated with the classical It\^o SDEs (see SDEs (\ref{sdet}) and (\ref{yx}) below), but not the mean-field type PDEs associated with the original non-linear system (\ref{sde1}) (whose optimal regularities for the solutions seem to be unknown).

Before proceeding further, we first provide the following chain rule formula on the Wasserstein space, which will be used several times below in our analysis.

\bl\label{chain}
Let $h:\sP_2(\mR^d)\to\mR$ admit a linear functional derivative.

\vspace{1mm}
\noindent(i) (Chain rule formula) Given function
$\Theta:\sP_2(\mR^d)\to\sP_2(\mR^d)$, then
we have
\begin{align}\label{lin2}
\frac{\delta}{\delta\mu}\big[h(\Theta(\mu))\big] (x)=\frac{\delta}{\delta\mu}\left(\int_{\mR^d}\frac{\delta h}{\delta\mu}\big(\Theta(\tilde\mu)\big)(y)\cdot\Theta(\mu)(\dif y)\right)(x)\bigg|_{\tilde\mu=\mu}.
\end{align}
In particular, if
$$
\Theta(\mu)(\dif y)=\varphi(\mu,y)\dif y,
$$
with $\varphi:\sP_2(\mR^d)\times\mR^d\to\mR$ having the linear functional derivative, then we have
\begin{align*}
\frac{\delta}{\delta\mu}\big[h(\Theta(\mu))\big] (x)=\int_{\mR^d}\frac{\delta h}{\delta\mu}\big(\Theta(\mu)\big)(y)\cdot\frac{\delta}{\delta\mu}\big[\varphi(\mu,y)\big](x)\dif y.
\end{align*}
\noindent(ii) For $i=1,2$, we define $\tilde h_i:\sP_2(\mR^d\times\mR^d)\to\mR$ by
$$
\tilde h_i(m):=h(\pi_i^* m),\quad\forall m\in \sP_2(\mR^d\times\mR^d),
$$
where $\pi_1^*m$ and $\pi_2^*m$ are the first and second marginal distributions of the probability measure $m$, respectively.
Then $\tilde h_i$ also has a linear functional derivative and
$$
\frac{\d\tilde h_i}{\d m}(m)(x_1,x_2)=\frac{\d h}{\d \mu}(\pi_i^*m)(x_i),\quad \forall x_1,x_2\in\mR^d.
$$
\el

\begin{proof}
The conclusion in $(i)$ can be proved similarly as in \cite[Proposition 2.2]{CF2}. In fact, by (\ref{lin}) we have for every $\mu_1,\mu_2\in\sP_2(\mR^{d})$,
\begin{align*}
&\theta^{-1}\Big[h\big(\Theta(\mu_1+\theta(\mu_2-\mu_1))\big)-h(\Theta(\mu_1))\Big]\\
&=\int_0^1\int_{\mR^d}\frac{\delta h}{\delta\mu}\Big(\Theta(\mu_1)+\tau\big[\Theta(\mu_1+\theta(\mu_2-\mu_1))-\Theta(\mu_1)\big]\Big)(y)\\
&\qquad\qquad\qquad\qquad\qquad\qquad \times\frac{\Theta(\mu_1+\theta(\mu_2-\mu_1))-\Theta(\mu_1)}{\theta}(\dif y)\dif \tau.
\end{align*}
Taking $\theta\to0$, by the continuity of $\mu\mapsto[\delta h/\delta\mu](\mu)(y)$ and the definition (\ref{lin0}), we obtain the formula (\ref{lin2}).
The statement of $(ii)$ is similar to \cite[Proposition A.5]{Z},  we omit the details of the proof.
\end{proof}

\subsection{Backward  Kolmogorov equation: short time estimates}
Fix $T>0$. Consider the following backward  Kolmogorov equation on $[0,T]\times\sP_2(\mR^{d_1})$:
\begin{equation}\label{ke0}
\left\{ \begin{aligned}
&\p_t u(t,\mu)+\int_{\mR^{d_1}}\Big[b(t,x)\cdot\p_\mu u(t,\mu)(x)\\
&\qquad\qquad\qquad\quad+\frac{1}{2}\T\big[\sigma\sigma^*(t,x)\cdot\p_{x}\big(\p_\mu u(t,\mu)(x)\big)\big]\Big]\mu(\dif x)=0,\\
&u(T,\mu)=\varphi(\mu),
\end{aligned} \right.
\end{equation}
where $b: [0,T]\times\mR^{d_1}\to\mR^{d_1}$, $\sigma: [0,T]\times\mR^{d_1}\to\mR^{d_1}\otimes\mR^{d_1}$ and $\varphi: \sP_2(\mR^{d_1})\to\mR$ are measurable functions. Throughout this subsection, we assume that $\sigma\sigma^*$ is uniformly elliptic, i.e., there exists a $\varrho>0$ such that for any $(t,x)\in[0,T]\times\mR^{d_1}$,
\begin{align}\label{uell}
|\sigma^*(t,x)z|^2\geq\varrho|z|^2,\qquad \forall z\in\mR^{d_1}.
\end{align}
We have the following result.

\bt\label{ke}
Assume that $\sigma, b\in C_b^{\alpha/2,\alpha}([0,T]\times\mR^{d_1})$ and $\varphi\in C_b^{(2,\a)}(\sP_2(\mR^{d_1}))$ with some $0<\a\leq2$, then there exists a unique solution $u\in C_{loc}^{1+\a/2,(2,2+\a)}\big([0,T)\times\sP_2(\mR^{d_1})\big)$ to the equation (\ref{ke0}) which is given by
\begin{align}\label{ke1}
u(t,\mu)=\varphi\big(\cL_{X_{t,T}(\xi)}\big),
\end{align}
where for $t\geq s$, $X_{s,t}(\xi)$ with $\cL_{\xi}=\mu$ is the unique weak solution of the following time inhomogeneous SDE in $\mR^{d_1}$:
\begin{align}\label{sdet}
\dif X_{s,t}=b(t,X_{s,t})\dif t+\sigma(t,X_{s,t})\dif W_t,\quad X_{s,s}=\xi.
\end{align}
Moreover, for every $0\leq s<t< T$, $\mu\in\sP_2(\mR^{d_1})$ and $x, \tilde x, \hat x\in\mR^{d_1}$, the following estimates hold: for $0<\a\leq1$,
\begin{equation} \label{esu}
		\begin{aligned}
|\p_\mu u(t,\mu)(x)|&\leq C_0\, (T-t)^{(\alpha-1)/2},\\
\big|\p_x\big(\p_\mu u(t,\mu)(x)\big)\big|&\leq C_0\, (T-t)^{\alpha/2-1},\\
|\p_\mu u(t,\mu)(x)-\p_\mu u(t,\mu)(\tilde x)|&\leq C_0\,(T-t)^{-1/2}\,\big(|x-\tilde x|^{\alpha}\wedge1\big),\\
|\p_\mu u(t,\mu)(x)-\p_\mu u(s,\mu)(x)|&\leq C_0\,(T-t)^{-1/2}\,(t-s)^{\alpha/2},\\
\bigg|\p_x\frac{\delta^2u}{\delta\mu^2}(t,\mu)(x,\tilde x)\bigg|&\leq C_0\, (T-t)^{(\alpha-1)/2},\\
\bigg|\p_x\frac{\delta^2u}{\delta\mu^2}(t,\mu)(x,\tilde x)-\p_x\frac{\delta^2u}{\delta\mu^2}(t,\mu)(x,\hat x)\bigg|&\leq C_0\, (T-t)^{-1/2}\,\big(|\tilde x-\hat x|^{\alpha}\wedge1\big),
\end{aligned}
	\end{equation}
where $C_0>0$ is a constant depending  on $T$, $\|\sigma\|_{C_b^{\a/2,\a}}$, $\|b\|_{C_b^{\a/2,\a}}$, $\|\varphi\|_{C_b^{(2,\a)}}$ and the uniformly elliptic lower bound $\varrho$ in (\ref{uell}).
\et

\begin{proof}
Let $u$ be defined by (\ref{ke1}), it is enough to prove the regularity properties of $u$, then $u$ solves equation (\ref{ke0}) follows by It\^o's formula (see e.g. \cite[Theorem 3.8]{CF2} or \cite[Theorem 7.2]{BLPR}). Let $X_{s,t}(x)$ denote the solution of the SDE (\ref{sdet}) starting from the fixed point $x\in\mR^{d_1}$ at time $s$, and $p(s,x;t,y)$ be the density function of $X_{s,t}(x)$. By the assumptions on the coefficients $\sigma$ and $b$, we have (see \cite[Chapter IV, Section 11]{LSU}) for $k=1,2$ and $0\leq s<t\leq T$,
\begin{align}\label{heat}
 |\p^k_xp(s,x;t,y)|\leq  C_1(t-s)^{-(d_1+k)/2}\exp\bigg\{-\lambda\frac{|x-y|^2}{t-s}\bigg\},
\end{align}
where $C_1, \lambda>0$ are constants depending  only on $T$, $\|\sigma\|_{C_b^{\a/2,\a}}$, $\|b\|_{C_b^{\a/2,\a}}$ and the uniformly elliptic lower bound $\varrho$ in (\ref{uell}).
Since the distribution of $\xi$ is $\mu$ and the system (\ref{sdet}) is a classical SDE, it holds that
$$
\cL_{X_{t,T}(\xi)}(\dif y)=\int_{\mR^{d_1}}\cL_{X_{t,T}(x)}(\dif y)\,\mu(\dif x)=\int_{\mR^{d_1}}p(t,x;T,y)\dif y\,\mu(\dif x).
$$
As a result, we have
\begin{align*}
u(t,\mu)=\varphi\left(\int_{\mR^{d_1}}\cL_{X_{t,T}(x)}\,\mu(\dif x)\right).
\end{align*}
Thus by the chain rule formula in Lemma \ref{chain}, we derive that
\begin{align}\label{x1}
\frac{\delta u}{\delta\mu}(t,\mu)(x_1)&=\int_{\mR^{d_1}}\frac{\delta\varphi}{\delta\mu} \big(\cL_{X_{t,T}(\xi)}\big)(x)\,\cL_{X_{t,T}(x_1)}(\dif x)
\end{align}
and
\begin{align*}
\frac{\delta^2 u}{\delta\mu^2}(t,\mu)(x_1,x_2)&=\int\!\!\!\int_{\mR^{d_1}\times\mR^{d_1}} \frac{\delta^2\varphi}{\delta\mu^2} \big(\cL_{X_{t,T}(\xi)}\big)(x,y)\,\cL_{X_{t,T}(x_2)}(\dif y)\,\cL_{X_{t,T}(x_1)}(\dif x).
\end{align*}
Define
$$
h_1(t,\mu,x):=\frac{\delta\varphi}{\delta\mu} \big(\cL_{X_{t,T}(\xi)}\big)(x)\quad\text{and}\quad h_2(t,\mu,x,y):=\frac{\delta^2\varphi}{\delta\mu^2} \big(\cL_{X_{t,T}(\xi)}\big)(x,y).
$$
Then we can write
$$
\frac{\delta u}{\delta\mu}(t,\mu)(x_1)=\mE\big[h_1\big(t,\mu,X_{t,T}(x_1)\big)\big]
$$
and
$$
\frac{\delta^2 u}{\delta\mu^2}(t,\mu)(x_1,x_2)=\mE\big[h_2\big(t,\mu,X_{t,T}(x_1),\tilde X_{t,T}(x_2)\big)\big],
$$
where $\tilde X_{t,T}$ is an independent copy of $X_{t,T}$.
By the assumption $\varphi\in C_b^{(2,\a)}(\sP_2(\mR^{d_1}))$, we have that the functions $h_1(t,\mu,\cdot)$ and $h_2(t,\mu,\cdot,\cdot)$ are bounded and $\alpha$-H\"older continuous. Thus, by the classical theory of heat kernel estimates (see  \cite[Chapter IV, Section 14]{LSU}), we have for every $0<\a\leq 1$, $t\in[0,T)$ and $\mu\in\sP_2(\mR^{d_1})$,
$$
\frac{\delta u}{\delta\mu}(t,\mu)(\cdot)\in C^{2+\alpha}_b(\mR^{d_1}),\quad\frac{\delta^2 u}{\delta\mu^2}(t,\mu)(\cdot,\cdot)\in C^{2+\alpha}_b(\mR^{d_1}\times\mR^{d_1}),
$$
which in turn implies that $u(t,\cdot)\in C_b^{(2,2+\alpha)}(\sP_2(\mR^{d_1}))$. Meanwhile, using (\ref{x1}) and estimate (\ref{heat}) we derive that
\begin{align*}
|\p_\mu u(t,\mu)(x_1)|&=\left|\p_{x_1}\frac{\delta u}{\delta\mu}(t,\mu)(x_1)\right|=\left|\int_{\mR^{d_1}}h_1(t,\mu,x)\cdot\p_{x_1}p(t,x_1;T,x)\dif x\right|\\
&\leq \int_{\mR^{d_1}}\Big|h_1(t,\mu,x)-h_1(t,\mu,x_1)\Big|\cdot\big|\p_{x_1}p(t,x_1;T,x)\big|\dif x\\
&\leq C_2\int_{\mR^{d_1}}|x-x_1|^\alpha \cdot(T-t)^{-(d_1+1)/2}\exp\bigg\{-\lambda\frac{|x-x_1|^2}{T-t}\bigg\}\dif x\\
&\leq  C_2(T-t)^{(\alpha-1)/2},
\end{align*}
and
\begin{align*}
|\p_{x_1}\p_\mu u(t,\mu)(x_1)|
&\leq \int_{\mR^{d_1}}\Big|h_1(t,\mu,x)-h_1(t,\mu,x_1)\Big|\cdot\big|\p^2_{x_1}p(t,x_1;T,x)\big|\dif x\\
&\leq C_2\int_{\mR^{d_1}}|x-x_1|^\alpha \cdot(T-t)^{-(d_1+2)/2}\exp\bigg\{-\lambda\frac{|x-x_1|^2}{T-t}\bigg\}\dif x\\
&\leq C_2(T-t)^{\alpha/2-1},
\end{align*}
where $C_2$ also depends on $\|\varphi\|_{C_b^{(1,\a)}}$. Similarly, we can prove the other estimates in (\ref{esu}) by using the H\"older continuous of the density function , see \cite[estimates (13.1)-(13.3)]{LSU}. The conclusions for $1<\a\leq2$  can be proved by the standard arguments of differentiating   the equation one time directly as in  \cite[Chapter IV, Section 5]{LSU}, the details are omitted.
\end{proof}

\br
The time singularities at the end point $T$ in estimates (\ref{esu}) are mainly caused by the low regularity assumption that $\varphi\in C_b^{(2,\alpha)}(\sP_2(\mR^{d_1}))$ (which is even not Lipschitz continuous with respect to the Wasserstein distance). If we assume that $\varphi\in C_b^{(2,2+\alpha)}(\sP_2(\mR^{d_1}))$, then we could get  $u\in C_b^{1+\alpha/2,(2,2+\alpha)}([0,T]\times\sP_2(\mR^{d_1}))$.
\er

Using Theorem \ref{ke}, we can  get the following continuous dependence with respect to the coefficients for the distributions of the solutions of SDE (\ref{sdet}), which will be used to take the limit of the averaged equations of the approximation systems in Section 6. Namely, let $X_{s,t}$ be the solution of SDE (\ref{sdet}), and $\tilde X_{s,t}$ satisfies SDE (\ref{sdet}) with coefficients $\tilde \sigma$ and $\tilde b$, i.e.
\begin{align*}
\dif \tilde X_{s,t}=\tilde b(t,\tilde X_{s,t})\dif t+\tilde \sigma(t,\tilde X_{s,t})\dif W_t,\quad \tilde X_{s,s}=\xi.
\end{align*}
We have the following result.

\bl\label{diff}
Assume that (\ref{uell}) hold, and $\sigma, \tilde\sigma, b, \tilde b\in C_b^{\alpha/2,\alpha}([0,T]\times\mR^{d_1})$ with $0<\a\leq 1$. Then for every $t\in[0,T]$ and $\varphi\in C_b^{(2,\a)}(\sP_2(\mR^{d_1}))$, we  have
\begin{align*}
\big|\varphi\big(\cL_{X_{t,T}(\xi)}\big)-\varphi\big(\cL_{\tilde X_{t,T}(\xi)}\big)\big|\leq \hat C_T\Big(\|b-\tilde b\|_\infty+\|\sigma-\tilde\sigma\|_\infty\Big),
\end{align*}
where $\hat C_T>0$ is a constant with $\lim_{T\to0}\hat C_T=0$.
\el

\begin{proof}
Let $u$ be defined by (\ref{ke1}), and $\tilde u$ satisfy the equation (\ref{ke0}) with coefficients $\tilde \sigma$ and $\tilde b$, i.e.,
\begin{equation*}
\left\{ \begin{aligned}
&\p_t \tilde u(t,\mu)+\int_{\mR^{d_1}}\Big[\tilde b(t,x)\cdot\p_\mu \tilde u(t,\mu)(x)\\
&\qquad\qquad\qquad\quad+\frac{1}{2}\T\big[\tilde \sigma\tilde \sigma^*(t,x)\cdot\p_{x}\big(\p_\mu \tilde u(t,\mu)(x)\big)\big]\Big]\mu(\dif x)=0,\\
&\tilde u(T,\mu)=\varphi(\mu).
\end{aligned} \right.
\end{equation*}
Define
$$
v(t,\mu):=u(t,\mu)-\tilde u(t,\mu).
$$
Then by (\ref{ke1}) we have
$$
v(t,\mu)=\varphi\big(\cL_{X_{t,T}(\xi)}\big)-\varphi\big(\cL_{\tilde X_{t,T}(\xi)}\big),
$$
and
\begin{equation*}
\left\{ \begin{aligned}
&\p_t v(t,\mu)+\int_{\mR^{d_1}}\Big[\frac{1}{2}\T\big[\sigma\sigma^*(t,x)\cdot\p_{x}\big(\p_\mu v(t,\mu)(x)\big)\big]\\
&\qquad\quad\,\,\,+b(t,x)\cdot\p_\mu v(t,\mu)(x)\Big]\mu(\dif x)=\int_{\mR^{d_1}}\Big[\big(\tilde b(t,x)-b(t,x)\big)\cdot\p_\mu \tilde u(t,\mu)(x)\\
&\qquad\quad\,\,\,+\frac{1}{2}\T\big[\big(\tilde\sigma\tilde\sigma^*(t,x)-\sigma\sigma^*(t,x)\big) \cdot\p_{x}\big(\p_\mu \tilde u(t,\mu)(x)\big)\big]\Big]\mu(\dif x)=:f(t,\mu),\\
&v(T,\mu)=0.
\end{aligned} \right.
\end{equation*}
By the assumptions on the coefficients and the regularities of $\tilde u$ obtained in Theorem \ref{ke}, one can check that for every $t\in[0,T)$, $f(t,\cdot)\in C_b^{(1,\alpha)}$.
Thus, by \cite[Theorem 3.8]{CF2} we have
\begin{align}\label{00}
v(t,\mu)=-\int_t^Tf\big(s,\cL_{X_{t,s}(\xi)}\big)\dif s,
\end{align}
where $X_{t,s}(\xi)$ satisfies the equation (\ref{sdet}) with $\cL_\xi=\mu$ and the initial time $t$. As a result of the first two estimates in (\ref{esu}), we deduce that
\begin{align}\label{002}
|v(t,\mu)|&\leq C_0\Big(\|b-\tilde b\|_\infty+\|\sigma-\tilde\sigma\|_\infty\Big)\int_t^T\!\!\Big(\|\p_\mu \tilde u(s,\cdot)(\cdot)\|_\infty+\|\p_{x}\big(\p_\mu \tilde u(s,\cdot)(\cdot)\big)\|_\infty\Big)\dif s\no\\
&\leq C_0\Big(\|b-\tilde b\|_\infty+\|\sigma-\tilde\sigma\|_\infty\Big)\int_t^T(T-s)^{\a/2-1}\dif s,
\end{align}
which in turn implies the desired result.
\end{proof}

\subsection{Forward Kolmogorov equation: exponential decay of the solution}

Consider the following parameterized SDE in $\mR^{d_2}$:
\begin{align} \label{yx}
\dif Y^{t,x}_s =F(t,x,Y^{t,x}_s)\dif s+G(t,x,Y^{t,x}_s)\dif W_s,\qquad   Y^{t,x}_0=y\in\mR^{d_2},
\end{align}
where $(t,x)\in[0,T]\times\mR^{d_1}$ are regarded as parameters, $F: [0,T]\times\mR^{d_1}\times \mR^{d_2}\to\mR^{d_2}$, $G: [0,T]\times\mR^{d_1}\times \mR^{d_2}\to\mR^{d_2}\otimes\mR^{d_2}$ are measurable functions. We  make the following assumption on the coefficients:

\vspace{0.1cm}
\begin{description}
\item[$(\mathbf{A}_{FG})$] the coefficient  $GG^*(t,x,y)$ is non-degenerate in the sense that there exist constants $k,\varrho>0$ such that for every $(t,x,y)\in [0,T]\times\mR^{d_1}\times \mR^{d_2}$,
	$$
	\varrho(1+|y|)^{-k}|z|^{2} \leq |G^*(t,x,y)z|^2, \quad \forall z\in \mathbb{R}^{d_{2}},
	$$
and
there exist constants $\Lambda_1, \Lambda_2>0$ such that for every $(t,x,y)\in[0,T]\times\mR^{d_1}\times \mR^{d_2}$,
\begin{align*}
2\<y,F(t,x,y)\>+\|G(t,x,y)\|^2\leq -\Lambda_1|y|^2+\Lambda_2.
\end{align*}
\end{description}

\noindent Under $(\mathbf{A}_{FG})$, there exists a  unique invariant measure $\zeta^{t,x}(\dif y)$ for  the system (\ref{yx}) which is $\cV$-exponential ergodic with $\cV(y)= 1+|y|^p$ for every $p\geq 1$ (see e.g. \cite[Theorem 7.4]{XZ} or \cite[Remark 2.1]{EGZ}), i.e., for every $s\geq 0$ and $(t,x)\in[0,T]\times\mR^{d_1}$,
\begin{align}\label{erg}
\rho_\cV\big(\cL_{Y^{t,x}_s(y)},\zeta^{t,x}\big)\leq \Lambda_0\big(1+\cV(y)\big)\,\e^{-\gamma s},
\end{align}
where $\gamma$ and $\Lambda_0$ are positive constants depending only on $k, \varrho, \Lambda_1$ and $\Lambda_2$ in the assumption $(\mathbf{A}_{FG})$, and $\rho_\cV$ is defined by (\ref{tv}).

To prove the convergence of the distribution of the fast process of the form (\ref{es2}), it turns out to be essential to consider the following Kolmogorov equation on the entire timeline and the product measure space $\mR_+\times\sP_2(\mR^{d_1}\times\mR^{d_2})\times[0,T]$:
\begin{equation}\label{cp0}
\left\{\begin{aligned}
\displaystyle
&\p_s V(s,m;t)-\int\!\!\!\int_{\mR^{d_1}\times\mR^{d_2}}\bigg[F(t,x,y)\cdot\p_y\frac{\delta V}{\delta m}(s,m;t)(x,y)\\
\displaystyle
&\qquad\qquad\quad\,+\frac{1}{2}\T\Big(GG^*(t,x,y)\cdot\p^2_y\frac{\delta V}{\delta m}(s,m;t)(x,y)\Big)\bigg]m(\dif x,\dif y)=0,\\
\displaystyle
&V(0,m;t)=\psi(\pi_2^*m)-\psi(\tilde \zeta^{t,\pi_1^*m}),
\end{aligned}\right.
\end{equation}
where $t\in[0,T]$ is regarded as a parameter, $\psi: \sP_2(\mR^{d_2})\to\mR$ is a  measurable function, $\pi_1^*m$ and $\pi_2^*m$ are the first and second marginal distributions of the probability measure $m\in\sP_2(\mR^{d_1}\times\mR^{d_2})$, respectively, and the measure $\tilde \zeta^{t,\mu}$ with $(t,\mu)\in[0,T]\times\sP_2(\mR^{d_1})$ is defined by
\begin{align}\label{zemu}
\tilde \zeta^{t,\mu}(\dif y)=\int_{\mR^{d_1}}\zeta^{t,x}(\dif y)\mu(\dif x),
\end{align}
where $\zeta^{t,x}(\dif y)$ is the unique invariant measure for the system (\ref{yx}). The key point is to study the long time behavior of the solution of the equation (\ref{cp0}).
We have the following result.

\bt\label{cp}
 Let $(\mathbf{A}_{FG}) $ hold. Assume that  $F,G\in C^{\a/2,\a,\b}_p([0,T]\times\mR^{d_1}\times\mR^{d_2})$ and $\psi\in C^{(2,\b)}_p(\sP(\mR^{d_2}))$ with some $0<\a,\b\leq2$, then there exists a unique solution $V(s,m;t)$ to the equation (\ref{cp0}) which is given by
\begin{align}\label{VV}
V(s,m;t):=\psi\big(\cL_{Y^{t,\xi}_s(\eta)}\big)-\psi\big(\tilde \zeta^{t,\mu}\big),
\end{align}
where $\xi$ and $\eta$ are two random variables with $\cL_{(\xi,\eta)}=m$, $\mu=\pi_1^*m$, $\tilde \zeta^{t,\mu}$ is defined by (\ref{zemu}), and $Y^{t,\xi}_s(\eta)$ is the unique weak solution of the following stochastic system:
\begin{align} \label{yxi1}
\dif Y^{t,\xi}_s =F(t,\xi,Y^{t,\xi}_s)\dif s+G(t,\xi,Y^{t,\xi}_s)\dif W_s,\qquad   Y^{t,\xi}_0=\eta.
\end{align}
Moreover, there exists $p\geq1$ such that for $0<\beta\leq 1$, the following estimates hold:

\vspace{1mm}
\noindent(i) (Estimates of the solution $V$) for every $s>0$, $m\in\sP_2(\mR^{d_1}\times\mR^{d_2})$, $t\in[0,T]$ and $(x,y)\in\mR^{d_1}\times\mR^{d_2}$, we have
\begin{align}
&|V(s,m;t)|\leq C_0\,\e^{-\gamma s},\label{V1}\\
&\Big|\p_y\frac{\delta V}{\delta m}(s,m;t)(x,y)\Big|\leq C_0(1+|y|^p)s^{\frac{\b-1}{2}}\e^{-\gamma s},\label{V2}\\
&\Big|\p_y^2\frac{\delta V}{\delta m}(s,m;t)(x,y)\Big|\leq C_0(1+|y|^p)s^{\frac{\b}{2}-1}\e^{-\gamma s},\label{V3}
\end{align}
(ii) (Estimates  of $\delta V/\delta m$ w.r.t. the variable $x$) for every $x_1,x_2\in\mR^{d_1}$, we have
\begin{align}\label{vx}
\begin{split}
&\Big|\frac{\delta V}{\delta m}(s,m;t)(x_1,y)-\frac{\delta V}{\delta m}(s,m;t)(x_2,y)\Big|\leq C_0(1+|y|^p)\e^{-\gamma s}|x_1-x_2|^\a,\\
&\Big|\p_y\frac{\delta^2 V}{\delta m^2}(s,m;t)(x,y)(x_1,\tilde y)\\
&\qquad\,\qquad-\p_y\frac{\delta^2 V}{\delta m^2}(s,m;t)(x,y)(x_2,\tilde y)\Big|\leq C_0(1+|y|^p)s^{\frac{\b-1}{2}}\e^{-\gamma s}|x_1-x_2|^\a,\\
&\Big|\p_y^2\frac{\delta^2 V}{\delta m^2}(s,m;t)(x,y)(\tilde x_1,\tilde y)\\
&\qquad\,\qquad-\p_y^2\frac{\delta^2 V}{\delta m^2}(s,m;t)(x,y)(\tilde x_2,\tilde y)\Big|\leq C_0(1+|y|^p)s^{\frac{\b}{2}-1}\e^{-\gamma s}|\tilde x_1-\tilde x_2|^\a,\\
&\Big|\p_s\frac{\delta V}{\delta m}(s,m;t)(x_1,y)\!-\!\p_s\frac{\delta V}{\delta m}(s,m;t)(x_2,y)\Big|\leq C_0(1+|y|^p)s^{\frac{\b}{2}-1}\e^{-\gamma s}|x_1-x_2|^\a,
\end{split}
\end{align}
(iii) (Estimates  of $V$ w.r.t. the parameter $t$) for every $t_1, t_2\in[0,T]$, we have
\begin{align}\label{vt}
\begin{split}
&|V(s,m;t_1)-V(s,m;t_2)|\leq C_0\,\e^{-\gamma s}|t_1-t_2|^\frac{\a}{2},\\
&|\p_sV(s,m;t_1)-\p_sV(s,m;t_2)|\leq C_0\,s^{\frac{\b}{2}-1}\e^{-\gamma s}|t_1-t_2|^\frac{\a}{2},\\
&\Big|\p_y\frac{\delta V}{\delta m}(s,m;t_1)(x_1,y)-\p_y\frac{\delta V}{\delta m}(s,m;t_2)(x_2,y)\Big|\\
&\qquad\,\qquad\leq C_0(1+|y|^p)s^{\frac{\b-1}{2}}\e^{-\gamma s}\big(|t_1-t_2|^\frac{\a}{2}+|x_1-x_2|^\a\big),\\
&\Big|\p_y^2\frac{\delta V}{\delta m}(s,m;t_1)(x_1,y)-\p_y^2\frac{\delta V}{\delta m}(s,m;t_2)(x_2,y)\Big|\\
&\qquad\,\qquad\leq C_0(1+|y|^p)s^{\frac{\b}{2}-1}\e^{-\gamma s}\big(|t_1-t_2|^\frac{\a}{2}+|x_1-x_2|^\a\big),
\end{split}
\end{align}
where $C_0, \gamma>0$ are constants depending on the norms of the coefficients and the constants in the assumption $(\mathbf{A}_{FG})$.
\et

\br\label{small}
i) We shall show in Lemma \ref{inv-mea} that the distribution of $Y_s^{t,\xi}(\eta)$ depends on $(\xi,\eta)$ only through their joint distribution $m=\cL_{(\xi,\eta)}$, thus the function $V$ in (\ref{VV}) is well-defined.

\vspace{1mm}
ii) The exponential decay in the $s$-variable of the estimates (\ref{V1})-(\ref{vt}) is mainly due to the dissipative assumption $(\mathbf{A}_{FG})$.  The higher order derivatives with respect to the $y$-variable in the estimates (\ref{V2}) and (\ref{V3}) comes from the local elliptic property of the differential operator in the integral part of the equation (\ref{cp0}).  As in Theorem \ref{ke}, the time singularities at the starting point in the above  estimates involving the derivatives with respect to the $y$-variable
are caused by the low regularity assumption that $\psi\in C^{(2,\b)}_p(\sP(\mR^{d_2}))$.

\vspace{1mm}
iii) Note that the differential operator in the equation (\ref{cp0}) can be viewed as degenerate with respect to the $x$-variable, thus no regularity improvement occurs to the $x$-variable. Meanwhile, $t$ is only a parameter in the equation. Consequently, the solution $V$ remains the same regularities with respect to $(x,t)$-variables as the coefficients, see estimates in (\ref{vx}) and (\ref{vt}).

iv) From the proof below, it can be seen that if $\psi$ satisfies (\ref{dissf2}), then the constant $C_0$ in the  estimates (\ref{V1})-(\ref{vt}) can be replaced by $\kappa\, C_0$.
\er

The proof of Theorem \ref{cp} replies on the study of the long time behavior of $Y_s^{t,\xi}(\eta)$  which satisfies the system (\ref{yxi1}). We point out that it is not enough to consider $Y_s^{t,x}(y)$ which satisfies (\ref{yx}) since
 $\xi$ and $Y^{t,\xi}_s(\eta)$ are obviously not independent ($\xi$ and $\eta$ are not independent).
Let us first establish the following result.

\bl\label{inv-mea}
Let $(\mathbf{A}_{FG})$ hold, and $F,G\in C^{\a/2,\a,\b}_p([0,T]\times\mR^{d_1}\times\mR^{d_2})$ with some $0<\a,\b\leq1$.
Then for every $(\xi,\eta)$, the system (\ref{yxi1}) has a unique weak solution, and the distribution of $Y_s^{t,\xi}(\eta)$ depends on $(\xi,\eta)$ only through their joint distribution $m=\cL_{(\xi,\eta)}$. Moreover,   we have for any $s\geq0$,
\begin{align*}
\rho_\cV\big(\cL_{Y_s^{t,\xi}(\eta)},\tilde\zeta^{t,\mu}\big)\leq \Lambda_0\,\e^{-\gamma s},
\end{align*}
where $\tilde \zeta^{t,\mu}$ is defined by (\ref{zemu}) with $\mu=\cL_\xi$,   $\cV(y)= 1+|y|^p$ with $p\geq 1$, and  $\Lambda_0, \gamma>0$ are constants  depending only on $k, \varrho, \Lambda_1$ and $\Lambda_2$ in the assumption $(\mathbf{A}_{FG})$.
\el
\begin{proof}
We rewrite the equation (\ref{yxi1}) into  the following coupled degenerate stochastic system: for $s>0$,
\begin{equation}\label{XY}
\left\{ \begin{aligned}
&\dif X_s=0,\qquad\qquad\qquad\qquad\qquad\qquad\qquad\,\, X_0=\xi,\\
&\dif Y_s =F(t,X_s,Y_s)\dif s+G(t,X_s,Y_s)\dif W_s,\quad\,   Y_0=\eta,
\end{aligned} \right.
\end{equation}
where $t\in[0,T]$ is a parameter. It is obvious that $Y_s^{t,\xi}(\eta)=Y_s$. Since the above system is a classical SDE,  it is enough to consider
the solution of the system (\ref{XY}) with determinate initial value $(x,y)$, which is denoted by $(X^t_s(x,y),Y^t_s(x,y))$, i.e.,
\begin{equation}\label{dec}
\left\{ \begin{aligned}
&\dif X^t_s(x,y)=0,\qquad\qquad\qquad\qquad\qquad\qquad\qquad\qquad\quad\quad X^t_0(x,y)=x,\\
&\dif Y^t_s(x,y) =F\big(t,X^t_s(x,y),Y^t_s(x,y)\big)\dif s\\
&\qquad\qquad\qquad\qquad+G\big(t,X^t_s(x,y),Y^t_s(x,y)\big)\dif W_s,\qquad\quad   Y^t_0(x,y)=y.
\end{aligned} \right.
\end{equation}
Note that $Y_s^t(x,y)$ is the same as $Y_s^{t,x}(y)$ which satisfies (\ref{yx}). Then by the theory of classical SDEs  we have
\begin{align*}
\cL_{(X_s,Y_s)}(\dif\tilde x,\dif\tilde y)&=\int\!\!\!\int_{\mR^{d_1}\times\mR^{d_2}}\cL_{(X^t_s(x,y),Y^t_s(x,y))}(\dif\tilde x,\dif\tilde y)m(\dif x,\dif y)\\
&=\int\!\!\!\int_{\mR^{d_1}\times\mR^{d_2}}\delta_x(\dif\tilde x)\times\cL_{Y_s^{t,x}(y)}(\dif\tilde y)m(\dif x,\dif y),
\end{align*}
where $(X_s,Y_s)$ is the solution of the system (\ref{XY}) and $m=\cL_{(\xi,\eta)}$. Thus the weak well-posedness of (\ref{XY}) follows directly. In particular, we have
\begin{align*}
\cL_{Y_s^{t,\xi}(\eta)}(\dif\tilde y)&=\cL_{Y_s}(\dif\tilde y)=\cL_{(X_s,Y_s)}(\mR^{d_1},\dif\tilde y)\\
&=\int\!\!\!\int_{\mR^{d_1}\times\mR^{d_2}}\cL_{Y_s^{t,x}(y)}(\dif\tilde y)m(\dif x,\dif y).
\end{align*}
Recall that $\zeta^{t,x}$ is the unique invariant measure of the equation (\ref{yx}) and satisfies (\ref{erg}).
By the Minkowski inequality and the estimate (\ref{erg}) we deduce that
\begin{align*}
\rho_\cV\big(\cL_{Y_s^{t,\xi}(\eta)},\tilde\zeta^{t,\mu}\big)
&=\rho_\cV\bigg(\int\!\!\!\int_{\mR^{d_1}\times\mR^{d_2}}\cL_{Y_s^{t,x}(y)}m(\dif x,\dif y),\int_{\mR^{d_1}}\zeta^{t,x}\mu(\dif x)\bigg)\\
&=\rho_\cV\bigg(\int\!\!\!\int_{\mR^{d_1}\times\mR^{d_2}}\cL_{Y_s^{t,x}(y)}m(\dif x,\dif y),\int\!\!\!\int_{\mR^{d_1}\times\mR^{d_2}}\zeta^{t,x}m(\dif x,\dif y)\bigg)\\
&\leq\int\!\!\!\int_{\mR^{d_1}\times\mR^{d_2}}\rho_\cV\big(\cL_{Y_s^{t,x}(y)},\zeta^{t,x}\big)m(\dif x,\dif y)\\
&\leq\int\!\!\!\int_{\mR^{d_1}\times\mR^{d_2}}\big(1+\cV(y)\big)m(\dif x,\dif y)\cdot\Lambda_0\,\e^{-\gamma s}\leq\Lambda_0\,\e^{-\gamma s}.
\end{align*}
The proof is finished.
\end{proof}

Now we proceed to give:

\begin{proof}[{\bf Proof of Theorem \ref{cp}}]
By regarding $t$ as a parameter, the operator in the equation (\ref{cp0}) can be viewed as the generator of the SDE (\ref{XY}), thus the solution $V$ should have the probability representation (see e.g. \cite[Theorem 3.8]{CF2} or \cite[Theorem 7.2]{BLPR}) that
\begin{align*}
V(s,m;t)=\psi\big(\pi_2^*\cL_{(X_s,Y_s)}\big)-\psi\big(\tilde\zeta^{t,\pi_1^*\cL_{(X_s,Y_s)}}\big),
\end{align*}
where $\cL_{(\xi,\eta)}=m$.
Note that the system (\ref{XY}) is equivalent to (\ref{yxi1}), i.e., we have
\begin{align*}
\pi_2^*\cL_{(X_s,Y_s)}=\cL_{Y_s}=\cL_{Y_s^{t,\xi}(\eta)}\quad\text{and}\quad
\pi_1^*\cL_{(X_s,Y_s)}=\cL_{X_s}=\cL_\xi=\mu.
\end{align*}
As a result, we have (\ref{VV}). Using the estimate in Lemma \ref{inv-mea} and the assumption on $\psi$, we derive that for some $\theta\in(0,1)$ and $p>1$,
\begin{align}\label{ss}
|V(s,m;t)|&=\bigg|\int_{\mR^{d_2}}\frac{\delta\psi}{\delta\nu}\Big(\cL_{Y_s^{t,\xi}(\eta)}+\theta\big(\tilde \zeta^{t,\mu}-\cL_{Y_s^{t,\xi}(\eta)}\big)\Big)(y)\big(\cL_{Y_s^{t,\xi}(\eta)}-\tilde \zeta^{t,\mu}\big)(\dif y)\bigg|\no\\
&\leq \|\delta\psi/\delta\nu\|_{L^\infty_p}\int_{\mR^{d_2}}(1+|y|^p)\big|\cL_{Y_s^{t,\xi}(\eta)}-\tilde \zeta^{t,\mu}\big|(\dif y)\no\\
&\leq C_0\,\rho_{\cV}\big(\cL_{Y_s^{t,\xi}(\eta)},\tilde\zeta^{t,\mu}\big)\leq C_0\,\e^{-\gamma s},
\end{align}
which yields the estimate (\ref{V1}).
Below, we proceed to prove the a-priori estimates for the solution.
Applying Lemma \ref{chain}, we get
\begin{align*}
&\quad \frac{\delta V}{\delta m}(0,m;t)(x,y)\\
&=\frac{\delta \psi}{\delta \nu}(\pi_2^*m)(y)-\frac{\delta}{\delta m}\bigg(\int_{\mR^{d_1}}\int_{\mR^{d_2}}\frac{\delta \psi}{\delta \nu}\big(\tilde \zeta^{t,\pi_1^*{\tilde m}}\big)(\tilde y)\zeta^{t,\tilde x}(\dif \tilde y)\pi_1^*m(\dif \tilde x)\bigg)(x,y)\bigg|_{\tilde m=m}\\
&=\frac{\delta \psi}{\delta \nu}(\pi_2^*m)(y)-\int_{\mR^{d_2}}\frac{\delta \psi}{\delta \nu}\big(\tilde \zeta^{t,\pi_1^*m}\big)(\tilde y)\zeta^{t,x}(\dif \tilde y).
  \end{align*}
Taking linear functional derivative with respect to $m$ from both sides of equation (\ref{cp0}), we have that $[\delta V/\delta m](s,m;t)(x,y)$ satisfies the following equation:
\begin{equation}\label{cp00}
\left\{\begin{aligned}
\displaystyle
&\p_s \frac{\delta V}{\delta m}(s,m;t)(x,y)-F(t,x,y)\cdot\p_y\frac{\delta V}{\delta m}(s,m;t)(x,y)\\
&\qquad\quad-\frac{1}{2}\T\Big(\cG(t,x,y)\cdot\p^2_y\frac{\delta V}{\delta m}(s,m;t)(x,y)\Big)\\
&\qquad\quad-\int\!\!\!\int_{\mR^{d_1}\times\mR^{d_2}}\bigg[F(t,\tilde x,\tilde y)\cdot\p_{\tilde y}\frac{\delta }{\delta m}\bigg(\frac{\delta V}{\delta m}\bigg)(s,m;t)(\tilde x,\tilde y)(x,y)\\
\displaystyle
&\qquad\qquad\quad\,+\frac{1}{2}\T\Big(\cG(t,\tilde x,\tilde y)\cdot\p^2_{\tilde y}\frac{\delta }{\delta m}\bigg(\frac{\delta V}{\delta m}\bigg)(s,m;t)(\tilde x,\tilde y)(x,y)\Big)\bigg]m(\dif \tilde x,\dif \tilde y)=0,\\
\displaystyle
&\frac{\delta V}{\delta m}(0,m;t)(x,y)=\frac{\delta \psi}{\delta \nu}(\pi_2^*m)(y)-\int_{\mR^{d_2}}\frac{\delta \psi}{\delta \nu}\big(\tilde \zeta^{t,\pi_1^*m}\big)(\tilde y)\zeta^{t,x}(\dif \tilde y).
\end{aligned}\right.
\end{equation}
The formulation (\ref{cp00}) falls into the equation  considered in \cite[Theorem 3.8]{CF2}. Consequently, by the assumptions on the coefficients, we have that $[\delta V/\delta m](s,m;t)(x,y)$ admits the probability representation that
$$
\frac{\delta V}{\delta m}(s,m;t)(x,y)=\mE\frac{\delta \psi}{\delta \nu}\big(\cL_{Y_s^{t,\xi}(\eta)}\big)\big(Y_s^t(x,y)\big)-\int_{\mR^{d_2}}\frac{\delta \psi}{\delta \nu}\big(\tilde \zeta^{t,\mu}\big)(\tilde y)\zeta^{t,x}(\dif \tilde y),
$$
where $Y_s^{t,\xi}(\eta)$ is the solution of the system (\ref{yxi1}), and  $Y_s^t(x,y)$ satisfies the decoupled equation (\ref{dec}). Writing
\begin{align*}
\Big|\frac{\delta V}{\delta m}(s,m;t)(x,y)\Big|&\leq \bigg|\mE\frac{\delta \psi}{\delta \nu}\big(\cL_{Y_s^{t,\xi}(\eta)}\big)\big(Y_s^t(x,y)\big)-\mE\frac{\delta \psi}{\delta \nu}\big(\tilde \zeta^{t,\mu}\big)\big(Y_s^t(x,y)\big)\bigg|\\
&\quad+\bigg|\mE\frac{\delta \psi}{\delta \nu}\big(\tilde \zeta^{t,\mu}\big)\big(Y_s^t(x,y)\big)-\int_{\mR^{d_2}}\frac{\delta \psi}{\delta \nu}\big(\tilde \zeta^{t,\mu}\big)(\tilde y)\zeta^{t,x}(\dif \tilde y)\bigg|,
\end{align*}
and using the similar argument as in the proof of (\ref{V1}), we have
$$
\Big|\frac{\delta V}{\delta m}(s,m;t)(x,y)\Big|\leq C_0(1+|y|^p)\e^{-\gamma s}.
$$
Furthermore, we write for every fixed $y\in\mR^{d_2}$,
$$
\tilde f(y):=\frac{\delta \psi}{\delta \nu}\big(\cL_{Y_s^{t,\xi}(\eta)}\big)(y)-\int_{\mR^{d_2}}\frac{\delta \psi}{\delta \nu}\big(\tilde \zeta^{t,\mu}\big)(\tilde y)\zeta^{t,x}(\dif \tilde y).
$$
Then we have
$$
\frac{\delta V}{\delta m}(s,m;t)(x,y)=\mE\tilde f\big(Y_s^{t,x}(y)\big).
$$
Since $\tilde f$ is H\"older continuous with respect to the $y$-variable, and $Y_s^{t,x}(y)$ satisfies the time-homogeneous equation (\ref{yx}), we have by \cite[Theorem 2.1, Corollary 2.2]{X} that estimates (\ref{V2}) and (\ref{V3}) hold. The regularities with respect to the $x$-variable and the parameter $t$ in estimates (\ref{vx}) and (\ref{vt}) can be proved by the same arguments as in \cite[Theorem 2.1]{RX1}, for the sake of simplicity, we omit the details.
\end{proof}

We shall also need  the following uniform in time continuous dependence with respect to the coefficients for the distributions of the solutions of time-homogeneous SDEs with dissipative coefficients, which is similar to Lemma \ref{diff} and will be used to take the limit of the frozen equations of the approximation systems in Section 6. Namely, let $Y_s(\eta)$   and $\tilde Y_s(\eta)$ satisfy
\begin{align*}
\dif Y_s=F(Y_s)\dif s+G(Y_s)\dif W_s,\quad Y_0=\eta,
\end{align*}
and
\begin{align*}
\dif \tilde Y_s=\tilde F(\tilde Y_s)\dif t+\tilde G(\tilde Y_s)\dif W_s,\quad \tilde Y_0=\eta.
\end{align*}
Assume that there exist constants $k,\varrho>0$ such that
	\begin{align}\label{32}
	\varrho(1+|y|)^{-k}|z|^{2} \leq |G^*(y)z|^2\wedge|\tilde G^*(y)z|^2, \quad \forall z\in \mathbb{R}^{d_{2}},
	\end{align}
and for any $q\geq2$, there exist constants $\Lambda_1, \Lambda_2>0$ such that
\begin{align}\label{33}
\Big(2\<y,F(y)\>+(q-1)\|G(y)\|^2\Big)\vee\Big(2\<y,\tilde F(y)\>+(q-1)\|\tilde G(y)\|^2\Big)\leq -\Lambda_1|y|^2+\Lambda_2.
\end{align}
We have the following result.

\bl\label{diff2}
Assume that (\ref{32}), (\ref{33}) hold, and $F, \tilde F, G, \tilde G\in C_p^{\beta}(\mR^{d_2})$ with $0<\beta\leq 1$. Then, for every $s>0$ and $\psi\in C_p^{(2,\beta)}(\sP_2(\mR^{d_2}))$, we have
\begin{align*}
\big|\psi(\cL_{Y_s(\eta)})-\psi(\cL_{\tilde Y_s(\eta)})\big|&\leq C_d\Big(\|F-\tilde F\|_{L_p^\infty}+\|G-\tilde G\|_{L_p^\infty}\Big),
\end{align*}
where $C_d>0$ is a constant independent of $s$.
\el
\begin{proof}
We give the main difference between the proof of Lemma \ref{diff}. Let us define
$$
\tilde u(s,\nu):=\psi(\cL_{\tilde Y_s(\eta)}),
$$
where $\nu=\cL_\eta$. Then we may argue as in the derivation of  (\ref{00}) and (\ref{002}) to get that there exists a $p\geq 1$ such that
\begin{align*}
\big|\psi(\cL_{Y_s(\eta)})-&\psi(\cL_{\tilde Y_s(\eta)})\big|\leq C_d\,\mE\bigg(\int_0^s\Big[\big|F(Y_r(\eta))-\tilde F(Y_r(\eta))\big|\cdot\|\p_\nu u(r,\cdot)(\cdot)\|_{L_p^\infty}\\
&+\big|G(Y_r(\eta))-\tilde G(Y_r(\eta))\big|\cdot\|\p_y\big(\p_\nu u(r,\cdot)(\cdot)\big)\|_{L_p^\infty}\Big]\cdot\big(1+|Y_r(\eta)|^p\big)\dif r\bigg).
\end{align*}
Under the assumption (\ref{33}), there exists a unique invariant measure $\zeta$ (which is independent of $\nu$) for the process $\tilde Y_s(\eta)$. Let
$$
\hat u(s,\nu):=\psi(\cL_{Y_s(\eta)})-\psi(\zeta).
$$
Then we have
$$
\p_\nu u(s,\nu)(y)=\p_\nu \hat u(s,\nu)(y).
$$
As in the proof of estimates (\ref{V2}) and (\ref{V3}) (see also \cite[Corollary 2.2]{X}), we have that there exist  $C_0, \gamma>0$ such that for every $r>0$,
$$
\|\p_\nu u(r,\cdot)(\cdot)\|_{L_p^\infty}=\|\p_\nu \hat u(r,\cdot)(\cdot)\|_{L_p^\infty}\leq C_0\,r^{\frac{\beta-1}{2}}\e^{-\gamma r},
$$
and
$$
\|\p_y\big(\p_\nu u(r,\cdot)(\cdot)\big)\|_{L_p^\infty}=\|\p_y\big(\p_\nu \hat u(r,\cdot)(\cdot)\big)\|_{L_p^\infty}\leq C_0\,r^{\frac{\beta}{2}-1}\e^{-\gamma r}.
$$
As a result, we further obtain
\begin{align*}
\big|\psi(\cL_{Y_s(\eta)})-\psi(\cL_{\tilde Y_s(\eta)})\big|
&\leq C_d\,\Big(\|F-\tilde F\|_{L_p^\infty}+\|G-\tilde G\|_{L_p^\infty}\Big)\\
&\qquad\times\bigg(\int_0^s\mE\big(1+|Y_r(\eta)|^{2p}\big)\cdot r^{\frac{\beta}{2}-1}\e^{-\gamma r}\dif r\bigg)
\end{align*}
Under (\ref{33}), we have that for every $q\geq 2$,
$$
\sup_{r\geq 0}\mE\big(1+|Y_r(\eta)|^q\big)\leq C_1<\infty.
$$
Thus we deduce that
\begin{align*}
\big|\psi(\cL_{Y_s(\eta)})-\psi(\cL_{\tilde Y_s(\eta)})\big|
&\leq C_d\,\Big(\|F-\tilde F\|_{L_p^\infty}+\|G-\tilde G\|_{L_p^\infty}\Big)\cdot\int_0^s r^{\frac{\beta}{2}-1}\e^{-\gamma r}\dif r,
\end{align*}
which in turn implies the desired result.
\end{proof}

\section{Poisson equation and mollifying on Wasserstein space}

In this section, we first provide the result of Poisson equation with parameters in the whole space in subsection 4.1.
Then we prepare some approximation results of functions on Wasserstein space in subsection 4.2.

\subsection{Poisson equation and regularities of averaged functions}
Let us consider the following parameterized SDE in $\mR^{d_2}$:
\begin{align} \label{ytx}
\dif Y^{t,x,\mu}_s =F(t,x,\mu,Y^{t,x,\mu}_s)\dif s+G(t,x,\mu,Y^{t,x,\mu}_s)\dif W_s,\qquad   Y^{t,x,\mu}_0=y\in\mR^{d_2},
\end{align}
where $(t,x,\mu)\in\mR_+\times\mR^{d_1}\times\sP_2(\mR^{d_1})$ are regarded as parameters. We  make the following assumption on the coefficients:

\vspace{0.1cm}
\begin{description}
\item[$(\mathbf{\tilde A}_{FG})$] the coefficient  $GG^*(t,x,\mu,y)$ is non-degenerate in the sense that there exist constants $k,\varrho>0$ such that for every $(t,x,\mu,y)\in \mR_+\times\mR^{d_1}\times\sP_2(\mR^{d_1})\times \mR^{d_2}$,
	$$
	\varrho(1+|y|)^{-k}|z|^{2} \leq |G^*(t,x,\mu,y)z|^2, \quad \forall z\in \mathbb{R}^{d_{2}},
	$$
and for every $q\geq2$, there exist constants $\Lambda_1, \Lambda_2>0$ such that for every $(t,x,\mu,y)\in\mR_+\times\mR^{d_1}\times\sP_2(\mR^{d_1})\times \mR^{d_2}$,
\begin{align*}
2\<y,F(t,x,\mu,y)\>+(q-1)\|G(t,x,\mu,y)\|^2\leq -\Lambda_1|y|^2+\Lambda_2.
\end{align*}
\end{description}

Given a function $f$ on $\mR_+\times\mR^{d_1}\times\sP_2(\mR^{d_1})\times\mR^{d_2}$, we consider the following Poisson equation in the whole space $\mR^{d_2}$:
\begin{align}\label{pof}
\sL_0(t,x,\mu,y)U(t,x,\mu,y)=-f(t,x,\mu,y),
\end{align}
where $(t,x,\mu)\in\mR_+\times\mR^{d_1}\times\sP_2(\mR^{d_1})$  are regarded as parameters, and  the operator $\sL_0$  is defined by
\begin{align}\label{ophat}
\sL_0\hat \varphi(y)&:=\sL_0(t,x,\mu,y)\hat \varphi(y)\no\\
&:=\frac{1}{2}\T\Big(\cG(t,x,\mu,y)\cdot\p_y^2\hat \varphi(y)\Big)
+F(t,x,\mu,y)\cdot\p_y\hat \varphi(y),
\end{align}
where $\cG(t,x,\mu,y)=GG^*(t,x,\mu,y)$. Note that $\sL_0$ is just the infinitesimal generator of $ Y^{t,x,\mu}_s$ given by (\ref{ytx}).
In order to ensure the well-posedness of the equation (\ref{pof}), we need to assume that $f$ satisfies the following centering condition:
\begin{align}\label{cenf}
\int_{\mR^{d_2}}f(t,x,\mu,y)\zeta^{t,x,\mu}(\dif y)=0,\quad\forall (t,x,\mu)\in\mR_+\times\mR^{d_1}\times\sP_2(\mR^{d_1}),
\end{align}
where $\zeta^{t,x,\mu}(dy)$ is the invariant measure of the equation (\ref{ytx}).
We have the following  result.

\bt\label{lem1}
Assume that $(\mathbf{\tilde A}_{FG})$ holds, and  $F,G\in C^{\a/2,\a,(\ell,\a),\b}_p$ with some $0<\a,\b\leq2$ and  $\ell=1,2$. Then for every function $f\in C_p^{\a/2,\a,(\ell,\a),\b}$ satisfying (\ref{cenf}), there exists a unique solution $U\in C_p^{\a/2,\a,(\ell,\a),2+\b}$ to (\ref{pof}) which also satisfies (\ref{cenf}) and is given by
\begin{align*}
U(t,x,\mu,y)=\int_0^\infty\mE f\Big(t,x,\mu,Y_s^{t,x,\mu}(y)\Big)\dif s,
\end{align*}
where $Y_s^{t,x,\mu}(y)$ satisfies the equation (\ref{ytx}).

Moreover, there exits a constant $p>0$ such that
for any $t\in\mR_+$, $x\in\mR^{d_1}$, $\mu\in\sP_2(\mR^{d_1})$ and $y\in\mR^{d_2}$,
\begin{align*}
&|U(t,x,\mu,y)|+|\p_yU(t,x,\mu,y)|+|\p_y^2U(t,x,\mu,y)|\leq C_0(1+|y|^p),\\
&\|U(\cdot,\cdot,\cdot,y)\|_{C_b^{\a/2,\a,(\ell,\a)}}\leq C_{0}(1+|y|^p),
\end{align*}
and for any $y_1,y_2\in\mR^{d_2}$,
\begin{align*}
|\p_y^2U(t,x,\mu,y_1)-\p_y^2U(t,x,\mu,y_2)|\leq C_0(1+|y_1|^p+|y_2|^p)|y_1-y_2|^\b,
\end{align*}
where $C_0>0$ is a constant depending on   $[F]_{C_p^{\a/2,\a,(\ell,\a),\b}}$, $[G]_{C_p^{\a/2,\a,(\ell,\a),\b}}$ and $[f]_{C_p^{\a/2,\a,(\ell,\a),\b}}$.
\et

\br
Theorem \ref{lem1} will be used to study the asymptotic behavior of the non-autonomous stochastic system (\ref{sde2}) in Section 5. In fact, we will only need this result when the coefficients of $\sL_0$ do not depend on the parameter $\mu$, i.e., SDE (\ref{ytx}) with
\begin{align*}
\dif Y^{t,x}_s =F(t,x,Y^{t,x}_s)\dif s+G(t,x,Y^{t,x}_s)\dif W_s,\qquad   Y^{t,x}_0=y\in\mR^{d_2}.
\end{align*}
The reason we consider (\ref{pof}) with an additional  parameter $\mu$ is that we need this to  prove the regularity of the averaged functions with respect to the $\mu$-variable, see Corollary \ref{cor} below, which will be further used to establish the regularity of the averaged functions of the form (\ref{barb}) in Lemma \ref{lem3} (this avoids to use the mean-field type Poisson equation associated with the non-linear stochastic system (\ref{frozen0})).
\er

\begin{proof}
The existence and uniqueness of the solution $U$ to (\ref{pof}), and the estimates with respect to the parameters $(t,x)$ and the variable $y$ can be proved similarly as in \cite[Theorem 2.1]{RX1}, so we omit the details here. Our task is to prove the estimate of $U$ with respect to $\mu$-variable. Let us first consider the case $\ell=1$. Since $U$ satisfies the Poisson equation (\ref{pof}), by the definition of the linear functional derivative,
we have for any $\mu,\mu'\in \sP_2(\mR^{d_1})$ and $\theta>0$ that
\begin{align*}
&\sL_0(t,x,\mu,y)\Big(\frac{U(t,x,\mu,y)-U(t,x,(1-\theta)\mu+\theta\mu',y)}{\theta}\Big)\\
&=\frac{f(t,x,(1-\theta)\mu+\theta\mu',y)-f(t,x,\mu,y)}{\theta}\\
&\quad+\frac{1}{\theta}\Big(\sL_0(t,x,(1-\theta)\mu+\theta\mu',y)-\sL_0(t,x,\mu,y)\Big)U(t,x,(1-\theta)\mu+\theta\mu',y)\\
&=:h_1^\theta(t,x,\mu,y),
\end{align*}
where
\begin{align*}
&\Big(\sL_0(t,x,(1-\theta)\mu+\theta\mu',y)-\sL_0(t,x,\mu,y)\Big)U(t,x,(1-\theta)\mu+\theta\mu',y)\\
&=\Big(F(t,x,(1-\theta)\mu+\theta\mu',y)-F(t,x,\mu,y)\Big) \cdot\p_yU(t,x,(1-\theta)\mu+\theta\mu',y)\\
&\quad+\frac{1}{2}\T\Big(\cG(t,x,(1-\theta)\mu+\theta\mu',y)-\cG(t,x,\mu,y)\Big) \cdot\p^2_yU(t,x,(1-\theta)\mu+\theta\mu',y).
\end{align*}
This together with $U\in C_p^{0,0,(0,0),2+\b}$ implies that the function $h_1^\theta(t,x,\mu,y)$ satisfies the centering condition, that is,
\begin{align}\label{u3}
\int_{\mR^{d_2}}h_1^\theta(t,x,\mu,y)\zeta^{t,x,\mu}(\dif y)=0,
\end{align}
and we further have
\begin{align}\label{u4}
\frac{U(t,x,(1-\theta)\mu+\theta\mu',y)-U(t,x,\mu,y)}{\theta}
=\int_0^\infty \mE h_1^\theta(t,x,\mu,Y_s^{t,x,\mu}(y))\dif s.
\end{align}
Note that
\begin{align*}
&\lim_{\theta\to0}h_1^\theta(t,x,\mu,y)\\
&=\int_{\mR^{d_1}}\Big[\frac{\d f}{\d\mu}(t,x,\mu,y)(\tilde x)+\frac{\d F}{\d\mu}(t,x,\mu,y)(\tilde x)\cdot\p_yU(t,x,\mu,y)\\
&\qquad\quad+\frac{1}{2}\T\Big(\frac{\d \cG}{\d\mu}(t,x,\mu,y)(\tilde x)\cdot\p^2_yU(t,x,\mu,y)\Big)\Big](\mu'-\mu)(\dif \tilde x)\no\\
&=:\int_{\mR^{d_1}}h_1(t,x,\mu,y)(\tilde x)(\mu'-\mu)(\dif \tilde x).
\end{align*}
Meanwhile, by the assumptions on $F,G,f$ and (\ref{lip}), there exists a constant $C_1>0$ such that
\begin{align*}
&|F(t,x,(1-\theta)\mu+\theta\mu',y)-F(t,x,\mu,y)|\leq C_1\theta(1+|y|^p)\cdot \|\mu'-\mu\|_{{\text{TV}}},\\
&|\cG(t,x,(1-\theta)\mu+\theta\mu',y)-\cG(t,x,\mu,y)|\leq C_1\theta(1+|y|^p)\cdot \|\mu'-\mu\|_{{\text{TV}}}
\end{align*}
and
\begin{align*}
|f(t,x,(1-\theta)\mu+\theta\mu',y)-f(t,x,\mu,y)|\leq C_1\theta(1+|y|^p)\cdot \|\mu'-\mu\|_{{\text{TV}}}.
\end{align*}
These together with (\ref{u3}) and the dominated convergence theorem yield that
\begin{align}\label{h1c}
\int_{\mR^{d_2}}h_1(t,x,\mu,y)(\tilde x)\zeta^{t,x,\mu}(\dif y)=0.
\end{align}
According to (\ref{u3}) and (\ref{erg}), there exist constants $C_2,\l_0>0$ such that
\begin{align*}
\big|\mE\big[h_1^\theta(t,x,\mu,Y_s^{t,x,\mu}(y))\big]\big|\leq C_2(1+|y|^p)e^{-\l_0s}.
\end{align*}
Thus, taking the limit $\theta\to0$ in (\ref{u4}) we obtain
$$
\frac{\d U}{\d\mu}(t,x,\mu,y)(\tilde x)=\int_0^\infty \mE h_1(t,x,\mu,Y_s^{t,x,\mu}(y))(\tilde x)\dif s.
$$
Furthermore, by the assumptions that $F,G\in C_p^{0,0,(1,\a),\b}$ and $f\in C_p^{0,0,(1,\a),\b}$, we deduce that for every $y_1,y_2\in\mR^{d_2}$,
\begin{align*}
&\big|h_1(t,x,\mu,y_1)(\tilde x)-h_1(t,x,\mu,y_2)(\tilde x)\big|\leq\Big|\frac{\d f}{\d\mu}(t,x,\mu,y_1)(\tilde x)-\frac{\d f}{\d\mu}(t,x,\mu,y_2)(\tilde x)\Big|\\
&\quad+\Big|\frac{\d F}{\d\mu}(t,x,\mu,y_1)(\tilde x)\cdot\p_yU(t,x,\mu,y_1)-\frac{\d F}{\d\mu}(t,x,\mu,y_2)(\tilde x)\cdot\p_yU(t,x,\mu,y_2)\Big|\\
&\quad+\frac{1}{2}\Big|\T\Big(\frac{\d \cG}{\d\mu}(t,x,\mu,y_1)(\tilde x)\cdot\p^2_yU(t,x,\mu,y_1)-\frac{\d \cG}{\d\mu}(t,x,\mu,y_2)(\tilde x)\cdot\p^2_yU(t,x,\mu,y_2)\Big)\Big|\\
&\leq C_3(1+|y_1|^p+|y_2|^p)|y_1-y_2|^\b,
\end{align*}
where $C_3>0$ is a constant depending on $[F]_{C_p^{0,0,(1,\a),\b}}$, $[G]_{C_p^{0,0,(1,\a),\b}}$ and $[f]_{C_p^{0,0,(1,\a),\b}}$. This means that $h_1$ is $\b$-H\"older continuous with respect to $y$-variable. As a result, we get that
\begin{align*}
\sL_0(t,x,\mu,y)\frac{\d U}{\d\mu}(t,x,\mu,y)(\tilde x)=-h_1(t,x,\mu,y)(\tilde x),
\end{align*}
which in turn implies that
\begin{align}\label{estu}
\Big|\frac{\d U}{\d\mu}(t,x,\mu,y)(\tilde x)\Big|\leq C_3(1+|y|^p).
\end{align}
Moreover, we have that for every $\tilde x_1,\tilde x_2\in\mR^{d_1}$,
\begin{align*}
&\sL_0(t,x,\mu,y)\Big(\frac{\d U}{\d\mu}(t,x,\mu,y)(\tilde x_1)-\frac{\d U}{\d\mu}(t,x,\mu,y)(\tilde x_2)\Big)\\
&=-\big[h_1(t,x,\mu,y)(\tilde x_1)-h_1(t,x,\mu,y)(\tilde x_2)\big].
\end{align*}
Consequently, we have
\begin{align*}
\Big|\frac{\d U}{\d\mu}(t,x,\mu,y)(\tilde x_1)-\frac{\d U}{\d\mu}(t,x,\mu,y)(\tilde x_2)\Big|\leq C_4(1+|y|^p)|\tilde x_1-\tilde x_2|^\a,
\end{align*}
where $C_4>0$ is a constant depending on $[F]_{C_p^{0,0,(1,\a),\b}}$, $[G]_{C_p^{0,0,(1,\a),\b}}$ and $[f]_{C_p^{0,0,(1,\a),\b}}$. This together with (\ref{estu}) means that $U(t,x,\cdot,y)\in C_b^{(1,\a)}$. Similarly, we can deduce that
$$
\frac{\d^2U}{\d\mu^2}(t,x,\mu,y)(\tilde x_1,\tilde x_2)=\int_0^\infty \mE h_2(t,x,\mu,Y_s^{t,x,\mu}(y))(\tilde x_1,\tilde x_2)\dif s,
$$
where $h_2$ satisfies the centering condition (\ref{cenf}) and is given by
\begin{align*}
h_2(t,x,\mu,y)(\tilde x_1,\tilde x_2)&=\frac{\d^2f}{\d\mu^2}(t,x,\mu,y)(\tilde x_1,\tilde x_2)+\frac{\d^2F}{\d\mu^2}(t,x,\mu,y)(\tilde x_1,\tilde x_2)\cdot\p_yU(t,x,\mu,y)\\
&\quad+\frac{\d F}{\d\mu}(t,x,\mu,y)(\tilde x_1)\cdot\p_y\Big[\frac{\d U}{\d\mu}(t,x,\mu,y)(\tilde x_2)\Big]\\
&\quad+\frac{\d F}{\d\mu}(t,x,\mu,y)(\tilde x_2)\cdot\p_y\Big[\frac{\d U}{\d\mu}(t,x,\mu,y)(\tilde x_1)\Big]\\
&\quad+\frac{1}{2}\T\Big(\frac{\d^2\cG}{\d\mu^2}(t,x,\mu,y)(\tilde x_1,\tilde x_2)\cdot\p^2_yU(t,x,\mu,y)\Big)\\
&\quad+\frac{1}{2}\T\Big(\frac{\d \cG}{\d\mu}(t,x,\mu,y)(\tilde x_1)\cdot\p^2_y\Big[\frac{\d U}{\d\mu}(t,x,\mu,y)(\tilde x_2)\Big)\Big]\\
&\quad+\frac{1}{2}\T\Big(\frac{\d \cG}{\d\mu}(t,x,\mu,y)(\tilde x_2)\cdot\p^2_y\Big[\frac{\d U}{\d\mu}(t,x,\mu,y)(\tilde x_1)\Big)\Big].
\end{align*}
Using the similar arguments as above, we get $U(t,x,\cdot,y)\in C_b^{(2,\a)}$. Thus the proof is completed.
\end{proof}

Given a function $f(t,x,\mu,y)$, we shall denote
\begin{align}\label{barf}
\bar f(t,x,\mu):=\int_{\mR^{d_2}}f(t,x,\mu,y)\zeta^{t,x,\mu}(\dif y).
\end{align}
As a direct application of Theorem \ref{lem1}, we have the following regularity results for the averaged function $\bar f$.

\bc\label{cor}
Assume that $(\mathbf{\tilde A}_{FG})$ holds, $F, G\in C^{\a/2,\a,(2,\a),\b}_p$ and $f\in C^{\a/2,\a,(2,\a),\b}_p$ with $0<\a,\b\leq2$. Let $\bar f(t,x,\mu)$ be defined by (\ref{barf}). Then we have $\bar f\in C_b^{\a/2,\a,(2,\a)}$. In particular, we have
\begin{align*}
\frac{\d\bar f}{\d\mu}(t,x,\mu)(\tilde x)&=\int_{\mR^{d_2}}\bigg[\frac{\d f}{\d\mu}(t,x,\mu,y)(\tilde x)+\frac{\d F}{\d\mu}(t,x,\mu,y)(\tilde x)\cdot\p_y\Phi(t,x,\mu,y)\\
&\qquad\quad+\frac{1}{2}\T\Big(\frac{\d \cG}{\d\mu}(t,x,\mu,y)(\tilde x)\cdot\p^2_y\Phi(t,x,\mu,y)\Big)\bigg]\zeta^{t,x,\mu}(\dif y),
\end{align*}
where $\Phi$ is the solution of the following Poisson equation:
\begin{align}\label{bposs}
\sL_0(t,x,\mu,y)\Phi(t,x,\mu,y)=-[f(t,x,\mu,y)-\bar f(t,x,\mu)].
\end{align}
\ec
\begin{proof}
The assertion that the function  $\bar f(\cdot,\cdot,\mu)\in C_b^{\a/2,\a}$ follows by \cite[Lemma 3.2]{RX1}. Here, we need only show the regularity of $\bar f$ with respect to the $\mu$-variable.
Note that the function
$$
(t,x,\mu,y)\mapsto f(t,x,\mu,y)-\bar f(t,x,\mu)=:\Delta f(t,x,\mu,y)
$$
always satisfies the centering condition. Thus under our assumptions and by Theorem \ref{lem1} there exists a unique solution $\Phi\in C_p^{0,0,(0,0),2+\b}$ to the Poisson equation (\ref{bposs}).
Following the same arguments as in (\ref{h1c}) we have
\begin{align}\label{bc}
&\int_{\mR^{d_2}}\bigg[\frac{\delta\Delta f}{\d\mu}(t,x,\mu,y)(\tilde x)+\frac{\d F}{\d\mu}(t,x,\mu,y)(\tilde x)\cdot\p_y\Phi(t,x,\mu,y)\no\\
&\qquad+\frac{1}{2}\T\Big(\frac{\d \cG}{\d\mu}(t,x,\mu,y)(\tilde x)\cdot\p^2_y\Phi(t,x,\mu,y)\Big)\bigg]\zeta^{t,x,\mu}(\dif y)=0.
\end{align}
Since
\begin{align*}
&\int_{\mR^{d_2}}\frac{\delta \Delta f}{\d\mu}(t,x,\mu,y)(\tilde x)\zeta^{t,x,\mu}(\dif y)\\
&=- \frac{\d\bar f}{\d\mu}(t,x,\mu)(\tilde x)+\int_{\mR^{d_2}} \frac{\d f}{\d\mu}(t,x,\mu,y)(\tilde x)\zeta^{t,x,\mu}(\dif y),
\end{align*}
by (\ref{bc}) we get
\begin{align*}
\frac{\d\bar f}{\d\mu}(t,x,\mu)(\tilde x)&=\int_{\mR^{d_2}}\bigg[\frac{\d f}{\d\mu}(t,x,\mu,y)(\tilde x)+\frac{\d F}{\d\mu}(t,x,\mu,y)(\tilde x)\cdot\p_y\Phi(t,x,\mu,y)\\
&\qquad\quad+\frac{1}{2}\T\Big(\frac{\d \cG}{\d\mu}(t,x,\mu,y)(\tilde x)\cdot\p^2_y\Phi(t,x,\mu,y)\Big)\bigg]\zeta^{t,x,\mu}(\dif y).
\end{align*}
This together with the assumption $F,G,f\in C_p^{0,0,(1,\a),\b}$ implies that $[\d\bar f/\d\mu](t,x,\mu)(\cdot)$ is $\a$-H\"older continuous. As a result, we have $\bar f(t,x,\cdot)\in C_b^{(1,\a)}$. Similarly, we have
\begin{align*}
\frac{\d^2\bar f}{\d\mu^2}(t,x,\mu)(\tilde x_1,\tilde x_2)&=\int_{\mR^{d_2}}\bigg[\frac{\d^2 f}{\d\mu^2}(t,x,\mu,y)(\tilde x_1,\tilde x_2)+\frac{\d^2 F}{\d\mu^2}(t,x,\mu,y)(\tilde x_1,\tilde x_2)\cdot\p_y\Phi(t,x,\mu,y)\\
&\qquad\quad+\frac{\d F}{\d\mu}(t,x,\mu,y)(\tilde x_1)\cdot\p_y\Big[\frac{\d\Phi}{\d\mu}(t,x,\mu,y)(\tilde x_2)\Big]\\
&\qquad\quad+\frac{\d F}{\d\mu}(t,x,\mu,y)(\tilde x_2)\cdot\p_y\Big[\frac{\d\Phi}{\d\mu}(t,x,\mu,y)(\tilde x_1)\Big]\\
&\qquad\quad+\frac{1}{2}\T\Big(\frac{\d^2 \cG}{\d\mu^2}(t,x,\mu,y)(\tilde x_1,\tilde x_2)\cdot\p^2_y\Phi(t,x,\mu,y)\Big)\\
&\qquad\quad+\frac{1}{2}\T\Big(\frac{\d \cG}{\d\mu}(t,x,\mu,y)(\tilde x_2)\cdot\p^2_y\Big[\frac{\d\Phi}{\d\mu}(t,x,\mu,y)(\tilde x_2)\Big]\Big)\\
&\qquad\quad+\frac{1}{2}\T\Big(\frac{\d \cG}{\d\mu}(t,x,\mu,y)(\tilde x_2)\cdot\p^2_y\Big[\frac{\d\Phi}{\d\mu}(t,x,\mu,y)(\tilde x_1)\Big]\Big)\bigg]\zeta^{t,x,\mu}(\dif y),
\end{align*}
which implies that $[\d^2\bar f/\d\mu^2](t,x,\mu)(\tilde x_1,\cdot)$ is $\a$-H\"older continuous and $\bar f(t,x,\cdot)\in C_b^{(2,\a)}$. Thus the proof is finished.
\end{proof}

\subsection{Mollifying approximation on Wasserstein space}
Due to the low regularity assumptions on the coefficients, we need some mollification arguments for both the space and the distribution variables. The  mollification for the space variable is classical, the main aim here is to construct smooth approximations of   functions $f: \sP_2(\mR^{d_1})\to\mR$. In particular, we show that when $f\in C_b^{(1,\a)}$ with some $0<\alpha\leq 1$ (which is only $\a$-H\"older continuous with respect to the Wasserstein distance and thus not Lions differentiable), there exists a sequence of functions $f_n\in C^{(1,\infty)}$ (which admit $1$-order Lions derivative) such that
$$
\|f_n-f\|_\infty\leq C_0\,n^{-\a},\quad\|\p_\mu f_n\|_\infty\leq C_0\,n^{1-\a}\quad\text{and}\quad\|\p_{\tilde x}\p_\mu f_n\|_\infty\leq C_0\,n^{2-\a},
$$
where $C_0>0$ is a constant independent of $n$. The explicit dependence on $n$ on the right hand sides of the above inequalities will play an important role in determining the rate of convergence of the multi-scales system (\ref{sde1}).

Let $\rho_1:\mR\to[0,1]$ and $\rho_2:\mR^{d_1}\to[0,1]$ be two smooth radial convolution kernel functions such that $\int_{\mR}\rho_1(r)\dif r=\int_{\mR^{d_1}}\rho_2(x)\dif x=1$, and for any $k\geq1$, $|\nabla^k\rho_1(r)|\leq C_k\rho_1(r)$ and $|\nabla^k\rho_2(x)|\leq C_k\rho_2(x)$ where $C_k$ are positive constants. For every $n\geq1$, let
$$\rho^n_1(r):=n^2\rho_1(n^2r)\qquad\mathrm{and}\qquad\rho^n_2(x):=n^{d_1}\rho_2(nx).$$
Given a function $f(t,x,\mu,y)$, we define the mollifying approximations of $f$ in $t,x$ and $\mu$ variables by
\begin{align}\label{fnn}
f_n(t,x,\mu,y)&:=f(\cdot,\cdot,\mu\ast\rho_2^n,y)\ast\rho_2^n\ast\rho_1^n\no\\
&:=\int_{\mR^{d_1+1}}f(t-s,x-z,\mu\ast\rho_2^n,y)\rho_2^n(z)\rho_1^n(s)\dif z\dif s,
\end{align}
where $\mu\ast\rho_2^n$ is defined by
\begin{align}\label{mun}
\mu\ast\rho_2^n(\cdot):=\int_{\cdot}\int_{\mR^{d_1}}\rho_2^n(x-z)\mu(\dif z)\dif x.
\end{align}
In particular, when $f$ depends only on the $\mu$-variable, we have
\begin{align*}
f_n(\mu)&:=f(\mu\ast\rho_2^n),
\end{align*}
which gives the mollifying approximation for the functions defined on Wasserstein space.
We have the following results.

\bl\label{fn}
Assume that $f\in C_p^{\a/2,\a,(1,\a),0}$ with $0<\a\leq2$ and let $f_n$ be defined by (\ref{fnn}). Then we have $f_n\in C_p^{1,2,(1,2),0}$, and
\begin{align}
&\|f_n(\cdot,\cdot,\cdot,y)-f(\cdot,\cdot,\cdot,y)\|_{\infty}\leq C_0\,(1+|y|^p)n^{-\a} ,\label{fn1}\\
&\|\p_xf_n(\cdot,\cdot,\cdot,y)\|_\infty+\|\p_\mu f_n(\cdot,\cdot,\cdot,y)(\cdot)\|_{\infty}\leq C_0\,(1+|y|^p)n^{1-\a\wedge1} ,\label{fn2}
\end{align}
and
\begin{align}\label{fn3}
\|\p_tf_n(\cdot,\cdot,\cdot,y)\|_{\infty}\!+\!\|\p^2_xf_n(\cdot,\cdot,\cdot,y)\|_{\infty}\!+\!\|\p_{\tilde x}\p_\mu f_n(\cdot,\cdot,\cdot,y)(\cdot)\|_{\infty}\leq C_0\,(1+|y|^p)n^{2-\a},
\end{align}
where $C_0>0$ is a constant independent of $n$.
\el
\begin{proof}
The estimates concerning the derivatives of $f_n$ with respect to $t$ and $x$ variables can be proved similarly as in \cite[Lemma 4.1]{RX1}, we omit the details here. In the following, we focus on the estimates of $f_n$ with respect to $\mu$ when $0<\a\leq 1$, the case $1<\a\leq 2$ can be proved similarly.
By the definition we have
\begin{align*}
&|f_n(t,x,\mu,y)-f(t,x,\mu,y)|\\
&\leq\int_{\mR^{d_1+1}}|f(t-s,x-z,\mu\ast\rho_2^n,y)-f(t,x,\mu,y)|\rho_2^n(z)\rho_1^n(s)\dif z\dif s\\
&\leq\int_{\mR^{d_1+1}}|f(t-s,x-z,\mu\ast\rho_2^n,y)-f(t,x,\mu\ast\rho_2^n,y)|\rho_2^n(z)\rho_1^n(s)\dif z\dif s\\
&\quad+|f(t,x,\mu\ast\rho_2^n,y)-f(t,x,\mu,y)|=:\cI_1+\cI_2.
\end{align*}
For $\cI_1$, by the assumption that $f\in C_p^{\a/2,\a,(1,\a),0}$, there exists a constant $p>0$ such that
\begin{align*}
\cI_1&\leq C_1\int_{\mR^{d_1+1}}(|s|^{\a/2}+|z|^\a)\cdot(1+|y|^p)\rho_2^n(z)\rho_1^n(s)\dif z\dif s\\
&\leq C_1\,n^{-\a}(1+|y|^p).
\end{align*}
As for $\cI_2$, by (\ref{mun}) and (\ref{lin}), we have
\begin{align*}
\cI_2&=\Big|\int_0^1\int_{\mR^{d_1}}\frac{\d f}{\d\mu}(t,x,(1-\theta)\mu+\theta\mu\ast\rho_2^n,y)(\tilde x)(\mu\ast\rho_2^n-\mu)(\dif \tilde x)\dif \theta\Big|\\
&=\Big|\int_0^1\int_{\mR^{d_1}}\int_{\mR^{d_1}}\frac{\d f}{\d\mu}(t,x,(1-\theta)\mu+\theta\mu\ast\rho_2^n,y)(\tilde x)\rho_2^n(z-\tilde x)\mu(\dif z)\dif \tilde x\dif \theta\\
&\quad-\int_0^1\int_{\mR^{d_1}}\frac{\d f}{\d\mu}(t,x,(1-\theta)\mu+\theta\mu\ast\rho_2^n,y)(z)\mu(\dif z)\dif \theta\Big|\\
&\leq\int_0^1\int_{\mR^{d_1}}\int_{\mR^{d_1}}\Big|\frac{\d f}{\d\mu}(t,x,(1-\theta)\mu+\theta\mu\ast\rho_2^n,y)(z-\tilde x)\\
&\qquad\quad-\frac{\d f}{\d\mu}(t,x,(1-\theta)\mu+\l\mu\ast\rho_2^n,y)(z)\Big|\rho_2^n(\tilde x)\dif \tilde x\mu(\dif z)\dif \theta\\
&\leq C_2\int_{\mR^{d_1}}|\tilde x|^\a\cdot(1+|y|^p)\rho_2^n(\tilde x)\dif \tilde x\leq C_2\,n^{-\a}(1+|y|^p).
\end{align*}
Combining the above computations, (\ref{fn1}) is true. Furthermore, by the chain rule formula in Lemma \ref{chain} we deduce that
\begin{align*}
\frac{\d f_n}{\d\mu}(t,x,\mu,y)(\tilde x)&=\int_{\mR^{d_1}}\int_{\mR^{d_1+1}}\frac{\d f}{\d\mu}(t-s,x-z,\mu\ast\rho_2^n,y)(\tilde z)\rho_2^n(z)\rho_1^n(s)\dif z\dif s\cdot\rho_2^n(\tilde x-\tilde z)\dif\tilde z\\
&=\int_{\mR^{d_1+1}}\int_{\mR^{d_1}}\frac{\d f}{\d\mu}(t-s,x-z,\mu\ast\rho_2^n,y)(\tilde z)\rho_2^n(\tilde x-\tilde z)\rho_2^n(z)\rho_1^n(s)\dif\tilde z\dif z\dif s,
\end{align*}
which in turn implies that
\begin{align*}
&|\p_{\tilde x}\p_\mu f_n(t,x,\mu,y)(\tilde x)|=\Big|\p_{\tilde x}^2\frac{\d f_n}{\d\mu}(t,x,\mu,y)(\tilde x)\Big|\\
&=\Big|\int_{\mR^{d_1+1}}\int_{\mR^{d_1}}\frac{\d f}{\d\mu}(t-s,x-z,\mu\ast\rho_2^n,y)(\tilde z)\cdot\nabla^2\rho_2^n(\tilde x-\tilde z)\rho_2^n(z)\rho_1^n(s)\dif\tilde z\dif z\dif s\Big|\\
&\leq\int_{\mR^{d_1+1}}\int_{\mR^{d_1}}\Big|\frac{\d f}{\d\mu}(t-s,x-z,\mu\ast\rho_2^n,y)(\tilde x-\tilde z)-\frac{\d f}{\d\mu}(t-s,x-z,\mu\ast\rho_2^n,y)(\tilde x)\Big|\\
&\qquad\quad\times|\nabla^2\rho_2^n(\tilde z)|\dif\tilde z\cdot\rho_2^n(z)\rho_1^n(s)\dif z\dif s\\
&\leq C_3\int_{\mR^{d_1}}|\tilde z|^\a\cdot(1+|y|^p)\cdot|\nabla^2\rho_2^n(\tilde z)|\dif \tilde z\leq C_3\,n^{2-\a}(1+|y|^p).
\end{align*}
Similarly, we get
\begin{align*}
|\p_\mu f_n(t,x,\mu,y)(\tilde x)|=\Big|\p_{\tilde x}\frac{\d f_n}{\d\mu}(t,x,\mu,y)(\tilde x)\Big|\leq C_4\,n^{1-\a}(1+|y|^p).
\end{align*}
Thus, the estimates (\ref{fn2}) and (\ref{fn3}) hold. The proof is finished.
\end{proof}

\section{Multi-scale non-autonomous SDEs}

In this section, we consider the non-autonomous multi-scale stochastic system (\ref{sde2}). Namely,
\begin{equation} \label{sde22}
\left\{ \begin{aligned}
&\dif \hat X^{\eps}_t =b_\eps(t,\hat X^{\eps}_t,\hat Y^{\eps}_t)\dif t
+\sigma_\eps(t,\hat X^{\eps}_t)\dif W^1_t,\qquad\qquad\qquad\, \hat X^{\eps}_0=\xi,\\
&\dif \hat Y^{\eps}_t =\frac{1}{\eps}F_\eps(t,\hat X^{\eps}_t,\hat Y^{\eps}_t)\dif t+\frac{1}{\sqrt{\eps}}G_\eps(t,\hat X^{\eps}_t,\hat Y^{\eps}_t)\dif W_t^2,\qquad   \hat Y^{\eps}_0=\eta,
\end{aligned} \right.
\end{equation}
where for every $0<\eps\ll1$, $b_\eps: \mR_+\times\mR^{d_1}\times\mR^{d_2}\to\mR^{d_1}$, $\sigma_\eps: \mR_+\times\mR^{d_1}\to\mR^{d_1}\otimes\mR^{d_1}$, $F_\eps: \mR_+\times\mR^{d_1}\times\mR^{d_2}\to\mR^{d_2}$ and $G_\eps: \mR_+\times\mR^{d_1}\times\mR^{d_2}\to\mR^{d_2}\otimes\mR^{d_2}$ are measurable functions.
The aim is to prove the strong  convergence in the averaging principle as well as  the convergence of the distribution for both the slow process $\hat X^{\eps}_t$ and the fast motion $\hat Y^{\eps}_t$ of the system (\ref{sde22}).  These will play an important role below to study the asymptotic limit for the non-linear stochastic system (\ref{sde1}).

To study the asymptotic behavior of the system (\ref{sde22}), our basic assumptions on the coefficients are as follows.

\vspace{2mm}
\begin{description}
\item[\bf (A)] There exist functions $\hat b(t,x,y)$, $\hat \sigma(t,x)$, $\hat F(t,x,y)$ and $\hat G(t,x,y)$  such that for every fixed $t\in\mR_+$, $x\in\mR^{d_1}$ and $y\in\mR^{d_2}$,
\begin{align*}
&\qquad\lim_{\eps\to0}b_\eps(t,x,y)=\hat b(t,x,y),\qquad\qquad \lim_{\eps\to0}\sigma_\eps(t,x)=\hat \sigma(t,x),\\
&\qquad\lim_{\eps\to0}F_\eps(t,x,y)=\hat F(t,x,y)\quad{\rm{and}}\quad \lim_{\eps\to0}G_\eps(t,x,y)=\hat G(t,x,y).
\end{align*}
Moreover, the function $\hat a=\hat\sigma\hat\sigma^*$ and $\hat\cG=\hat G\hat G^*$ are  non-degenerate in the sense that there exists  constants $k, \varrho>0$ such that for any $t>0$, $x, z_1\in\mR^{d_1}$ and $y, z_2\in\mR^{d_2}$,
\begin{align*}
\qquad\varrho|z_1|^2\leq |\hat\sigma^*(t,x)z_1|^2\quad \text{and}\quad \varrho(1+|y|)^{-k}|z_2|^2\leq |\hat G^*(t,x,y)z_2|^2,
\end{align*}
and for any $q\geq2$, there exist constants $\Lambda_1, \Lambda_2>0$ such that for any $\eps\in(0,1)$,
\begin{align}\label{hatF}
2\<y,F_\eps(t,x,y)\>+(q-1)\|G_\eps(t,x,y)\|^2\leq -\Lambda_1|y|^2+\Lambda_2.
\end{align}
\end{description}

\vspace{2mm}
We shall show that as $\eps\to0$, the slow component $\hat X_t^\eps$ in system (\ref{sde22}) will converge (both strongly and in distribution) to $\bar{\hat X}_t$  which satisfies the following averaged equation:
\begin{align}\label{ave1}
\dif \bar{\hat X}_t=\bar{\hat b}(t,\bar{\hat X}_t)\dif t+\hat \sigma(t,\bar{\hat X}_t)\dif W^1_t,\qquad  \bar{\hat X}_0=\xi,
\end{align}
where the averaged drift coefficient is defined by
\begin{align}\label{barhatb}
\bar{\hat b}(t,x):=\int_{\mR^{d_2}}\hat b(t,x,y)\hat\zeta^{t,x}(\dif y),
\end{align}
and $\hat\zeta^{t,x}(\dif y)$ is the unique invariant measure of the following frozen equation: for fixed $t\in\mR_+$ and $x\in\mR^{d_1}$,
\begin{align}\label{frozen1}
\dif \hat Y_s^{t,x}=\hat F(t,x,\hat Y_s^{t,x})\dif s+\hat G(t,x,\hat Y_s^{t,x})\dif \hat W_s,\quad \hat Y_0^{t,x}=\eta,
\end{align}
and $\hat W_s $ is a new standard Brownian motion. Moreover, we give explicit characterization for the limit of the  fast motion $\hat Y_t^\eps$, i.e., the distribution of $\hat Y_t^\eps$ will converge to
$\mE\,\hat\zeta^{t,\bar{\hat X}_t}$, where the expectation is taken with respect to $\bar{\hat X}_t$.
At the same time, we obtain rates of convergence depending explicitly on the convergence of $b_\eps$, $\sigma_\eps$, $F_\eps$ and $G_\eps$ to $\hat b$, $\hat \sigma$, $\hat F$ and $\hat G$.

To shorten the notation, for every $\mu\in\sP_2(\mR^{d_1})$ we denote
\begin{align}\label{zetamu}
\tilde{\hat\zeta}^{t,\mu}(\dif y):=\int_{\mR^{d_1}}\hat\zeta^{t,x}(\dif y)\mu(\dif x).
\end{align}
Fix $T>0$. The following is the main result of this section.

\bt\label{non-aut}
Let $(\hat X_t^\eps,\hat Y_t^\eps)$ and $\bar{\hat X}_t$  satisfy the equation (\ref{sde22}) and (\ref{ave1}), respectively. Assume that {\bf(A)}   holds, $\hat b,\hat F,\hat G\in C^{\a/2,\a,\b}_p$ and $\hat\sigma\in C^{\a/2,\a}_b$ with some $0<\a,\b\leq 2$. Then for every $t\in[0,T]$, we have:

\vspace{1mm}
\noindent(i) (strong convergence of $\hat X_t^\eps$) assume further that $\hat\sigma\in C^{\a/2,1}_b$, then
\begin{align}\label{st-non}
\mE\big|\hat X_t^\eps-\bar{\hat X}_t\big|^2&\leq C_{1}\,\bigg(\eps^{\a\wedge1}+\int_0^t\Big[\|b_\eps(s,\cdot,\cdot)-\hat b(s,\cdot,\cdot)\|^2_{L_p^\infty}+\|\sigma_\eps(s,\cdot)-\hat \sigma(s,\cdot)\|^2_{L^\infty}\no\\
&\quad+\|F_\eps(s,\cdot,\cdot)-\hat F(s,\cdot,\cdot)\|^2_{L_p^\infty}+\|G_\eps(s,\cdot,\cdot)-\hat G(s,\cdot,\cdot)\|^2_{L_p^\infty}\Big]\dif s\bigg);
\end{align}

\vspace{1mm}
\noindent(ii) (convergence of the distribution of $\hat X_t^\eps$) for every $\varphi\in C_b^{(2,\alpha)}(\sP_2(\mR^{d_1}))$,
\begin{align}\label{we-non-x}
\Big|\varphi\big(\cL_{\hat X_t^\eps}\big)-\varphi\big(\cL_{\bar{\hat X}_t}\big)\Big|&\leq C_{2}\,\bigg(\eps^{\frac{\a}{2}}+\int_0^t\,(t-s)^{\frac{\alpha}{2}-1}\cdot\Big[\|b_\eps(s,\cdot,\cdot)-\hat b(s,\cdot,\cdot)\|_{L_p^\infty}\no\\
&\quad+\|\sigma_\eps(s,\cdot)-\hat \sigma(s,\cdot)\|_{L^\infty}+\|F_\eps(s,\cdot,\cdot)-\hat F(s,\cdot,\cdot)\|_{L_p^\infty}\no\\
&\quad+\|G_\eps(s,\cdot,\cdot)-\hat G(s,\cdot,\cdot)\|_{L_p^\infty}\Big]\dif s\bigg);
\end{align}

\vspace{1mm}
\noindent(iii) (convergence of the distribution of $\hat Y_t^\eps$)
for every $\psi\in C_p^{(2,\beta)}(\sP_2(\mR^{d_2}))$,
\begin{align}\label{we-non-y}
\big|\psi(\cL_{\hat Y_t^\eps})-\psi(\tilde{\hat \zeta}^{t,\cL_{\bar{\hat{X}}_t}})\big|&\leq  C_{2}\,\bigg(\eps^{\frac{\a}{2}}+\int_0^t\,(t-s)^{\frac{\alpha}{2}-1}\cdot\Big[\|b_\eps(s,\cdot,\cdot)-\hat b(s,\cdot,\cdot)\|_{L_p^\infty}\no\\
&\quad+\|\sigma_\eps(s,\cdot)-\hat \sigma(s,\cdot)\|_{L^\infty}+\|F_\eps(s,\cdot,\cdot)-\hat F(s,\cdot,\cdot)\|_{L_p^\infty}\no\\
&\quad+\|G_\eps(s,\cdot,\cdot)-\hat G(s,\cdot,\cdot)\|_{L_p^\infty}\Big]\dif s\bigg)+C_3\bigg(e^{-\frac{\gamma t}{\eps}}\no\\
&\quad+\frac{1}{\eps}\int_0^t\Big(\frac{t-s}{\eps}\Big)^{\frac{\b}{2}-1}\cdot e^{-\frac{\gamma(t-s)}{\eps}}\Big[\|F_\eps(s,\cdot,\cdot)-\hat F(s,\cdot,\cdot)\|_{L_p^\infty}\no\\
&\quad+\|G_\eps(s,\cdot,\cdot)-\hat G(s,\cdot,\cdot)\|_{L_p^\infty}\Big]\dif s\bigg),
\end{align}
where $\tilde{\hat \zeta}^{t,\mu}$ is defined by (\ref{zetamu}), $C_{1}, C_{2}>0$ are constants depending on $T$, the norms of the coefficients and the constants in condition {\bf(A)}, and  $C_3, \gamma>0$ are constants independent of $T$.
\et

We shall prove the strong and weak convergence results in Theorem \ref{non-aut} separately in the following three subsections. Here, we provide the following important comments for the above results.
\br\label{resde}
(i) The explicit dependence  on the convergence of $b_\eps$, $\sigma_\eps$, $F_\eps$ and $G_\eps$ to $\hat b$, $\hat \sigma$, $\hat F$ and $\hat G$ in the estimates (\ref{st-non})-(\ref{we-non-y}) will play an important role below in studying the asymptotic limit for the non-autonomous approximation systems (\ref{sde4}). The time-singular terms $(t-s)^{\alpha/2-1}$ and $(t-s)^{\beta/2-1}$ in the estimates (\ref{we-non-x}) and (\ref{we-non-y}) are due to the low regularity of the test functions.

 \vspace{1mm}
(ii) The estimate (\ref{we-non-x}) for the convergence of the distribution of the slow process is more general than the weak convergence in the averaging principle of the classical multi-scale SDEs, i.e., when $b_\eps=\hat b$, $\sigma_\eps=\hat \sigma$, $F_\eps=\hat F$ and $G_\eps=\hat G$. In this case, estimate (\ref{we-non-x}) reduces to: for every $\varphi\in C_b^{(2,\alpha)}(\sP_2(\mR^{d_1}))$ and $t\in[0,T]$,
\begin{align}\label{sdex}
\big|\varphi\big(\cL_{\hat X_t^\eps}\big)-\varphi\big(\cL_{\bar{\hat X}_t}\big)\big|\leq C_{2,T}\,\eps^{\frac{\a}{2}}.
\end{align}
The novelty is the low regularity assumption that $\varphi\in C_b^{(2,\alpha)}(\sP_2(\mR^{d_1}))$, which is not  Lipschitz continuous with respect to the Wasserstein distance and thus not differentiable in the sense of Lions. This is one of the key points that why we could only require the $C_b^{(2,\alpha)}$-regularity on the coefficients of the non-linear system (\ref{sde1}) with respect to the distribution of the slow process. As a direct result of estimate (\ref{sdex}), we have that for every $\hat\varphi\in C_b^\alpha(\mR^{d_1})$  and $t\in[0,T]$,
\begin{align*}
\big|\mE\hat\varphi(\hat X_t^\eps)-\mE\hat\varphi(\bar{\hat{X}}_t)\big|\leq C_{2,T}\,\eps^{\frac{\a}{2}}.
\end{align*}
Existing results in the literature require at least $\hat\varphi\in C_b^{2+\alpha}(\mR^{d_1})$ in the above estimate, see e.g. \cite{RX2}.

  \vspace{1mm}
(iii) The estimate (\ref{we-non-y}) is, to the best of our knowledge,  the first result established for the convergence of the distribution of the fast motion even for classical multi-scale SDEs, i.e., we have for every $\psi\in C_p^{(2,\beta)}(\sP_2(\mR^{d_2}))$  and $t\in[0,T]$,
\begin{align*}
\big|\psi(\cL_{\hat Y_t^\eps})-\psi(\tilde{\hat \zeta}^{t,\cL_{\bar{\hat{X}}_t}})\big|\leq \tilde C_{2,T}\,\eps^{\frac{\a}{2}}+C_3\,e^{-\frac{\gamma t}{\eps}}.
\end{align*}
In particular, for every $\hat\psi\in C_p^\beta(\mR^{d_2})$, we have
\begin{align*}
& \bigg|\mE \hat\psi(\hat Y_t^\eps)-\int_{\mR^{d_1}}\!\!\int_{\mR^{d_2}}\hat\psi(y)\hat\zeta^{t,x}(\dif y)\cL_{\bar{\hat{X}}_t}(\dif x)\bigg]\bigg|\\
&=
\bigg|\mE \hat\psi(\hat Y_t^\eps)-\mE\bigg[\int_{\mR^{d_2}}\hat\psi(y)\hat\zeta^{t,\bar{\hat{X}}_t}(\dif y)\bigg]\bigg|\leq \tilde C_{2,T}\,\eps^{\frac{\a}{2}}+C_3\,e^{-\frac{\gamma t}{\eps}}.
\end{align*}
Besides the low regularity assumption on the test function $\psi$, the above estimate implies that for every $t>0$, the distribution of the fast variable $\hat Y_t^\eps$ will converge weakly to $\mE\hat\zeta^{t,\bar{\hat{X}}_t}$ (where the expectation is taken with respect to $\bar{\hat{X}}_t$) as $\eps\to0$, and the rate of convergence is $\eps^{\frac{\a}{2}\wedge1}+\e^{-\frac{\gamma t}{\eps}}$, which is independent of the index $\b$ of the coefficients (the regularity of the coefficients with respect to the fast motion).

Note that the constant $C_3$ is independent of the time variable.                                         The exponential decay term in the rate is natural since even if the multi-scale system is not fully coupled, i.e., the fast motion does not depend on the slow process, we would have that for every $t>0$,
\begin{align*}
\bigg|\mE \hat\psi(\hat Y_t^\eps)-\int_{\mR^{d_2}}\hat\psi(y)\zeta(\dif y)\bigg|\leq C_3\,\e^{-\frac{\gamma t}{\eps}},
\end{align*}
where $\zeta$ is the unique invariant measure for the fast motion.

  \vspace{1mm}
iv) Similar as in Remark \ref{small} (iv), it can be seen that if $\psi$ satisfies (\ref{dissf2}), then the constant $C_3$ in the  estimate (\ref{we-non-y}) can be replaced by $\kappa\, C_3$.
\er

\subsection{Strong convergence in the averaging principle}

It seems to be difficult to prove the strong convergence of $\hat X_t^\eps$ to $\bar{\hat X}_t$ in Theorem \ref{non-aut} {\it (i)} directly due to the low regularity (only H\"older continuous) of the coefficients. For this reason, we shall use the Zvonkin's transformation as in  \cite{RX2,V0} to transform the equations of $\hat X_t^\eps$ and $\bar{\hat X}_t$ into new ones.

For $T>0$ and $\l>0$, consider the following backward PDE on $[0,T]\times\mR^{d_1}$:
\begin{equation}\label{pde}
\left\{ \begin{aligned}
&\p_t v(t,x)+\bar{\hat b}(t,x)\cdot\p_x v(t,x)+\frac{1}{2}\T\big(\hat a(t,x)\cdot\p_x^2 v(t,x)\big)\\
&\qquad\qquad\qquad\qquad\qquad\qquad\!+\bar{\hat b}(t,x)=\l v(t,x),\quad t\in[0,T),\\
&v(T,x)=0,
\end{aligned} \right.
\end{equation}
where $\hat a(t,x):=\hat\sigma\hat\sigma^*(t,x)$. Under our assumptions on the coefficients and according to Corollary \ref{cor}, we have $\bar{\hat b}\in C_b^{\a/2,\a}$. Thus, there exists a unique solution $v\in C_b^{1+\a/2,2+\a}([0,T]\times\mR^{d_1})$ to the equation (\ref{pde}), see e.g. \cite[Chapter IV, Section 5]{LSU}. Moreover, we  have $\p_xv\in C_b^{(1+\a)/2,1+\a}([0,T]\times\mR^{d_1})$, and we can choose $\l$  large enough such that
\begin{align}\label{v}
|\p_xv(t,x)|\leq\frac{1}{2},\quad\forall (t,x)\in[0,T]\times\mR^{d_1}.
\end{align}
Define the transformation function by
$$
\Gamma(t,x):=x+v(t,x),
$$
and then define two new processes by
\begin{align}\label{newp}
\bar V_t:=\Gamma(t,\bar{\hat X}_t)\quad{\rm{and}}\quad \hat V_t^\eps:=\Gamma(t,\hat X_t^\eps).
\end{align}
We have the following result.

\bl
Let $\bar V_t$ and $\hat V_t^\eps$ be defined by (\ref{newp}). Then we have for every $t\in[0,T]$,
\begin{align}\label{v1}
\dif \bar V_t=\l v(t,\bar{\hat X}_t)\dif t+\p_x\Gamma(t,\bar{\hat X}_t)\cdot\hat\sigma(t,\bar{\hat X}_t)\dif W_t^1, \qquad \bar V_0=\Gamma(0,\xi),
\end{align}
and
\begin{align}\label{v2}
\dif \hat V_t^\eps&=\l v(t,\hat X^\eps_t)\dif t+\p_x\Gamma(t,\hat X_t^\eps)\cdot\sigma_\eps(t,\hat X_t^\eps)\dif W_t^1\no\\
&\quad+\frac{1}{2}\T\Big(\big[a_\eps(t,\hat X_t^\eps)-\hat a(t,\hat X_t^\eps)\big]\cdot\p_x^2v(t,\hat X_t^\eps)\Big)\dif t\no\\
&\quad+\big[b_\eps(t,\hat X_t^\eps,\hat Y_t^\eps)-\bar{\hat b}(t,\hat X_t^\eps)\big]\cdot\p_x\Gamma(t,\hat X_t^\eps)\dif t,\qquad \hat V_0^\eps=\Gamma(0,\xi).
\end{align}
\el
\begin{proof}
We only proof (\ref{v2}) since the proof of (\ref{v1}) is easier and follows by the same argument. Using It\^{o}'s formula for $v(t,\hat X_t^\eps)$, we deduce that
\begin{align*}
v(t,\hat X_t^\eps)&=v(0,\xi)+\int_0^t\p_xv(s,\hat X_s^\eps)\cdot\sigma_\eps(s,\hat X_s^\eps)\dif W_s^1+\int_0^t\Big(\p_sv(s,\hat X_s^\eps)\\
&\quad+b_\eps(s,\hat X_s^\eps,\hat Y_s^\eps)\cdot\p_xv(s,\hat X_s^\eps)+\frac{1}{2}\T\big(a_\eps(s,\hat X_s^\eps)\cdot\p^2_xv(s,\hat X_s^\eps)\big)\Big)\dif s\\
&=v(0,\xi)+\int_0^t\p_xv(s,\hat X_s^\eps)\cdot\sigma_\eps(s,\hat X_s^\eps)\dif W_s^1+\int_0^t\Big(\p_sv(s,\hat X_s^\eps)\\
&\quad+\bar{\hat b}(s,\hat X_s^\eps)\cdot\p_xv(s,\hat X_s^\eps)+\frac{1}{2}\T\big(\hat a(s,\hat X_s^\eps)\cdot\p^2_xv(s,\hat X_s^\eps)\big)\Big)\dif s\\
&\quad+\int_0^t\big[b_\eps(s,\hat X_s^\eps,\hat Y_s^\eps)-\bar{\hat b}(s,\hat X_s^\eps)\big]\cdot\p_xv(s,\hat X_s^\eps)\dif s\\
&\quad+\frac{1}{2}\int_0^t\T\Big(\big[a_\eps(s,\hat X_s^\eps)-\hat a(s,\hat X_s^\eps)\big]\cdot\p_x^2v(s,\hat X_s^\eps)\Big)\dif s.
\end{align*}
Since $v$ satisfies the equation (\ref{pde}), we further have
\begin{align*}
v(t,\hat X_t^\eps)
&=v(0,\xi)+\int_0^t\p_xv(s,\hat X_s^\eps)\cdot\sigma_\eps(s,\hat X_s^\eps)\dif W_s^1+\int_0^t \big(\l v(s,\hat X_s^\eps)-\bar{\hat b}(s,\hat X_s^\eps)\big)\dif s\\
&\quad+\int_0^t\big[b_\eps(s,\hat X_s^\eps,\hat Y_s^\eps)-\bar{\hat b}(s,\hat X_s^\eps)\big]\cdot\p_xv(s,\hat X_s^\eps)\dif s\\
&\quad+\frac{1}{2}\int_0^t\T\Big(\big[a_\eps(s,\hat X_s^\eps)-\hat a(s,\hat X_s^\eps)\big]\cdot\p_x^2v(s,\hat X_s^\eps)\Big)\dif s.
\end{align*}
Pulsing this with the equation for $\hat X_t^\eps$, we obtain (\ref{v2}).
\end{proof}

Next, we  derive a strong fluctuation estimate by using the technique of the Poisson equation.
Let $f(t,x,y)$ be a function satisfying the centering condition, i.e.,
\begin{align}\label{cenf1}
\int_{\mR^{d_2}}f(t,x,y)\hat\zeta^{t,x}(\dif y)=0,\qquad \forall  (t,x)\in\mR_+\times\mR^{d_1},
\end{align}
where $\hat\zeta^{t,x}(\dif y)$ is the unique invariant measure of the equation (\ref{frozen1}). We introduce the following Poisson equation in $\mR^{d_2}$:
\begin{align}\label{pss}
\hat\sL_0(t,x,y)\hat U(t,x,y)=-f(t,x,y),
\end{align}
where $(t,x)\in\mR_+\times\mR^{d_1}$ are regarded as parameters, and the operator $\hat\sL_0$  is the infinitesimal generator of $\hat Y^{t,x}_s$ satisfying the equation (\ref{frozen1}), i.e.,
\begin{align}\label{ophat}
\hat\sL_0\varphi(y):=\hat\sL_0(t,x,y)\varphi(y):=\frac{1}{2}\T\Big(\hat\cG(t,x,y)\cdot\p_y^2\varphi(y)\Big)
+\hat F(t,x,y)\cdot\p_y\varphi(y),
\end{align}
where $\hat\cG(t,x,y):=\hat G\hat G^*(t,x,y)$ and $\varphi(y)$ sufficiently smooth test functions  on $\mR^{d_2}$. The following result gives an estimate for the fluctuations of the process $f(s,\hat X_s^\eps,\hat Y_s^\eps)$ over the time interval $[0, t]$.

\bl\label{flu}
Assume that {\bf(A)} holds,  $\hat F,\hat G\in C^{\a/2,\a,\b}_p$ with $0<\a,\b\leq2$, $b_\eps,F_\eps, G_\eps\in L_p^\infty$ and $\sigma_\eps\in L^\infty$. Then for every $t\geq 0$ and $f\in C_p^{\a/2,\a,\b}$ satisfying (\ref{cenf1}), we have
\begin{align*}
\mE\left|\int_0^tf(s,\hat X_s^\eps,\hat Y_s^\eps)\dif s\right|^2&\leq C_t\bigg(\eps^{\a\wedge1}+\int_0^t \Big[\|F_\eps(s,\cdot,\cdot)-\hat F(s,\cdot,\cdot)\|^2_{L_p^\infty}\no\\
&\qquad\qquad+\|G_\eps(s,\cdot,\cdot)-\hat G(s,\cdot,\cdot)\|^2_{L_p^\infty}\Big]\dif s\bigg),
\end{align*}
where $C_t>0$ is a constant independent of $\eps$ and $\beta$.
\el
\begin{proof}
By the assumptions that $f\in C_p^{\a/2,\a,\b}$ satisfying (\ref{cenf1}), $\hat F, \hat G\in C_p^{\a/2,\a,\b}$ and according to Theorem \ref{lem1}, there is a unique solution $\hat U\in C_p^{\a/2,\a,2+\b}$ to the Poisson equation (\ref{pss}). Let $\hat U_n$ be the mollifying approximation of $\hat U$ defined as in (\ref{fnn}) (which does not depend on the $\mu$-variable here). Then by It\^o's formula, we have
\begin{align*}
&\hat U_n(t,\hat X_t^\eps,\hat Y_t^\eps)=\hat U_n(0,\xi,\eta)+\int_0^t\p_s\hat U_n(s,\hat X_s^\eps,\hat Y_s^\eps)\dif s+M_n^1(t)+\frac{1}{\sqrt{\eps}}M_n^2(t)\\
&\quad+\int_0^t\[b_\eps(s,\hat X_s^\eps,\hat Y_s^\eps)\cdot\p_x\hat U_n(s,\hat X_s^\eps,\hat Y_s^\eps)
+\frac{1}{2}\T\big(a_\eps(s,\hat X_s^\eps)\cdot\p^2_x\hat  U_n(s,\hat X_s^\eps,\hat Y_s^\eps)\big)\]\dif s\\
&\quad+\frac{1}{\eps}\int_0^t\[F_\eps(s,\hat X_s^\eps,\hat Y_s^\eps)\cdot\p_y\hat U_n(s,\hat X_s^\eps,\hat Y_s^\eps)
+\frac{1}{2}\T\big(\cG_\eps(s,\hat X_s^\eps,\hat Y_s^\eps)\cdot\p^2_y\hat U_n(s,\hat X_s^\eps,\hat Y_s^\eps)\big)\]\dif s,
\end{align*}
where $a_\eps(t,x):=\sigma_\eps\sigma_\eps^*(t,x)$, $\cG_\eps(t,x,y):=G_\eps G_\eps^*(t,x,y)$, and $M_n^1(t)$, $M_n^2(t)$ are martingales given by
\begin{align*}
&M_n^1(t):=\int_0^t\p_x\hat U_n(s,\hat X_s^\eps,\hat Y_s^\eps)\cdot\sigma_\eps(s,\hat X_s^\eps)\dif W_s^1,\\
&M_n^2(t):=\int_0^t\p_y\hat U_n(s,\hat X_s^\eps,\hat Y_s^\eps)\cdot G_\eps(s,\hat X_s^\eps,\hat Y_s^\eps)\dif W_s^2.
\end{align*}
This together with (\ref{pss}) yields that
\begin{align}\label{flu2}
&\int_0^tf(s,\hat X_s^\eps,\hat Y_s^\eps)\dif s
=-\int_0^t\hat\sL_0(s,\hat X_s^\eps,\hat Y_s^\eps)\hat U_n(s,\hat X_s^\eps,\hat Y_s^\eps)\dif s\no\\
&\quad+\int_0^t\hat\sL_0(s,\hat X_s^\eps,\hat Y_s^\eps)\[\hat U_n(s,\hat X_s^\eps,\hat Y_s^\eps)-\hat U(s,\hat X_s^\eps,\hat Y_s^\eps)\]\dif s\no\\
&=\eps\[\hat U_n(0,\xi,\eta)-\hat U_n(t,\hat X_t^\eps,\hat Y_t^\eps)\]+\eps\int_0^t\p_s\hat U_n(s,\hat X_s^\eps,\hat Y_s^\eps)\dif s+\eps M_n^1(t)+\sqrt{\eps}M_n^2(t)\no\\
&\quad+\eps\int_0^t\[b_\eps(s,\hat X_s^\eps,\hat Y_s^\eps)\cdot\p_x\hat U_n(s,\hat X_s^\eps,\hat Y_s^\eps)
+\frac{1}{2}\T\big(a_\eps(s,\hat X_s^\eps)\cdot\p^2_x\hat U_n(s,\hat X_s^\eps,\hat Y_s^\eps)\big)\]\dif s\no\\
&\quad+\int_0^t\[F_\eps(s,\hat X_s^\eps,\hat Y_s^\eps)-\hat F(s,\hat X_s^\eps,\hat Y_s^\eps)\]\cdot\p_y\hat U_n(s,\hat X_s^\eps,\hat Y_s^\eps)\dif s\no\\
&\quad+\frac{1}{2}\int_0^t\T\([\cG_\eps(s,\hat X_s^\eps,\hat Y_s^\eps)-\hat\cG(s,\hat X_s^\eps,\hat Y_s^\eps)]\cdot\p^2_y\hat U_n(s,\hat X_s^\eps,\hat Y_s^\eps)\)\dif s\no\\
&\quad+\int_0^t\hat\sL_0(s,\hat X_s^\eps,\hat Y_s^\eps)\[\hat U_n(s,\hat X_s^\eps,\hat Y_s^\eps)
-\hat U(s,\hat X_s^\eps,\hat Y_s^\eps)\]\dif s.
\end{align}
Taking  absolute value and  expectation on both sides of (\ref{flu2}), we have
\begin{align*}
\cU(\eps)&:=\mE\left|\int_0^tf(s,\hat X_s^\eps,\hat Y_s^\eps)\dif s\right|^2\leq C_t\Big[\eps^2\mE|\hat U_n(0,\xi,\eta)|^2+\eps^2\mE|\hat U_n(t,\hat X_t^\eps,\hat Y_t^\eps)|^2\\
&\quad+\eps^2\mE|M_n^1(t)|^2+\eps\mE|M_n^2(t)|^2\Big]+C_t\,\eps^2\int_0^t\mE|\p_s\hat U_n(s,\hat X_s^\eps,\hat Y_s^\eps)|^2\dif s\\
&\quad+C_t\,\eps^2\mE\bigg|\int_0^t\[b_\eps(s,\hat X_s^\eps,\hat Y_s^\eps)\cdot\p_x\hat U_n(s,\hat X_s^\eps,\hat Y_s^\eps)\\
&\qquad\quad+\frac{1}{2}\T\big(a_\eps(s,\hat X_s^\eps)\cdot\p^2_x\hat U_n(s,\hat X_s^\eps,\hat Y_s^\eps)\big)\]\dif s\bigg|^2\\
&\quad+C_t\,\mE\left|\int_0^t\[F_\eps(s,\hat X_s^\eps,\hat Y_s^\eps)-\hat F(s,\hat X_s^\eps,\hat Y_s^\eps)\]\cdot\p_y\hat U_n(s,\hat X_s^\eps,\hat Y_s^\eps)\dif s\right|^2\\
&\quad+C_t\,\mE\left|\int_0^t\T\([\cG_\eps(s,\hat X_s^\eps,\hat Y_s^\eps)-\hat\cG(s,\hat X_s^\eps,\hat Y_s^\eps)]\cdot\p^2_y\hat U_n(s,\hat X_s^\eps,\hat Y_s^\eps)\)\dif s\right|^2\\
&\quad+C_t\,\mE\left|\int_0^t\hat\sL_0(s,\hat X_s^\eps,\hat Y_s^\eps)\[\hat U_n(s,\hat X_s^\eps,\hat Y_s^\eps)
-\hat U(s,\hat X_s^\eps,\hat Y_s^\eps)\]\dif s\right|^2=:\sum_{i=1}^6\cU_i(\eps).
\end{align*}
Note that under (\ref{hatF}), we have for any $q\geq2$,
\begin{align}\label{bouyy}
\sup_{t\geq0}\mE|\hat Y_t^\eps|^q\leq C_0\,(1+\mE|\eta|^q).
\end{align}
Thus by Theorem \ref{lem1} we  derive that
\begin{align*}
\mE|\hat U_n(0,\xi,\eta)|^2+\mE|\hat U_n(t,\hat X_t^\eps,\hat Y_t^\eps)|^2\leq C_1\big(1+\mE|\eta|^{2p}+\mE|\hat Y_t^\eps|^{2p}\big)<\infty.
\end{align*}
At the same time, using the Burkholder-Davis-Gundy inequality, we get
\begin{align*}
\mE|M_n^2(t)|^2\leq C_1\int_0^t\big(1+\mE|\hat Y_s^\eps|^{4p}\big)\dif s<\infty,
\end{align*}
and in view of (\ref{fn2}),
\begin{align*}
\mE|M_n^1(t)|^2\leq C_1\,n^{2(1-\a)}\int_0^t\big(1+\mE|\hat Y_s^\eps|^{2p}\big)\dif s\leq C_1\,n^{2(1-\a)}.
\end{align*}
Consequently, we have
\begin{align*}
\cU_1(\eps)\leq C_1\,(\eps+\eps^2n^{2(1-\a)}).
\end{align*}
Similarly, by (\ref{fn2}), (\ref{fn3}), (\ref{bouyy}) and the assumptions that $b_{\eps}\in L_p^{\infty}$ and $\sigma_\eps\in L^\infty$, we can get
\begin{align*}
\cU_2(\eps)+\cU_3(\eps)&\leq C_2\,\eps^2n^{2(2-\a)}\int_0^t\big(1+\mE|\hat Y_s^\eps|^{4p}\big)\dif s\leq C_3\,\eps^2n^{2(2-\a)}.
\end{align*}
As for $\cU_4(\eps)$, applying Theorem \ref{lem1} again we have
\begin{align*}
\cU_4(\eps)&\leq C_4\int_0^t \|F_\eps(s,\cdot,\cdot)-\hat F(s,\cdot,\cdot)\|^2_{L_p^\infty}\cdot\big(1+\mE|\hat Y_s^\eps|^{4p}\big)\dif s\\
&\leq C_4\int_0^t \|F_\eps(s,\cdot,\cdot)-\hat F(s,\cdot,\cdot)\|^2_{L_p^\infty}\dif s.
\end{align*}
At the same time, we also have
\begin{align*}
\cU_5(\eps)&\leq C_5\int_0^t \|G_\eps(s,\cdot,\cdot)-\hat G(s,\cdot,\cdot)\|^2_{L_p^\infty}\cdot\big(1+\mE|\hat Y_s^\eps|^{6p}\big)\dif s\\
&\leq C_5\int_0^t \|G_\eps(s,\cdot,\cdot)-\hat G(s,\cdot,\cdot)\|^2_{L_p^\infty}\dif s.
\end{align*}
Finally, due to $\p_y^2\hat U\in C_p^{\a/2,\a,\b}$ and the fact that
$$
\p_y^2\hat U_n(t,x,y)=\p_y^2\hat U(\cdot,\cdot,y)\ast\rho^n_1\ast\rho^n_2,
$$
we derive by (\ref{fn1}) that
\begin{align*}
\cU_6(\eps)&\leq C_6\,\mE\bigg(\int_0^t\sum_{\ell=1,2}\big\|\p_y^\ell\hat U_n(s,\hat X_s^\eps,\hat Y_s^\eps)-\p_y^\ell\hat U(s,\hat X_s^\eps,\hat Y_s^\eps)\big\|^2\cdot\big(1+|\hat Y_s^\eps|^{4p}\big)\dif s\bigg)\\
&\leq C_6\,n^{-2\a}\int_0^t\big(1+\mE|\hat Y_s^\eps|^{6p}\big)\dif s\leq C_6\, n^{-2\a}.
\end{align*}
Combining the above computations, we arrive at
\begin{align*}
\cU(\eps)&\leq C_t\bigg(\eps+\eps^2n^{2(2-\a)}+n^{-2\a}\\
&\quad+\int_0^t \big[\|F_\eps(s,\cdot,\cdot)-\hat F(s,\cdot,\cdot)\|^2_{L_p^\infty}+\|G_\eps(s,\cdot,\cdot)-\hat G(s,\cdot,\cdot)\|^2_{L_p^\infty}\big]\dif s\bigg).
\end{align*}
Taking $n=\eps^{-1/2}$ we get
$$
\cU(\eps)\leq C_t\bigg(\eps^{\a\wedge1}+\int_0^t\big[\|F_\eps(s,\cdot,\cdot)-\hat F(s,\cdot,\cdot)\|^2_{L_p^\infty}+\|G_\eps(s,\cdot,\cdot)-\hat G(s,\cdot,\cdot)\|^2_{L_p^\infty}\big]\dif s\bigg).
$$
The proof is finished.
\end{proof}

Now, we are in the position to give:
\begin{proof}[{\bf Proof of Theorem \ref{non-aut} (i)}]
Taking $\l$ large enough such that (\ref{v}) holds and by the definition (\ref{newp}), we obtain that for every $t\in[0,T]$,
\begin{align}\label{xv}
\mE\big|\hat X_t^\eps-\bar{\hat X}_t\big|^2\leq4\mE\big|\hat V_t^\eps-\bar V_t\big|^2.
\end{align}
In view of (\ref{v1}) and (\ref{v2}), we have
\begin{align*}
\hat V_t^\eps-\bar V_t&=\int_0^t\l\big[v(s,\hat X_s^\eps)-v(s,\bar{\hat X}_s)\big]\dif s+\int_0^t\big[\hat b(s,\hat X_s^\eps,\hat Y_s^\eps)-\bar{\hat b}(s,\hat X_s^\eps)\big]\cdot\p_x\Gamma(s,\hat X_s^\eps)\dif s\\
&\quad+\int_0^t\big[b_\eps(s,\hat X_s^\eps,\hat Y_s^\eps)-\hat b(s,\hat X_s^\eps,\hat Y_s^\eps)\big]\cdot\p_x\Gamma(s,\hat X_s^\eps)\dif s\\
&\quad+\frac{1}{2}\int_0^t\T\(\big[a_\eps(s,\hat X_s^\eps)-\hat a(s,\hat X_s^\eps)\big]\cdot\p_x^2v(s,\hat X_s^\eps)\)\dif s\\
&\quad+\int_0^t\[\p_x\Gamma(s,\hat X_s^\eps)\cdot\sigma_\eps(s,\hat X_s^\eps)-\p_x\Gamma(s,\bar{\hat X}_s)\cdot\hat\sigma(s,\bar{\hat X}_s)\]\dif W_s^1.
\end{align*}
Taking expectation from both sides of the above equality, we further have that there exists a constant $C_t>0$ such that
\begin{align*}
\mE\big|\hat V_t^\eps-\bar V_t\big|^2&\leq C_t\int_0^t\mE\big|v(s,\hat X_s^\eps)-v(s,\bar{\hat X}_s)\big|^2\dif s\\
&\quad+C_t\,\mE\left|\int_0^t\big[\hat b(s,\hat X_s^\eps,\hat Y_s^\eps)-\bar{\hat b}(s,\hat X_s^\eps)\big]\cdot\p_x\Gamma(s,\hat X_s^\eps)\dif s \right|^2\\
&\quad+C_t\int_0^t\mE\Big|\big[b_\eps(s,\hat X_s^\eps,\hat Y_s^\eps)-\hat b(s,\hat X_s^\eps,\hat Y_s^\eps)\big]\cdot\p_x\Gamma(s,\hat X_s^\eps)\Big|^2\dif s\\
&\quad+C_t\int_0^t\mE\Big|\T\(\big[a_\eps(s,\hat X_s^\eps)-\hat a(s,\hat X_s^\eps)\big]\cdot\p_x^2v(s,\hat X_s^\eps)\)\Big|^2\dif s\\
&+C_t\int_0^t\mE\Big|\p_x\Gamma(s,\hat X_s^\eps)\cdot\sigma_\eps(s,\hat X_s^\eps)-\p_x\Gamma(s,\bar{\hat X}_s)\cdot\hat\sigma(s,\bar{\hat X}_s)\Big|^2\dif s=:\sum_{i=1}^5\cI_i(\eps).
\end{align*}
Below, we estimate these terms one by one. For the first term, since $v(t,\cdot)\in C_b^{2+\a}(\mR^{d_1})$, we have
\begin{align*}
\cI_1(\eps)\leq C_1\int_0^t\mE\big|\hat X_s^\eps-\bar{\hat X}_s\big|^2\dif s,
\end{align*}
To control the second term,  note that by the definition (\ref{barhatb}), the function $\big[\hat b(t,x,y)-\bar{\hat b}(t,x)\big]\cdot\p_x\Gamma(t,x)$ satisfies the centering condition (\ref{cenf1}) and belongs to $C_p^{\a/2,\a,\b}$.
Thus by Lemma \ref{flu} we obtain
\begin{align*}
\cI_2(\eps)\leq C_2\bigg(\eps^{\a\wedge1}+\int_0^t\big[\|F_\eps(s,\cdot,\cdot)-\hat F(s,\cdot,\cdot)\|^2_{L_p^\infty}+\|G_\eps(s,\cdot,\cdot)-\hat G(s,\cdot,\cdot)\|^2_{L_p^\infty}\big]\dif s\bigg).
\end{align*}
As for $\cI_3(\eps)$, by (\ref{v}) and the fact that $\p_x\Gamma(t,x)=\mI_{d_1}+\p_xv(t,x)$, we have
\begin{align*}
\cI_3(\eps)&\leq C_3\int_0^t\|b_\eps(s,\cdot,\cdot)-\hat b(s,\cdot,\cdot)\|^2_{L_p^\infty}\cdot\big(1+\mE|\hat Y_s^\eps|^{2p}\big)\dif s\\
&\leq C_3\int_0^t\|b_\eps(s,\cdot,\cdot)-\hat b(s,\cdot,\cdot)\|^2_{L_p^\infty}\dif s.
\end{align*}
Similarly, since $\sigma_\eps\in L^\infty$ and the function $\p_x^2v(t,\cdot)\in C_b^\a(\mR^{d_1})$, we get
\begin{align*}
\cI_4(\eps)\leq C_4\int_0^t\|\sigma_\eps(s,\cdot)-\hat \sigma(s,\cdot)\|^2_{L^\infty}\dif s.
\end{align*}
Finally, by the assumption that $\hat\sigma(t,\cdot)\in C_b^1(\mR^{d_1})$ we deduce
\begin{align*}
\cI_5(\eps)\leq C_5\int_0^t\mE\big|\hat X_s^\eps-\bar{\hat X}_s\big|^2\dif s+C_5\int_0^t\|\sigma_\eps(s,\cdot)-\hat \sigma(s,\cdot)\|^2_{L^\infty}\dif s.
\end{align*}
Combining the above estimates and in view of (\ref{xv}), we arrive at
\begin{align*}
&\mE\big|\hat X_t^\eps-\bar{\hat X}_t\big|^2\leq C_t\int_0^t\mE\big|\hat X_s^\eps-\bar{\hat X}_s\big|^2\dif s+C_t\,\eps^{\a\wedge1}\\
&\qquad+C_t\int_0^t\|b_\eps(s,\cdot,\cdot)-\hat b(s,\cdot,\cdot)\|^2_{L_p^\infty}\dif s+C_t\int_0^t\|\sigma_\eps(s,\cdot)-\hat \sigma(s,\cdot)\|^2_{L^\infty}\dif s\\
&\qquad+C_t\int_0^t\big[\|F_\eps(s,\cdot,\cdot)-\hat F(s,\cdot,\cdot)\|^2_{L_p^\infty}+\|G_\eps(s,\cdot,\cdot)-\hat G(s,\cdot,\cdot)\|^2_{L_p^\infty}\big]\dif s.
\end{align*}
This in turn implies the desired result by Gronwall's inequality. Thus the proof is completed.
\end{proof}

\subsection{Convergence of the distribution of  the slow process}

Fix $T>0$ below. To prove the convergence of the distribution of  $\hat X_t^\eps$ of the form (\ref{we-non-x}) in Theorem \ref{non-aut} $(ii)$, we need to establish a weak fluctuation estimate for the system (\ref{sde22}) with a function $f$ involving the distribution of  $\hat X_t^\eps$. Namely, let $f(t,x,\mu,y)$ be a function satisfying the centering condition, i.e.,
\begin{align}\label{cenfloc1}
\int_{\mR^{d_2}}f(t,x,\mu,y)\hat\zeta^{t,x}(\dif y)=0,\quad\forall (t,x,\mu)\in \mR_+\times\mR^{d_1}\times\sP_2(\mR^{d_1}),
\end{align}
where $\hat\zeta^{t,x}(\dif y)$ is the unique invariant measure of the frozen equation  (\ref{frozen1}). For simplicity, we shall say  that $f\in C_{p,loc}^{\a/2,\a,(1,\a),\b}([0,T))=C_{p,loc}^{\a/2,\a,(1,\a),\b}\big([0,T)\times\mR^{d_1} \times\sP_2(\mR^{d_1})\times\mR^{d_2}\big)$ with $0<\a\leq 2$ (here and below, local is  regarded to the $t$-variable at point $T$) if there exist constants $C_T, p>0$  such that  for $0<\a\leq1$ and every $0\leq s\leq t<T$,
\begin{align}\label{ff}
\begin{split}
&|f(t,x,\mu,y)|+\bigg|\frac{\delta f}{\delta\mu}(t,x,\mu,y)(\tilde x)\bigg|\leq C_T(1+|y|^p)(T-t)^{(\alpha-1)/2},\\
&|f(t,x,\mu,y)-f(s,x,\mu,y)|\leq C_T(1+|y|^p)(T-t)^{-\frac{1}{2}}|t-s|^\frac{\a}{2},\\
&|f(t,x_1,\mu,y)-f(t,x_2,\mu,y)|+\bigg|\frac{\delta f}{\delta\mu}(t,x,\mu,y)(x_1)-\frac{\delta f}{\delta\mu}(t,x,\mu,y)(x_2)\bigg|\\
&\qquad\qquad\qquad\qquad\qquad\qquad\leq C_T(1+|y|^p)(T-t)^{-\frac{1}{2}}|x_1-x_2|^\a,\\
&|f(t,x,\mu,y_1)-f(t,x,\mu,y_2)|\leq C_T(1+|y_1|^p+|y_2|^p)(T-t)^{(\alpha-1)/2}|y_1-y_2|^\b.
\end{split}
\end{align}
Similar to Lemma \ref{fn}, we provide the following result concerning the mollifying approximation of $f\in C_{p,loc}^{\a/2,\a,(1,\a),\b}([0,T))$, the proof is more or less standard and thus omitted.

\bl\label{fnloc}
Given a function $f\in C_{p,loc}^{\a/2,\a,(1,\a),2}([0,T))$ with $0<\a\leq2$, there exists a sequence of functions $f_n\in C_{p}^{1,2,(1,2+\alpha),2}([0,T))$ such that for every $t\in[0,T)$,
\begin{align}
&\|\p_xf_n(t,\cdot,\cdot,y)\|_\infty+\|\p_\mu f_n(t,\cdot,\cdot,y)(\cdot)\|_\infty\leq C_T(T-t)^{-1/2}n^{1-\a\wedge1}(1+|y|^p),\label{fnloc2}\\
&\|\p_tf_n(t,\cdot,\cdot,y)\|_\infty+\|\p^2_xf_n(t,\cdot,\cdot,y)\|_\infty +\|\p_{\tilde x}\p_\mu f_n(t,\cdot,\cdot,y)(\cdot)\|_\infty\no\\
&\qquad\qquad\qquad\qquad\qquad\qquad\qquad\quad\leq C_T(T-t)^{-1/2}n^{2-\a}(1+|y|^p),\label{fnloc3}\\
&|f_n(0,x,\mu,y)|+|f_n(T,x,\mu,y)|\leq n^{1-\alpha\wedge1}(1+|y|^p),\label{fnloc5}
\end{align}
and
\begin{align}
|\p_y^\ell f_n(t,x,\mu,y)-\p_y^\ell f(t,x,\mu,y)|\leq C_T(T-t)^{-1/2}n^{-\a}(1+|y|^p),\quad\ell=0,1,2,\label{fnloc4}
\end{align}
where $C_T>0$ is a constant independent of $n$.
\el

We have the following result for the fluctuations of the process $f(t,\hat X_t^\eps,\cL_{\hat X_t^\eps},\hat Y_t^\eps)$ over the time interval $[0, T]$.

\bl\label{fluloc}
Assume that {\bf(A)} holds, $\hat F, \hat G\in C^{\a/2,\a,\b}_p$ with $0<\a,\b\leq2$, $b_\eps,F_\eps, G_\eps\in L_p^\infty$ and $\sigma_\eps\in L^\infty$. Then for every $f\in C_{p,loc}^{\a/2,\a,(1,\a),\b}([0,T))$ satisfying (\ref{cenfloc1}), we have
\begin{align*}
\bigg|\mE\int_0^Tf(s,\hat X_s^\eps,\cL_{\hat X_s^\eps},\hat Y_s^\eps)\dif s\bigg|&\leq C_T\bigg(\eps^{\tfrac{\a}{2}}+\int_0^T(T-s)^{\frac{\a-1}{2}}\Big[\|F_\eps(s,\cdot,\cdot)-\hat F(s,\cdot,\cdot)\|_{L_p^\infty}\\
&\qquad\quad+\|G_\eps(s,\cdot,\cdot)-\hat G(s,\cdot,\cdot)\|_{L_p^\infty}\Big]\dif s\bigg),
\end{align*}
where $C_T>0$ is a constant independent of $\eps$ and $\beta$.
\el

\begin{proof}
By the assumptions that $f\in C_{p,loc}^{\a/2,\a,(1,\a),\b}([0,T))$ satisfying (\ref{cenfloc1}), $\hat F, \hat G\in C^{\a/2,\a,\b}_p$ and using Theorem \ref{lem1}, there is a unique solution $\hat U\in C_{p,loc}^{\a/2,\a,(1,\a),2+\b}([0,T))$ to the following Poisson equation:
\begin{align}\label{pssloc}
\hat\sL_0(t,x,y)\hat U(t,x,\mu,y)=-f(t,x,\mu,y),
\end{align}
where $(t,x,\mu)\in[0,T)\times\mR^{d_1}\times\sP_2(\mR^{d_1})$ are regarded as parameters, and the operator $\hat\sL_0$ is defined by (\ref{ophat}). Let $\hat U_n$ be the mollifying approximation of $\hat U$ given as in Lemma \ref{fnloc}.
Then by It\^o's formula (see e.g. \cite[Theorem 7.1]{BLPR} or \cite[Proposition 2.1]{CF2}                  ), we have
\begin{align*}
\hat U_n(T,\hat X_T^\eps,\cL_{\hat X_T^\eps},\hat Y_T^\eps)&=\hat U_n(0,\xi,\mu,\eta)+\!\int_0^T\!\p_s\hat U_n(s,\hat X_s^\eps,\cL_{\hat X_s^\eps},\hat Y_s^\eps)\dif s+M_n^1(T)\\
&\quad+\frac{1}{\sqrt{\eps}}M_n^2(T)+\int_0^T\[b_\eps(s,\hat X_s^\eps,\hat Y_s^\eps)\cdot\p_x\hat U_n(s,\hat X_s^\eps,\cL_{\hat X_s^\eps},\hat Y_s^\eps)\\
&\qquad\qquad\quad+\frac{1}{2}\T\big(a_\eps(s,\hat X_s^\eps)\cdot\p^2_x\hat U_n(s,\hat X_s^\eps,\cL_{\hat X_s^\eps},\hat Y_s^\eps)\big)\]\dif s\\
&\quad+\int_0^T\tilde\mE\[b_\eps(s,\tilde{\hat X}_s^\eps,\tilde{\hat Y}_s^\eps)\cdot\p_\mu\hat U_n(s,\hat X_s^\eps,\cL_{\hat X_s^\eps},\hat Y_s^\eps)(\tilde{\hat X}_s^\eps)\\
&\qquad\qquad\quad+\frac{1}{2}\T\big(a_\eps(s,\tilde{\hat X}_s^\eps)\cdot\p_{\tilde x}\p_\mu\hat U_n(s,\hat X_s^\eps,\cL_{\hat X_s^\eps},\hat Y_s^\eps)(\tilde{\hat X}_s^\eps)\big)\]\dif s\\
&\quad+\frac{1}{\eps}\int_0^T\[F_\eps(s,\hat X_s^\eps,\hat Y_s^\eps)\cdot\p_y\hat U_n(s,\hat X_s^\eps,\cL_{\hat X_s^\eps},\hat Y_s^\eps)\\
&\qquad\qquad\quad+\frac{1}{2}\T\big(\cG_\eps(s,\hat X_s^\eps,\hat Y_s^\eps)\cdot\p^2_y\hat U_n(s,\hat X_s^\eps,\cL_{\hat X_s^\eps},\hat Y_s^\eps)\big)\]\dif s,
\end{align*}
where the process ($\tilde{\hat X}^{\eps}_s,\tilde{\hat Y}^{\eps}_s$) is a copy of the original process $(\hat X^{\eps}_s,\hat Y^{\eps}_s)$ defined on a copy $(\tilde\Omega,\tilde\sF,\tilde\mP)$ of the original probability space $(\Omega,\sF,\mP)$, and for $i=1,2$, $M_{n}^i(T)$ are two martingales defined by
\begin{align*}
&M_n^1(T):=\int_0^T\p_x\hat U_n(s,\hat X_s^\eps,\cL_{\hat X_s^\eps},\hat Y_s^\eps)\cdot \sigma_\eps(s,\hat X_s^\eps)\dif W_s^1,\\
&M_n^2(T):=\int_0^T\p_y\hat U_n(s,\hat X_s^\eps,\cL_{\hat X_s^\eps},\hat Y_s^\eps)\cdot G_\eps(s,\hat X_s^\eps,\hat Y_s^\eps)\dif W_s^2.
\end{align*}
This together with (\ref{pssloc}) yields that
\begin{align*}
&\quad\left|\mE\bigg(\int_0^Tf(s,\hat X_s^\eps,\cL_{\hat X_s^\eps},\hat Y_s^\eps)\dif s\bigg)\right|\leq C_0\left|\mE\bigg(\int_0^T-\hat\sL_0(s,\hat X_s^\eps,\hat Y_s^\eps)\hat U_n(s,\hat X_s^\eps,\cL_{\hat X_s^\eps},\hat Y_s^\eps)\dif s\bigg)\right|\\
&\qquad\quad+ C_0\left|\mE\bigg(\int_0^T\hat\sL_0(s,\hat X_s^\eps,\hat Y_s^\eps)\big[\hat U_n(s,\hat X_s^\eps,\cL_{\hat X_s^\eps},\hat Y_s^\eps)-\hat U(s,\hat X_s^\eps,\cL_{\hat X_s^\eps},\hat Y_s^\eps)\big]\dif s\bigg)\right|\\
&\leq C_0\,\eps\Big[\mE\big|\hat U_n(0,\xi,\mu,\eta)\big|+\mE\big|\hat U_n(T,\hat X_T^\eps,\cL_{\hat X_T^\eps},\hat Y_T^\eps)\big|\Big]+C_0\,\eps\left|\int_0^T\mE\Big[\p_s\hat U_n(s,\hat X_s^\eps,\cL_{\hat X_s^\eps},\hat Y_s^\eps)\Big]\dif s\right|\\
&\quad+C_0\,\eps\bigg|\int_0^T\mE\[b_\eps(s,\hat X_s^\eps,\hat Y_s^\eps)\cdot\p_x\hat U_n(s,\hat X_s^\eps,\cL_{\hat X_s^\eps},\hat Y_s^\eps)\\
&\qquad\qquad+\frac{1}{2}\T\big(a_\eps(s,\hat X_s^\eps)\cdot\p^2_x\hat U_n(s,\hat X_s^\eps,\cL_{\hat X_s^\eps},\hat Y_s^\eps)\big)\]\dif s\bigg|\\
&\quad+C_0\,\eps\bigg|\int_0^T\mE\tilde\mE\[b_\eps(s,\tilde{\hat X}_s^\eps,\tilde{\hat Y}_s^\eps)\cdot\p_\mu\hat U_n(s,\hat X_s^\eps,\cL_{\hat X_s^\eps},\hat Y_s^\eps)(\tilde{\hat X}_s^\eps)\\
&\qquad\qquad+\frac{1}{2}\T\big(a_\eps(s,\tilde{\hat X}_s^\eps)\cdot\p_{\tilde x}\p_\mu\hat U_n(s,\hat X_s^\eps,\cL_{\hat X_s^\eps},\hat Y_s^\eps)(\tilde{\hat X}_s^\eps)\big)\]\dif s\bigg|\\
&\quad+C_0\left|\mE\bigg(\int_0^T\big[F_\eps(s,\hat X_s^\eps,\hat Y_s^\eps)-\hat F(s,\hat X_s^\eps,\hat Y_s^\eps)\big]\cdot\p_y\hat U_n(s,\hat X_s^\eps,\cL_{\hat X_s^\eps},\hat Y_s^\eps)\dif s\bigg)\right|\\
&\quad+C_0\left|\mE\bigg(\int_0^T\T\big([\cG_\eps(s,\hat X_s^\eps,\hat Y_s^\eps)-\hat \cG(s,\hat X_s^\eps,\hat Y_s^\eps)]\cdot\p^2_y\hat U_n(s,\hat X_s^\eps,\cL_{\hat X_s^\eps},\hat Y_s^\eps)\big)\dif s\bigg)\right|\no\\
&\quad+C_0\left|\mE\bigg(\int_0^T\hat\sL_0(s,\hat X_s^\eps,\hat Y_s^\eps)\big[\hat U_n(s,\hat X_s^\eps,\cL_{\hat X_s^\eps},\hat Y_s^\eps)-\hat U(s,\hat X_s^\eps,\cL_{\hat X_s^\eps},\hat Y_s^\eps)\big]\dif s\bigg)\right|=:\sum_{i=1}^7\cV_i(\eps).
\end{align*}
By estimate (\ref{fnloc5}) in Lemma \ref{fnloc} and estimate (\ref{bouyy}), we derive that there exists a constant $C_1>0$ independent of $n$ such that
\begin{align*}
\cV_1(\eps)\leq C_1\,n^{1-\a\wedge1}\eps\big(1+\mE|\eta|^{p}+\mE|\hat Y_T^\eps|^{p}\big)\leq C_1\,n^{1-\a\wedge1}\eps.
\end{align*}
Similarly, by (\ref{fnloc2}), (\ref{fnloc3}) and the assumptions that $b_{\eps}\in L_p^{\infty}$ and $\sigma_\eps\in L^\infty$, we get
\begin{align*}
\cV_2(\eps)+\cV_3(\eps)+\cV_4(\eps)&\leq C_2\,\eps n^{2-\a}\int_0^T(T-s)^{-\frac{1}{2}}\big(1+\mE|\hat Y_s^\eps|^{2p}\big)\dif s\\
&\leq C_2\,\eps n^{2-\a}.
\end{align*}
As for $\cV_5(\eps)$, by Theorem \ref{lem1} and Lemma \ref{fnloc} we have
\begin{align*}
\cV_5(\eps)&\leq C_3\int_0^T(T-s)^{\frac{\a-1}{2}}\|F_\eps(s,\cdot,\cdot)-\hat F(s,\cdot,\cdot)\|_{L_p^\infty}\big(1+\mE|\hat Y_s^\eps|^{2p}\big)\dif s\\
&\leq C_3 \int_0^T(T-s)^{\frac{\a-1}{2}}\|F_\eps(s,\cdot,\cdot)-\hat F(s,\cdot,\cdot)\|_{L_p^\infty}\dif s.
\end{align*}
At the same time, by $G_{\eps}\in L_p^{\infty}$ we deduce that
\begin{align*}
\cV_6(\eps)&\leq C_4\int_0^T(T-s)^{\frac{\a-1}{2}}\|G_\eps(s,\cdot,\cdot)-\hat G(s,\cdot,\cdot)\|_{L_p^\infty}\big(1+\mE|\hat Y_s^\eps|^{3p}\big)\dif s\\
&\leq C_4 \int_0^T(T-s)^{\frac{\a-1}{2}}\|G_\eps(s,\cdot,\cdot)-\hat G(s,\cdot,\cdot)\|_{L_p^\infty}\dif s.
\end{align*}
Finally, using estimate (\ref{fnloc4})
we derive that
\begin{align*}
\cV_7(\eps)&\leq C_5\,\mE\bigg(\int_0^T\sum_{\ell=1,2}\big\|\p_y^\ell\hat U_n(s,\hat X_s^\eps,\cL_{\hat X_s^\eps},\hat Y_s^\eps)-\p_y^\ell\hat U(s,\hat X_s^\eps,\cL_{\hat X_s^\eps},\hat Y_s^\eps)\big\|\cdot\big(1+|\hat Y_s^\eps|^{2p}\big)\dif s\bigg)\\
&\leq C_5\,n^{-\a}\int_0^T(T-s)^{-\frac{1}{2}}\big(1+\mE|\hat Y_s^\eps|^{3p}\big)\dif s\leq C_5\,n^{-\a}.
\end{align*}
Combining the above computations, we arrive at
\begin{align*}
&\left|\mE\bigg(\int_0^Tf(s,\hat X_s^\eps,\cL_{\hat X_s^\eps},\hat Y_s^\eps)\dif s\bigg)\right|\leq C_6\bigg(\eps n^{1-\a\wedge1}+\eps n^{2-\a}+n^{-\a}\\
&\qquad+\int_0^T(T-s)^{\frac{\a-1}{2}}\Big[\|F_\eps(s,\cdot,\cdot)-\hat F(s,\cdot,\cdot)\|_{L_p^\infty}+\|G_\eps(s,\cdot,\cdot)-\hat G(s,\cdot,\cdot)\|_{L_p^\infty}\Big]\dif s\bigg).
\end{align*}
Taking $n=\eps^{-1/2}$ we get
\begin{align*}
&\quad\left|\mE\bigg(\int_0^Tf(s,\hat X_s^\eps,\cL_{\hat X_s^\eps},\hat Y_s^\eps)\dif s\bigg)\right|\leq C_7\bigg(\eps^\frac{\a}{2}\\
&\qquad\qquad+\int_0^T(T-s)^{\frac{\a-1}{2}}\Big[\|F_\eps(s,\cdot,\cdot)-\hat F(s,\cdot,\cdot)\|_{L_p^\infty}+\|G_\eps(s,\cdot,\cdot)-\hat G(s,\cdot,\cdot)\|_{L_p^\infty}\Big]\dif s\bigg).
\end{align*}
The proof is finished.
\end{proof}

In order to prove the weak convergence of $\hat X_t^\eps$ to $\bar{\hat X}_t$, we need to consider the following backward Kolmogorov equation in $[0,T]\times\sP_2(\mR^{d_1})$:
\begin{equation}\label{cp1}
\left\{ \begin{aligned}
&\p_t u(t,\mu)+\int_{\mR^{d_1}}\Big[\bar{\hat b}(t,x)\!\cdot\!\p_\mu u(t,\mu)(x)\\
&\qquad\qquad\qquad\quad+\frac{1}{2}\T\big(\hat a(t,x)\!\cdot\!\p_{x}(\p_\mu u(t,\mu)(x))\big)\Big]\mu(\dif x)=0,\quad t\in[0,T),\\
&u(T,\mu)=\varphi(\mu),
\end{aligned} \right.
\end{equation}
where $\varphi\in C_b^{(2,\a)}(\sP_2(\mR^{d_1}))$ is a given function. Recall that we have $\bar{\hat b}, \hat \sigma\in C_b^{\a/2,\a}$. According to
Theorem \ref{ke}, there exists a unique solution $u\in C_{loc}^{1+\a/2,(2,2+\a)}([0, T)\times\sP_2(\mR^{d_1}))$ to equation (\ref{cp1}) which is given by
\begin{align}\label{cp2}
u(t,\mu)=\varphi\big(\cL_{\bar {\hat X}_{t,T}(\xi)}\big),\quad \forall t\in[0,T],
\end{align}
where $\cL_\xi=\mu$, and for $0\leq s<t$, the process $\bar {\hat X}_{s,t}(\xi)$ satisfies SDE (\ref{ave1}) with initial value $\xi$ at time $s$, i.e.,
\begin{align*}
\dif \bar{\hat X}_{s,t}(\xi)=\bar{\hat b}(t,\bar{\hat X}_{s,t}(\xi))\dif t+\hat \sigma(t,\bar{\hat X}_{s,t}(\xi))\dif W^1_t,\qquad  \bar{\hat X}_{s,s}(\xi)=\xi.
\end{align*}
For simplicity, we shall write $\bar{\hat X}_{t}:=\bar{\hat X}_{t}(\xi):=\bar{\hat X}_{0,t}(\xi)$.

Now, we proceed  to give:

\begin{proof}[{\bf Proof of Theorem \ref{non-aut} (ii)}]
For every $\varphi\in C_b^{(2,\a)}(\sP_2(\mR^{d_1}))$, let $u(t,\mu)$ be defined by (\ref{cp2}). Then we have
$$
u(T,\cL_{\hat X_T^\eps})=\varphi(\cL_{\hat X_T^\eps})\quad\text{and} \quad u(0,\mu)=\varphi(\cL_{\bar{\hat X}_T}).
$$
By It\^o's formula, we obtain
\begin{align*}
\sK(\eps)&:=|\varphi(\cL_{\hat X_T^\eps})-\varphi(\cL_{\bar{\hat X}_T})|=|u(T,\cL_{\hat X_T^\eps})-u(0,\mu)|\\
&=\bigg|\mE\bigg(\int_0^T\p_su(s,\cL_{\hat X_s^\eps}) +b_\eps(s,\hat X^{\eps}_s,\hat Y^{\eps}_s)\cdot\p_{\mu}u(s,\cL_{\hat X_s^\eps})(\hat X^{\eps}_s)\\
&\quad+\frac{1}{2}\T\Big(a_\eps(s,\hat X^{\eps}_s)\cdot\p_{x}\big(\p_{\mu}u(s,\cL_{\hat X_s^\eps})(\hat X^{\eps}_s)\big)\Big)\dif s\bigg)\bigg|.
\end{align*}
In view of the equation (\ref{cp1}), we further obtain that
\begin{align*}
\sK(\eps)
&\leq\frac{1}{2}\bigg|\mE\bigg(\int_0^T\T\(\big[a_\eps(s,\hat X^{\eps}_s)-\hat a(s,\hat X^{\eps}_s)\big]\cdot\p_x\big(\p_{\mu}u(s,\cL_{\hat X_s^\eps})(\hat X^{\eps}_s)\big)\)\dif s\bigg)\bigg|\\
&\quad+\bigg|\mE\bigg(\int_0^T\big[b_\eps(s,\hat X^{\eps}_s,\hat Y^{\eps}_s)-\hat b(s,\hat X^{\eps}_s,\hat Y^{\eps}_s)\big]\cdot\p_{\mu}u(s,\cL_{\hat X_s^\eps})(\hat X^{\eps}_s)\dif s\bigg)\bigg|\\
&\quad+\bigg|\mE\bigg(\int_0^T\[\hat b(s,\hat X^{\eps}_s,\hat Y^{\eps}_s)-\bar {\hat b}(s,\hat X^{\eps}_s)\]\cdot\p_{\mu}u(s,\cL_{\hat X_s^\eps})(\hat X^{\eps}_s)\dif s\bigg)\bigg|\\
&=:\sK_1(\eps)+\sK_2(\eps)+\sK_3(\eps).
\end{align*}
According to Theorem \ref{ke}, we have
$$
|\p_{\mu}u(t,\mu)(x)|\leq C_0(T-t)^{(\alpha-1)/2}\quad\text{and}\quad \big|\p_x\big(\p_{\mu}u(t,\mu)(x)\big)\big|\leq C_0(T-t)^{\alpha/2-1}.
$$
As a result, we have
\begin{align*}
\sK_1(\eps)\leq C_1\int_0^T (T-s)^{\frac{\a}{2}-1}\|\sigma_\eps(s,\cdot)-\hat \sigma(s,\cdot)\|_{L^\infty}\dif s,
\end{align*}
and
\begin{align*}
\sK_2(\eps)\leq C_2\int_0^T (T-s)^{\frac{\alpha-1}{2}}\|b_\eps(s,\cdot,\cdot)-\hat b(s,\cdot,\cdot)\|_{L_p^\infty}\dif s.
\end{align*}
It remains to control $\sK_3(\eps)$. Note that by the assumptions on the coefficients, the definition (\ref{barhatb}) and Theorem \ref{ke}, the function $\big[\hat b(t,x,y)-\bar {\hat b}(t,x)\big]\cdot\p_{\mu}u(t,\mu)(x)$ satisfies the centering condition (\ref{cenfloc1}) and belongs to $C_{p,loc}^{\a/2,\a,(1,\a),\b}([0,T))$ (i.e., satisfies the estimates in (\ref{ff})).
As a direct consequence of Lemma \ref{fluloc}, we have
\begin{align*}
\sK_3(\eps)&\leq C_T\bigg(\eps^{\tfrac{\a}{2}}+\int_0^T(T-s)^{\frac{\a-1}{2}}\Big[\|F_\eps(s,\cdot,\cdot)-\hat F(s,\cdot,\cdot)\|_{L_p^\infty}\\
&\qquad\qquad\qquad\quad+\|G_\eps(s,\cdot,\cdot)-\hat G(s,\cdot,\cdot)\|_{L_p^\infty}\Big]\dif s\bigg).
\end{align*}
Consequently, we arrive at
\begin{align*}
\sK(\eps)&\leq C_T\bigg(\eps^{\frac{\a}{2}}+\int_0^T (T-s)^{\frac{\a}{2}-1}\|\sigma_\eps(s,\cdot)-\hat \sigma(s,\cdot)\|_{L^\infty}\dif s\\
&\quad+\int_0^T (T-s)^{\frac{\alpha-1}{2}}\|b_\eps(s,\cdot,\cdot)-\hat b(s,\cdot,\cdot)\|_{L_p^\infty}\dif s\\
&\quad+\int_0^T(T-s)^{\frac{\a-1}{2}}\big[\|F_\eps(s,\cdot,\cdot)-\hat F(s,\cdot,\cdot)\|_{L_p^\infty}+\|G_\eps(s,\cdot,\cdot)-\hat G(s,\cdot,\cdot)\|_{L_p^\infty}\big]\dif s\Big)\bigg),
\end{align*}
which in turn implies the desired result.
\end{proof}

\subsection{Limit for the distribution of the fast motion}

Recall that $\hat Y_s^{t,x}$ satisfies the frozen equation (\ref{frozen1}), $\hat\zeta^{t,x}$ is the unique invariant measure of $\hat Y_s^{t,x}$, and $\tilde{\hat \zeta}^{t,\mu}$ is defined by (\ref{zetamu}).
To prove  the convergence of distribution of the fast process $\hat Y_t^\eps$, we  consider the following Kolmogorov equation on $\mR_+\times\sP_2(\mR^{d_1}\times\mR^{d_2})\times[0,T]$:
\begin{equation}\label{eqV}
\left\{\begin{aligned}
\displaystyle
&\p_s V(s,m;t)-\int_{\mR^{d_1}\times\mR^{d_2}}\bigg[\hat F(t,x,y)\cdot\p_y\frac{\delta V}{\delta m}(s,m;t)(x,y)\\
\displaystyle
&\qquad\qquad\quad+\frac{1}{2}\T\Big(\hat\cG(t,x,y)\cdot\p^2_y\frac{\delta V}{\delta m}(s,m;t)(x,y)\Big)\bigg]m(\dif x,\dif y)=0,\\
\displaystyle
&V(0,m;t)=\psi(\pi_2^*m)-\psi(\tilde{\hat \zeta}^{t,\pi_1^*m}),
\end{aligned}\right.
\end{equation}
where $t\in[0,T]$ is a parameter, $\psi\in C_p^{(2,\b)}(\sP_2(\mR^{d_2}))$ is a given function. Under our assumptions and by Theorem \ref{cp}, there exists a unique solution $V(s,m;t)$ to the equation (\ref{eqV}) which is given by
\begin{align}\label{V}
V(s,m;t):=\psi\big(\cL_{\hat Y^{t,\xi}_s(\eta)}\big)-\psi\big(\tilde{\hat \zeta}^{t,\mu}\big),
\end{align}
where $m=\cL_{(\xi,\eta)}$, $\mu=\pi_1^*m=\cL_\xi$ and $\hat Y^{t,\xi}_s=\hat Y^{t,\xi}_s(\eta)$  satisfies the following equation: for $t\in[0,T]$,
\begin{align*}
&\dif \hat Y^{t,\xi}_s =\hat F(t,\xi,\hat Y^{t,\xi}_s)\dif s+\hat G(t,\xi,\hat Y^{t,\xi}_s)\dif \hat W_s^2,\qquad   \hat Y^{t,\xi}_0=\eta.
\end{align*}
Moreover,  we have by Lemma \ref{inv-mea} that there exist constants $C_0,\gamma>0$ such that for any $s\geq0$, $t\in[0,T]$ and $\cV(y)=1+|y|^p$ with $p\geq 1$,
\begin{align}\label{haty}
\rho_\cV(\cL_{\hat Y_s^{t,\xi}(\eta)},\tilde{\hat\zeta}^{t,\mu})\leq C_0 \,\e^{-\gamma s},
\end{align}
where $\rho_\cV$ is defined by (\ref{tv}).

Now, we proceed to give:
\begin{proof}[{\bf Proof of Theorem \ref{non-aut} (iii)}]
For every $t\geq0$ and $\psi\in C_p^{(2,\b)}$, we write
\begin{align*}
\sJ(\eps)&:=\Big|\psi\big(\cL_{\hat Y_t^\eps}\big)-\psi\big(\tilde{\hat \zeta}^{t,\cL_{\bar{\hat{X}}_t}}\big)\Big|\\
&\leq\Big|\psi\big(\tilde{\hat \zeta}^{t,\cL_{\hat{X}_t^\eps}}\big)-\psi\big(\tilde{\hat \zeta}^{t,\cL_{\bar{\hat{X}}_t}}\big)\Big|+\Big|\psi\big(\cL_{\hat Y^{0,\xi}_{t/\eps}(\eta)}\big)-\psi\big(\tilde{\hat \zeta}^{0,\mu}\big)\Big|\\
&\quad+\Big|\psi\big(\cL_{\hat Y_t^\eps}\big)-\psi\big(\tilde{\hat \zeta}^{t,\cL_{\hat{X}_t^\eps}}\big)-\psi\big(\cL_{\hat Y^{0,\xi}_{t/\eps}(\eta)}\big)+\psi\big(\tilde{\hat \zeta}^{0,\mu}\big)\Big|=:\sum_{i=1}^3\sJ_i(\eps),
\end{align*}
where $\xi$ and $\eta$ are the initial value of $\hat X_t^\eps$ and $\hat Y_t^\eps$ respectively, and $\cL_\xi=\mu$. We proceed to control each term separately.

\vspace{1mm}
\noindent{\it  i) (Control of the first term)} For every $\mu\in\sP_2(\mR^{d_1})$, we  let
\begin{align*}
\tilde \psi(t,\mu):=\psi\big(\tilde{\hat \zeta}^{t,\mu}\big),
\end{align*}
where $\tilde{\hat \zeta}^{t,\mu}$ is defined by (\ref{zetamu}). Then for every $t\geq0$, by the chain rule formula in Lemma \ref{chain} we have
\begin{align*}
\frac{\delta\tilde \psi}{\delta\mu}(t,\mu)( x)=\int_{\mR^{d_2}}\frac{\delta\psi}{\delta\nu}\big(\tilde{\hat \zeta}^{t,\mu}\big)(y_1)\hat\zeta^{t,x}(\dif y_1)
\end{align*}
and
\begin{align*}
\frac{\delta^2\tilde \psi}{\delta\mu^2}(t,\mu)(x,\tilde x)=\int_{\mR^{d_2}}\int_{\mR^{d_2}}\frac{\delta^2\psi}{\delta\nu^2}\big(\tilde{\hat \zeta}^{t,\nu}\big)(y_1,y_2)\hat\zeta^{t,x}(\dif y_1)\hat\zeta^{t,\tilde x}(\dif y_2).
\end{align*}
These together with $\psi\in C^{(2,\b)}_p$ imply that $\tilde \psi(t,\cdot)\in C_b^{(2,0)}$.  Furthermore, by Corollary \ref{cor} (see also \cite[Lemma 3.2]{RX1}), we have that $\tilde \psi(t,\cdot)\in C_b^{(2,\a)}$. As a direct result of the convergence of the distribution for the slow process obtained in (\ref{we-non-x}), we deduce that for every fixed $t\geq 0$,
\begin{align*}
\sJ_1(\eps)&=\Big|\tilde \psi\big(t,\cL_{\hat X_t^\eps}\big)-\tilde \psi\big(t,\cL_{\bar{\hat X}_t}\big)\Big|\\
&\leq C_1\bigg(\eps^{\frac{\a}{2}}+\int_0^t\,(t-s)^{\frac{\a}{2}-1}\cdot\Big[\|b_\eps(s,\cdot,\cdot)-\hat b(s,\cdot,\cdot)\|_{L_p^\infty}+\|\sigma_\eps(s,\cdot)-\hat \sigma(s,\cdot)\|_{L^\infty}\no\\
&\qquad\quad+\|F_\eps(s,\cdot,\cdot)-\hat F(s,\cdot,\cdot)\|_{L_p^\infty}+\|G_\eps(s,\cdot,\cdot)-\hat G(s,\cdot,\cdot)\|_{L_p^\infty}\Big]\dif s\bigg),
\end{align*}
where $C_1>0$ is a constant depending on $T$.

\vspace{1mm}
\noindent{\it ii) (Control of the second term)} For the second term $\sJ_2(\eps)$, by (\ref{lin}) and the assumption that $\psi\in C^{(2,\b)}_p$,  we have
\begin{align*}
\sJ_2(\eps)&=\left|
\int_0^1\int_{\mR^{d_2}}\frac{\delta\psi}{\delta\nu}\Big(\tilde{\hat \zeta}^{0,\mu}+\theta\big(\cL_{\hat Y^{0,\xi}_{t/\eps}(\eta)}-\tilde{\hat \zeta}^{0,\mu}\big)\Big)(y)\big(\cL_{\hat Y^{0,\xi}_{t/\eps}(\eta)}-\tilde{\hat \zeta}^{0,\mu}\big)(\dif y)\dif \theta\right|\\
&\leq C_2\int_{\mR^{d_2}}\big(1+\cV(y)\big)\big|\cL_{\hat Y^{0,\xi}_{t/\eps}(\eta)}-\tilde{\hat \zeta}^{0,\mu}\big|(\dif y)\leq C_2\cdot\rho_{\cV}\big(\cL_{\hat Y^{0,\xi}_{t/\eps}(\eta)},\tilde{\hat \zeta}^{0,\mu}\big).
\end{align*}
Consequently, using (\ref{haty}) with $s=t/\eps$ we get
\begin{align*}
\sJ_2(\eps)\leq C_3\,\e^{-\frac{\gamma t}{\eps}},
\end{align*}
where $C_3>0$ is a constant independ of $t\in[0,T]$.

\vspace{1mm}
\noindent
{\it iii) (Control of the third term)} To control the last term, let $V(s,m;t)$ be defined by (\ref{V}), and for $T\geq s$, define
\begin{align*}
\tilde V(s,m;t)=V(T-s,m;t).
\end{align*}
Then we have
\begin{align*}
&\tilde V(T,m;t)=V(0,m;t)=\psi(\pi_2^*m)-\psi\big(\tilde{\hat \zeta}^{t,\pi_1^*m}\big),\\
&\tilde V(0,m;0)=V(T,m;0)=\psi\big(\cL_{\hat Y^{0,\xi}_T(\eta)}\big)-\psi\big(\tilde{\hat \zeta}^{0,\mu}\big).
\end{align*}
As a result, taking $T=t/\eps$ we arrive at
\begin{align*}
\sJ_3(\eps)=\Big|\tilde V\Big(\frac{t}{\eps},\cL_{(\hat X_t^\eps,\hat Y_t^\eps)};t\Big)-\tilde V(0,\cL_{(\xi,\eta)};0)\Big|.
\end{align*}
Due to the low regularity of the function $\tilde V$ with respect to $m$ and $t$ variables, we define the mollifying approximation $\tilde V_n$ by (\ref{fnn}), i.e.,
\begin{align*}
\tilde V_n(s,m;t)=\tilde V(s,m*\rho_2^n;\cdot)*\rho_1^n,
\end{align*}
where
\begin{align*}
m*\rho_2^n(\cdot):=\int_\cdot\int_{\mR^{d_1}}\rho_2^n(x-\tilde x)m(\dif\tilde x,\dif y)\dif x.
\end{align*}
Then by Lemma \ref{fn} we have
\begin{align}\label{3n}
\sJ_3(\eps)&\leq\Big|\tilde V_n\Big(\frac{t}{\eps},\cL_{(\hat X_t^\eps,\hat Y_t^\eps)};t\Big)-\tilde V\Big(\frac{t}{\eps},\cL_{(\hat X_t^\eps,\hat Y_t^\eps)};t\Big)\Big|+\big|\tilde V_n(0,\cL_{(\xi,\eta)};0)-\tilde V(0,\cL_{(\xi,\eta)};0)\big|\no\\
&\quad+\Big|\tilde V_n\Big(\frac{t}{\eps},\cL_{(\hat X_t^\eps,\hat Y_t^\eps)};t\Big)-\tilde V_n(0,\cL_{(\xi,\eta)};0)\Big|\no\\
&\leq C_3\, n^{-\a}+\Big|\tilde V_n\Big(\frac{t}{\eps},\cL_{(\hat X_t^\eps,\hat Y_t^\eps)};t\Big)-\tilde V_n(0,\cL_{(\xi,\eta)};0)\Big|.
\end{align}
Using It\^{o}'s formula for the system (\ref{sde22}), we deduce that
\begin{align*}
{\tilde\sJ}_3(\eps):&=\Big|\tilde V_n\Big(\frac{t}{\eps},\cL_{(\hat X_t^\eps,\hat Y_t^\eps)};t\Big)-\tilde V_n(0,\cL_{(\xi,\eta)};0)\Big|\\
&=\bigg|\frac{1}{\eps}\int_0^t\p_s\tilde V_n\Big(\frac{r}{\eps},\cL_{(\hat X_r^\eps,\hat Y_r^\eps)};r\Big)\dif r+\int_0^t\p_t\tilde V_n\Big(\frac{r}{\eps},\cL_{(\hat X_r^\eps,\hat Y_r^\eps)};r\Big)\dif r\\
&\quad+\frac{1}{\eps}\int_0^t\mE\Big[F_\eps(r,\hat X_r^\eps,\hat Y_r^\eps)\cdot\p_y\frac{\delta\tilde V_n}{\delta m}\Big(\frac{r}{\eps},\cL_{(\hat X_r^\eps,\hat Y_r^\eps)};r\Big)(\hat X_r^\eps,\hat Y_r^\eps)\\
&\qquad\qquad\quad+\frac{1}{2}\T\Big(\cG_\eps(r,\hat X_r^\eps,\hat Y_r^\eps)\cdot\p_y^2\frac{\delta\tilde V_n}{\delta m}\Big(\frac{r}{\eps},\cL_{(\hat X_r^\eps,\hat Y_r^\eps)};r\Big)(\hat X_r^\eps,\hat Y_r^\eps)\Big)\Big]\dif r\\
&\quad+\int_0^t\mE\Big[b_\eps(r,\hat X_r^\eps,\hat Y_r^\eps)\cdot\p_x\frac{\delta\tilde V_n}{\delta m}\Big(\frac{r}{\eps},\cL_{(\hat X_r^\eps,\hat Y_r^\eps)};r\Big)(\hat X_r^\eps,\hat Y_r^\eps)\\
&\qquad\qquad\quad+\frac{1}{2}\T\Big(a_\eps(r,\hat X_r^\eps)\cdot\p_x^2\frac{\delta\tilde V_n}{\delta m}\Big(\frac{r}{\eps},\cL_{(\hat X_r^\eps,\hat Y_r^\eps)};r\Big)(\hat X_r^\eps,\hat Y_r^\eps)\Big)\Big]\dif r\bigg|.
\end{align*}
This together with the equation (\ref{eqV}) implies that
\begin{align*}
{\tilde\sJ}_3(\eps)
&\leq\bigg|\frac{1}{\eps}\int_0^t\Big[\p_s\tilde V_n\Big(\frac{r}{\eps},\cL_{(\hat X_r^\eps,\hat Y_r^\eps)};r\Big)-\p_s\tilde V\Big(\frac{r}{\eps},\cL_{(\hat X_r^\eps,\hat Y_r^\eps)};r\Big)\Big]\dif r\bigg|\\
&\quad+\bigg|\frac{1}{\eps}\int_0^t\mE\Big[F_\eps(r,\hat X_r^\eps,\hat Y_r^\eps)\cdot\Big(\p_y\frac{\delta\tilde V_n}{\delta m}\Big(\frac{r}{\eps},\cL_{(\hat X_r^\eps,\hat Y_r^\eps)};r\Big)(\hat X_r^\eps,\hat Y_r^\eps)\\
&\qquad\qquad\qquad\qquad-\p_y\frac{\delta\tilde V}{\delta m}\Big(\frac{r}{\eps},\cL_{(\hat X_r^\eps,\hat Y_r^\eps)};r\Big)(\hat X_r^\eps,\hat Y_r^\eps)\Big)\Big]\dif r\bigg|\\
&\quad+\bigg|\frac{1}{2\eps}\int_0^t\mE\Big[\T\Big(\cG_\eps(r,\hat X_r^\eps,\hat Y_r^\eps)\cdot\Big(\p_y^2\frac{\delta\tilde V_n}{\delta m}\Big(\frac{r}{\eps},\cL_{(\hat X_r^\eps,\hat Y_r^\eps)};r\Big)(\hat X_r^\eps,\hat Y_r^\eps)\\
&\qquad\qquad\qquad\qquad-\p_y^2\frac{\delta\tilde V}{\delta m}\Big(\frac{r}{\eps},\cL_{(\hat X_r^\eps,\hat Y_r^\eps)};r\Big)(\hat X_r^\eps,\hat Y_r^\eps)\Big)\Big)\Big]\dif r\bigg|\\
&\quad+\bigg|\frac{1}{\eps}\int_0^t\mE\Big[\big[F_\eps(r,\hat X_r^\eps,\hat Y_r^\eps)-\hat F(r,\hat X_r^\eps,\hat Y_r^\eps)\big]\cdot\p_y\frac{\delta\tilde V}{\delta m}\Big(\frac{r}{\eps},\cL_{(\hat X_r^\eps,\hat Y_r^\eps)};r\Big)(\hat X_r^\eps,\hat Y_r^\eps)\Big]\dif r\bigg|\\
&\quad+\bigg|\frac{1}{2\eps}\int_0^t\mE\Big[\T\Big(\big[\cG_\eps(r,\hat X_r^\eps,\hat Y_r^\eps)-\hat\cG(r,\hat X_r^\eps,\hat Y_r^\eps)\big]\\
&\qquad\qquad\qquad\qquad\qquad\cdot\p_y^2\frac{\delta\tilde V}{\delta m}\Big(\frac{r}{\eps},\cL_{(\hat X_r^\eps,\hat Y_r^\eps)};r\Big)(\hat X_r^\eps,\hat Y_r^\eps)\Big)\Big]\dif r\bigg|\\
&\quad+\bigg|\int_0^t\mE\Big[b_\eps(r,\hat X_r^\eps,\hat Y_r^\eps)\cdot\p_x\frac{\delta\tilde V_n}{\delta m}\Big(\frac{r}{\eps},\cL_{(\hat X_r^\eps,\hat Y_r^\eps)};r\Big)(\hat X_r^\eps,\hat Y_r^\eps)\\
&\qquad\qquad\quad+\frac{1}{2}\T\Big(a_\eps(r,\hat X_r^\eps)\cdot\p_x^2\frac{\delta\tilde V_n}{\delta m}\Big(\frac{r}{\eps},\cL_{(\hat X_r^\eps,\hat Y_r^\eps)};r\Big)(\hat X_r^\eps,\hat Y_r^\eps)\Big)\Big]\dif r\bigg|\\
&\quad+\bigg|\int_0^t\p_t\tilde V_n\Big(\frac{r}{\eps},\cL_{(\hat X_r^\eps,\hat Y_r^\eps)};r\Big)\dif r\bigg|=:\sum_{i=1}^7{\tilde\sJ}_{3,i}(\eps).
\end{align*}
For the first term, we have by the second estimate in (\ref{vt}) and estimate (\ref{fn1}) in Lemma \ref{fn} that
\begin{align*}
{\tilde\sJ}_{3,1}(\eps)&\leq C_4\, n^{-\a}\cdot\frac{1}{\eps}\int_0^t\Big(\frac{t-r}{\eps}\Big)^{\frac{\b}{2}-1} \cdot\e^{-\gamma\frac{t-r}{\eps}}\dif r\leq C_4\, n^{-\a}.
\end{align*}
Similarly, by the last two estimates in (\ref{vx}) and (\ref{vt}), and using estimate (\ref{fn1}) in Lemma \ref{fn} again, we have
\begin{align*}
{\tilde\sJ}_{3,2}(\eps)&\leq C_5\, \frac{1}{\eps}\int_0^t\mE\bigg[\Big|\p_y\frac{\delta\tilde V_n}{\delta m}\Big(\frac{r}{\eps},\cL_{(\hat X_r^\eps,\hat Y_r^\eps)};r\Big)(\hat X_r^\eps,\hat Y_r^\eps)\\
&\qquad\qquad\qquad-\p_y\frac{\delta\tilde V}{\delta m}\Big(\frac{r}{\eps},\cL_{(\hat X_r^\eps,\hat Y_r^\eps)};r\Big)(\hat X_r^\eps,\hat Y_r^\eps)\Big|\cdot\big(1+|\hat Y_r^\eps|^p\big)\bigg]\dif r\\
&\leq C_5\,n^{-\a}\cdot\frac{1}{\eps}\int_0^t\Big(\frac{t-r}{\eps}\Big)^{\frac{\b-1}{2}} \cdot\e^{-\gamma\frac{t-r}{\eps}}\dif r\leq C_5\, n^{-\a},
\end{align*}
and
\begin{align*}
{\tilde\sJ}_{3,3}(\eps)&\leq C_6\,\frac{1}{\eps}\int_0^t\mE\bigg[\Big|\p_y^2\frac{\delta\tilde V_n}{\delta m}\Big(\frac{r}{\eps},\cL_{(\hat X_r^\eps,\hat Y_r^\eps)};r\Big)(\hat X_r^\eps,\hat Y_r^\eps)\\
&\qquad\qquad\qquad-\p_y^2\frac{\delta\tilde V}{\delta m}\Big(\frac{r}{\eps},\cL_{(\hat X_r^\eps,\hat Y_r^\eps)};r\Big)(\hat X_r^\eps,\hat Y_r^\eps)\Big|\big(1+|\hat Y_r^\eps|^{2p}\big)\bigg]\dif r\\
&\leq C_6\,n^{-\a}\cdot\frac{1}{\eps}\int_0^t\Big(\frac{t-r}{\eps}\Big)^{\frac{\b}{2}-1} \cdot\e^{-\gamma\frac{t-r}{\eps}}\dif r\leq C_6\, n^{-\a}.
\end{align*}
To control the forth and fifth terms, we use the estimates (\ref{V2}) and (\ref{V3}) to deduce that
\begin{align*}
{\tilde\sJ}_{3,4}(\eps)&\leq C_7\,\frac{1}{\eps}\int_0^t\big\|F_\eps(r,\cdot,\cdot)-\hat F(r,\cdot,\cdot)\big\|_{L_p^\infty}\\
&\qquad\qquad\quad\cdot\mE\Big[\Big|\p_y\frac{\delta\tilde V}{\delta m}\Big(\frac{r}{\eps},\cL_{(\hat X_r^\eps,\hat Y_r^\eps)};r\Big)(\hat X_r^\eps,\hat Y_r^\eps)\Big|\big(1+|\hat Y_r^\eps|^p\big)\Big]\dif r\\
&\leq C_7\,\frac{1}{\eps}\int_0^t\|F_\eps(r,\cdot,\cdot)-\hat F(r,\cdot,\cdot)\|_{L_p^\infty}\cdot\Big(\frac{t-r}{\eps}\Big)^{\frac{\b-1}{2}}\cdot e^{-\frac{\gamma(t-r)}{\eps}}\dif r,
\end{align*}
and
\begin{align*}
{\tilde\sJ}_{3,5}(\eps)&\leq C_8\,\frac{1}{\eps}\int_0^t\big\|G_\eps(r,\cdot,\cdot)-\hat G(r,\cdot,\cdot)\big\|_{L_p^\infty}\\
&\qquad\qquad\quad\cdot\mE\Big[\Big|\p^2_y\frac{\delta\tilde V}{\delta m}\Big(\frac{r}{\eps},\cL_{(\hat X_r^\eps,\hat Y_r^\eps)};r\Big)(\hat X_r^\eps,\hat Y_r^\eps)\Big|\big(1+|\hat Y_r^\eps|^{2p}\big)\Big]\dif r\\
&\leq C_8\,\frac{1}{\eps}\int_0^t\|G_\eps(r,\cdot,\cdot)-\hat G(r,\cdot,\cdot)\|_{L_p^\infty}\cdot\Big(\frac{t-r}{\eps}\Big)^{\frac{\b}{2}-1}\cdot e^{-\frac{\gamma(t-r)}{\eps}}\dif r.
\end{align*}
Finally, by the first estimate in (\ref{vx}), the first estimate in (\ref{vt}) and using estimate (\ref{fn3}) in Lemma \ref{fn}, we have
\begin{align*}
{\tilde\sJ}_{3,6}(\eps)+{\tilde\sJ}_{3,7}(\eps)&\leq C_9\,n^{2-\a}\int_0^t\e^{-\frac{t-r}{\eps}}\dif r\leq C_9\,\eps\,n^{2-\a},
\end{align*}
In addition, one can check that the above constants $C_i, i=4,\cdots,9$ are independent of $t\in[0,T]$. Taking the above computations back into (\ref{3n}), we arrive at
\begin{align*}
\sJ_3(\eps)&\leq \hat C_3\bigg(n^{-\a}+\eps\,n^{2-\a}+\frac{1}{\eps}\int_0^t\|F_\eps(r,\cdot,\cdot)-\hat F(r,\cdot,\cdot)\|_{L_p^\infty}\cdot\Big(\frac{t-r}{\eps}\Big)^{\frac{\b-1}{2}}\cdot e^{-\frac{\gamma(t-r)}{\eps}}\dif r\\
&\quad+\frac{1}{\eps}\int_0^t\|G_\eps(r,\cdot,\cdot)-\hat G(r,\cdot,\cdot)\|_{L_p^\infty}\cdot\Big(\frac{t-r}{\eps}\Big)^{\frac{\b}{2}-1}\cdot e^{-\frac{\gamma(t-r)}{\eps}}\dif r\bigg).
\end{align*}
Choosing $n=\eps^{-1/2}$ we obtain
\begin{align*}
\sJ_3(\eps)&\leq \hat C_3\bigg(\eps^{\frac{\a}{2}}+\frac{1}{\eps}\int_0^t\|F_\eps(r,\cdot,\cdot)-\hat F(r,\cdot,\cdot)\|_{L_p^\infty}\cdot\Big(\frac{t-r}{\eps}\Big)^{\frac{\b-1}{2}}\cdot e^{-\frac{\gamma(t-r)}{\eps}}\dif r\\
&\quad+ \frac{1}{\eps}\int_0^t\|G_\eps(r,\cdot,\cdot)-\hat G(r,\cdot,\cdot)\|_{L_p^\infty}\cdot\Big(\frac{t-r}{\eps}\Big)^{\frac{\b}{2}-1}\cdot e^{-\frac{\gamma(t-r)}{\eps}}\dif r\bigg).
\end{align*}
Combining the estimates in {\it i)-iii)}, the proof is finished.
\end{proof}

\section{Multi-scale non-linear stochastic systems}

Throughout this section, we assume that the assumptions in Theorem \ref{main1} hold. We shall first study the asymptotic behavior of the non-autonomous approximation systems (\ref{sde4}) in subsection 6.1. Then, we seek the limit for the averaged systems  to derive the limits for the distributions of the slow process and fast motion of the non-linear system (\ref{sde1}) in subsection 6.2. Finally, we show in subsection 6.3 that the strong convergence in the averaging principle for the system (\ref{sde1}) follows directly from the convergence of the distributions of the slow and fast processes.

\subsection{The non-autonomous approximation systems}
To study the asymptotic behavior for the non-linear system (\ref{sde1}), we consider the non-autonomous approximation systems in (\ref{sden}). Namely,
for every $t\geq 0$ and $\eps>0$,
$$
X_t^{0,\eps}:= \xi\quad\text{and}\quad Y_t^{0,\eps}:= \eta,
$$
and for $n\geq 1$,
\begin{equation} \label{sden0}
\left\{ \begin{aligned}
&\dif X^{n,\eps}_t =b\big(X^{n,\eps}_t,\cL_{X_t^{n-1,\eps}},Y^{n,\eps}_t,\cL_{Y_t^{n-1,\eps}}\big)\dif t\\
&\qquad\qquad\quad+\sigma\big(X^{n,\eps}_t,\cL_{X_t^{n-1,\eps}},\cL_{Y_t^{n-1,\eps}}\big)\dif W^1_t, \qquad\quad\quad\quad X^{n,\eps}_0=\xi,\\
&\dif Y^{n,\eps}_t =\frac{1}{\eps}F\big(X^{n,\eps}_t,\cL_{X_t^{n-1,\eps}},Y^{n,\eps}_t,\cL_{Y_t^{n-1,\eps}}\big)\dif t\\
&\qquad\qquad\quad+\frac{1}{\sqrt{\eps}}G\big(X^{n,\eps}_t,\cL_{X_t^{n-1,\eps}},Y^{n,\eps}_t,\cL_{Y_t^{n-1,\eps}}\big)\dif W_t^2,\quad Y^{n,\eps}_0=\eta.
\end{aligned} \right.
\end{equation}
Each approximation system is a linear one
since the distributions appearing in the coefficients are not $\cL_{X_t^{n,\eps}}$ and $\cL_{Y_t^{n,\eps}}$ but rather $\cL_{X_t^{n-1,\eps}}$ and $\cL_{Y_t^{n-1,\eps}}$, i.e., the distributions of
the solutions of the approximations at the previous step. As stated in Section 2, for every $n\geq 1$,
we define
$$
b_{n,\eps}(t,x,y):=b\big(x,\cL_{X_t^{n-1,\eps}},y,\cL_{Y_t^{n-1,\eps}}\big), \quad \sigma_{n,\eps}(t,x):=\sigma\big(x,\cL_{X_t^{n-1,\eps}},\cL_{Y_t^{n-1,\eps}}\big)
$$
and
$$
F_{n,\eps}(t,x,y):=F\big(x,\cL_{X_t^{n-1,\eps}},y,\cL_{Y_t^{n-1,\eps}}\big),\quad G_{n,\eps}(t,x,y):=G\big(x,\cL_{X_t^{n-1,\eps}},y,\cL_{Y_t^{n-1,\eps}}\big).
$$
Then, each of the approximation system (\ref{sden0}) can be viewed as a non-autonomous system of the form (\ref{sde22}), i.e.,
\begin{equation} \label{sde40}
\left\{ \begin{aligned}
&\dif X^{n,\eps}_t =b_{n,\eps}(t,X^{n,\eps}_t,Y^{n,\eps}_t)\dif t
+\sigma_{n,\eps}(t,X^{n,\eps}_t)\dif W^1_t,\qquad\qquad\qquad X^{n,\eps}_0=\xi,\\
&\dif Y^{n,\eps}_t =\frac{1}{\eps}F_{n,\eps}(t,X^{n,\eps}_t,Y^{n,\eps}_t)\dif t+\frac{1}{\sqrt{\eps}}G_{n,\eps}(t,X^{n,\eps}_t,Y^{n,\eps}_t)\dif W_t^2,\quad\,   Y^{n,\eps}_0=\eta.
\end{aligned} \right.
\end{equation}
Therefore, by Theorem \ref{non-aut}, we have the following result for the approximation systems (\ref{sde40}).
For consistency, for every $t\geq 0$ and $x\in\mR^{d_1}$, we denote
$$
\bar X_t^0\equiv\xi\quad\text{and}\quad\zeta_0^{t,x}:=\cL_\eta.
$$

\bt\label{weakn}
For every $n\geq 1$, let $(X_t^{n,\eps},Y_t^{n,\eps})$ satisfy the system (\ref{sde40}). Assume that {\bf ($\bf H_{1}$)}  and {\bf ($\bf H_{2}$)} hold, $b,F,G\in C_p^{\a,(2,\a),\b,(2,\b)}$ and $\sigma\in C_p^{\a,(2,\a),(2,\b)}$ with $0<\a,\b\leq2$. Then we have for every $t\in[0,T]$  and $\varphi\in C_b^{(2,\alpha)}(\sP_2(\mR^{d_1}))$,
\begin{align}\label{weakxn}
\big|\varphi(\cL_{X_t^{n,\eps}})-\varphi(\cL_{\bar X_t^{n}})\big|\leq C_T\,\eps^{\frac{\a}{2}},
\end{align}
and for every $\psi\in C_p^{(2,\beta)}(\sP_2(\mR^{d_2}))$,
\begin{align}\label{weakyn}
\Big|\psi(\cL_{Y_t^{n,\eps}})-\psi\Big(\int_{\mR^{d_1}}\zeta_{n}^{t,x}\cL_{\bar X_t^{n}}(\dif x)\Big)\Big|\leq \tilde C_T\,\eps^{\frac{\a}{2}}+C_0\,e^{-\frac{\gamma_0t}{\eps}},
\end{align}
where $C_T,\tilde C_T, C_0>0$ and $\gamma_0>0$ are constants independent of $n$ and $\eps$,  and $\bar X_t^n$ satisfies the following averaged equation:
\begin{align}\label{aven}
\dif \bar X_t^{n}=\bar b_n(t,\bar X_t^{n})\dif t+\sigma_n(t,\bar X_t^{n})\dif W_t^1,\quad \bar X_0^n=\xi,
\end{align}
where the averaged drift coefficient is given by
\begin{align}\label{barbn}
\bar b_n(t,x):=\int_{\mR^{d_2}}b_n(t,x,y)\zeta_n^{t,x}(\dif y),
\end{align}
and $b_n$ and $\sigma_n$ are defined recursively by
\begin{eqnarray}
\begin{split}\label{bn}
&b_n(t,x,y):=b\Big(x,\cL_{\bar X_t^{n-1}},y,\int_{\mR^{d_1}}\zeta_{n-1}^{t,\tilde x}\cL_{\bar X_t^{n-1}}(\dif\tilde x)\Big),\\
&\sigma_n(t,x):=\sigma\Big(x,\cL_{\bar X_t^{n-1}},\int_{\mR^{d_1}}\zeta_{n-1}^{t,\tilde x}\cL_{\bar X_t^{n-1}}(\dif\tilde x)\Big),
\end{split}
\end{eqnarray}
and $\zeta_n^{t,x}(\dif y)$ is the unique invariant measure of the frozen equation
\begin{align}\label{Yn}
\dif Y_s^{n,(t,x)}=F_n(t,x,Y_s^{n,(t,x)})\dif s+G_n(t,x,Y_s^{n,(t,x)})\dif \tilde W_s,
\end{align}
where $(t,x)$ are parameters, $\tilde W_s$ is a standard Brownian motion and $F_n$ and $G_n$ are defined recursively by
\begin{eqnarray}
\begin{split}\label{Fn}
&F_n(t,x,y):=F\Big(x,\cL_{\bar X_t^{n-1}},y,\int_{\mR^{d_1}}\zeta_{n-1}^{t,\tilde x}\cL_{\bar X_t^{n-1}}(\dif\tilde x)\Big),\\
&G_n(t,x,y):=G\Big(x,\cL_{\bar X_t^{n-1}},y,\int_{\mR^{d_1}}\zeta_{n-1}^{t,\tilde x}\cL_{\bar X_t^{n-1}}(\dif\tilde x)\Big).
\end{split}
\end{eqnarray}
\et

\begin{proof}
For clarity, we divide the proof into the following steps.

\vspace{1mm}
\noindent{\it i)} For $n=1$,  note that we have
\begin{align*}
b_{1,\eps}(t,x,y)=b(x,\cL_\xi,y,\cL_\eta)=b_1(t,x,y),\quad  \sigma_{1,\eps}(t,x)=\sigma(x,\cL_\xi,\cL_\eta)=\sigma_1(t,x),
\end{align*}
and
\begin{align*}
F_{1,\eps}(t,x,y)&=F(x,\cL_\xi,y,\cL_\eta)=F_1(t,x,y),\\ G_{1,\eps}(t,x,y)&=G(x,\cL_\xi,y,\cL_\eta)=G_1(t,x,y),
\end{align*}
which are independent of $t$ and $\eps$. Thus the system (\ref{sde40}) reduces to the following classical autonomous multi-scale SDE:
\begin{equation*}
\left\{ \begin{aligned}
&\dif X^{1,\eps}_t =b\big(X^{1,\eps}_t,\cL_{\xi},Y^{1,\eps}_t,\cL_{\eta}\big)\dif t
+\sigma\big(X^{1,\eps}_t,\cL_{\xi},\cL_{\eta}\big)\dif W^1_t,\qquad\qquad\quad\,\,\,\, X^{1,\eps}_0=\xi,\\
&\dif Y^{1,\eps}_t =\frac{1}{\eps}F\big(X^{1,\eps}_t,\cL_{\xi},Y^{1,\eps}_t,\cL_{\eta}\big)\dif t+\frac{1}{\sqrt{\eps}}G\big(X^{1,\eps}_t,\cL_{\xi},Y^{1,\eps}_t,\cL_{\eta}\big)\dif W_t^2,\quad   Y^{1,\eps}_0=\eta.
\end{aligned} \right.
\end{equation*}
Since the coefficients satisfy the conditions in Theorem \ref{non-aut}, as a direct consequence we have that
\begin{align}\label{n1}
\begin{split}
&\big|\varphi(\cL_{X_t^{1,\eps}})-\varphi(\cL_{\bar X_t^{1}})\big|\leq C_2\,\eps^{\frac{\a}{2}}=:\mC_1(t,\eps),\\
&\Big|\psi(\cL_{Y_t^{1,\eps}})-\psi\Big(\int_{\mR^{d_1}}\zeta_{1}^{t,x}\cL_{\bar X_t^{1}}(\dif x)\Big)\Big|\leq C_2\,\eps^{\frac{\a}{2}}+C_3\,e^{-\frac{\gamma t}{\eps}}=:\mC_1(t,\eps)+\hat \mC_1(t,\eps),
\end{split}
\end{align}
where $\bar X_t^1$ satisfies the following averaged equation:
\begin{align}\label{aven1}
\dif \bar X_t^{1}=\bar b_1(\bar X_t^{1})\dif t+\sigma_1(\bar X_t^{1})\dif W_t^1,\quad \bar X_0^1=\xi,
\end{align}
with
\begin{align*}
\bar b_1(x)=\int_{\mR^{d_2}}b(x,\cL_\xi,y,\cL_\eta)\zeta_1^{x}(\dif y),\quad \sigma_1(x)=\sigma(x,\cL_\xi,\cL_\eta),
\end{align*}
and $\zeta_1^{x}(\dif y)$ is the unique invariant measure of the frozen equation
\begin{align}\label{Yn1}
\dif Y_s^{1,x}=F(x,\cL_\xi,Y_s^{1,x},\cL_\eta)\dif s+G(x,\cL_\xi,Y_s^{1,x},\cL_\eta)\dif\tilde W_s.
\end{align}
We remark that the constants $C_2,  C_3$ and $\gamma>0$ can be taken as the same as those in Theorem \ref{non-aut}.

\vspace{1mm}
\noindent {\it ii)} For $n=2$, by the definitions (\ref{bn}) and (\ref{Fn}), we have
\begin{align*}
b_2(t,x,y)&=b\Big(x,\cL_{\bar X_t^{1}},y,\int_{\mR^{d_1}}\zeta_{1}^{\tilde x}\cL_{\bar X_t^{1}}(\dif\tilde x)\Big),\\
\sigma_2(t,x)&=\sigma\Big(x,\cL_{\bar X_t^{1}},\int_{\mR^{d_1}}\zeta_{1}^{\tilde x}\cL_{\bar X_t^{1}}(\dif\tilde x)\Big),\\
F_2(t,x,y)&=F\Big(x,\cL_{\bar X_t^{1}},y,\int_{\mR^{d_1}}\zeta_{1}^{\tilde x}\cL_{\bar X_t^{1}}(\dif\tilde x)\Big),\\
G_2(t,x,y)&=G\Big(x,\cL_{\bar X_t^{1}},y,\int_{\mR^{d_1}}\zeta_{1}^{\tilde x}\cL_{\bar X_t^{1}}(\dif\tilde x)\Big).
\end{align*}
By the assumptions on the coefficients,  taking $\varphi(\cdot)=b(x,\cdot,y,\nu)$ and $\psi(\cdot)=b(x,\mu,y,\cdot)$ in (\ref{n1}), and  arguing as in the proof of (\ref{n1}) we obtain  that
\begin{align*}
|b_{2,\eps}(t,x,y)-b_2(t,x,y)|&=\Big|b\big(x,\cL_{X_t^{1,\eps}},y,\cL_{Y_t^{1,\eps}}\big)-b\Big(x,\cL_{\bar X_t^{1}},y,\int_{\mR^{d_1}}\zeta_{1}^{\tilde x}\cL_{\bar X_t^{1}}(\dif\tilde x)\Big)\Big|\\
&\leq C_0(1+|y|^p)\big[2\mC_1(t,\eps)+\hat \mC_1(t,\eps)\big].
\end{align*}
This in turn implies that
\begin{align*}
&\|b_{2,\eps}(t,\cdot,\cdot)-b_2(t,\cdot,\cdot)\|_{L_p^\infty}\leq C_0\big[2\mC_1(t,\eps)+\hat \mC_1(t,\eps)\big].
\end{align*}
Similarly, we can deduce that
\begin{align*}
&\|\sigma_{2,\eps}(t,\cdot)-\sigma_2(t,\cdot)\|_{L^\infty}\leq C_0\big[2\mC_1(t,\eps)+\hat \mC_1(t,\eps)\big],\\
&\|F_{2,\eps}(t,\cdot,\cdot)-F_2(t,\cdot,\cdot)\|_{L_p^\infty}\leq C_0\big[2\mC_1(t,\eps)+\kappa\,\hat \mC_1(t,\eps)\big],\\
&\|G_{2,\eps}(t,\cdot,\cdot)-G_2(t,\cdot,\cdot)\|_{L_p^\infty}\leq C_0\big[2\mC_1(t,\eps)+\kappa\,\hat \mC_1(t,\eps)\big],
\end{align*}
where $\kappa$ is the constant in (\ref{dissf2}) (see Remark \ref{resde} (iv)). Furthermore, we can derive that $b_2(t,\cdot,\cdot)$, $F_2(t,\cdot,\cdot)$, $G_2(t,\cdot,\cdot)\in C_p^{\a,\b}$ and $\sigma_2(t,\cdot)\in C_b^\a$.
Meanwhile, by ({\ref{lin}}) we have that for every $t_1,t_2\geq0$,
\begin{align*}
&|b_{2}(t_1,x,y)-b_{2}(t_2,x,y)|\\
&=\Big|b\Big(x,\cL_{\bar X_{t_1}^{1}},y,\int_{\mR^{d_1}}\zeta_{1}^{\tilde x}\cL_{\bar X_{t_1}^{1}}(\dif\tilde x)\Big)-b\Big(x,\cL_{\bar X_{t_2}^{1}},y,\int_{\mR^{d_1}}\zeta_{1}^{\tilde x}\cL_{\bar X_{t_2}^{1}}(\dif\tilde x)\Big)\Big|\\
&\leq C_0(1+|y|^p)\cW_2(\cL_{\bar X_{t_1}^{1}},\cL_{\bar X_{t_2}^{1}})^\a\leq C_0(1+|y|^p)|t_1-t_2|^\frac{\a}{2},
\end{align*}
and
\begin{align*}
&|\sigma_{2}(t_1,x)-\sigma_{2}(t_2,x)|+|F_{2}(t_1,x,y)-F_{2}(t_2,x,y)|\\
&+|G_{2}(t_1,x,y)-G_{2}(t_2,x,y)|\leq C_0(1+|y|^p)|t_1-t_2|^\frac{\a}{2}.
\end{align*}
Consequently, using Theorem \ref{non-aut} again we have
\begin{align*}
\big|\varphi(\cL_{X_t^{2,\eps}})-\varphi(\cL_{\bar X_t^{2}})\big|&\leq C_2\,\eps^{\frac{\a}{2}}+C_2\int_0^t\,(t-s)^{\frac{\a}{2}-1}\cdot\Big[\|b_{2,\eps}(s,\cdot,\cdot)- b_2(s,\cdot,\cdot)\|_{L_p^\infty}\\
&\quad+\|\sigma_{2,\eps}(s,\cdot)-\sigma_2(s,\cdot)\|_{L^\infty}
+\|F_{2,\eps}(s,\cdot,\cdot)-F_2(s,\cdot,\cdot)\|_{L_p^\infty}\\
&\quad+\|G_{2,\eps}(s,\cdot,\cdot)- G_2(s,\cdot,\cdot)\|_{L_p^\infty}\Big]\dif s\\
&\leq C_2\,\eps^{\frac{\a}{2}}+2\,C_0C_2\int_0^t(t-s)^{\frac{\a}{2}-1}\cdot\big[2\mC_1(s,\eps)+\hat \mC_1(s,\eps)\big]\dif s\\
&\quad+2\,C_0C_2\int_0^t(t-s)^{\frac{\a}{2}-1}\cdot\big[2\mC_1(s,\eps)+\kappa\,\hat \mC_1(s,\eps)\big]\dif s\\
&=:\mC_2(t,\eps),
\end{align*}
and
\begin{align*}
&\Big|\psi(\cL_{Y_t^{2,\eps}})-\psi\Big(\int_{\mR^{d_1}}\zeta_2^{t,x}\cL_{\bar X_t^{2}}(\dif x)\Big)\Big|\\
&\leq \mC_2(t,\eps)+C_3\,e^{-\frac{\gamma t}{\eps}}+\frac{C_3}{\eps}\int_0^t\,\Big(\frac{t-s}{\eps}\Big)^{\frac{\b}{2}-1}\cdot
e^{-\frac{\gamma(t-s)}{\eps}}\Big[\|F_{2,\eps}(s,\cdot,\cdot)-F_2(s,\cdot,\cdot)\|_{L_p^\infty}\\
&\qquad\quad+\|G_{2,\eps}(s,\cdot,\cdot)- G_2(s,\cdot,\cdot)\|_{L_p^\infty}\Big]\dif s\\
&\leq \mC_2(t,\eps)+C_3\,e^{-\frac{\gamma t}{\eps}}+\frac{2\,C_0C_3}{\eps}\int_0^t\,\Big(\frac{t-s}{\eps}\Big)^{\frac{\b}{2}-1}\cdot
e^{-\frac{\gamma(t-s)}{\eps}}\cdot\big[2\mC_1(s,\eps)+\kappa\,\hat \mC_1(s,\eps)\big]\dif s\\
&=:\mC_2(t,\eps)+\hat \mC_2(t,\eps),
\end{align*}
where $\bar X_t^{2}$ satisfies the averaged equation (\ref{aven}) with $n=2$, and $\zeta_2^{t,x}$ is the invariant measure for the frozen equation (\ref{Yn}) with $n=2$.

\vspace{1mm}
\noindent{\it iii)}
For $n\geq2$, we assume that the following estimates hold:
\begin{align}\label{psi1}
\big|\varphi(\cL_{X_t^{n,\eps}})-\varphi(\cL_{\bar X_t^{n}})\big|\leq \mC_{n}(t,\eps),
\end{align}
and
\begin{align}\label{psi2}
\Big|\psi(\cL_{Y_t^{n,\eps}})-\psi\Big(\int_{\mR^{d_1}}\zeta_{n}^{t,x}\cL_{\bar X_t^{n}}(\dif x)\Big)\Big|\leq \mC_n(t,\eps)+\hat \mC_{n}(t,\eps),
\end{align}
where $\mC_{n}(t,\eps)$ and $\hat \mC_{n}(t,\eps)$ are defined recursively by
\begin{align}
&\mC_{n}(t,\eps):=C_2\,\eps^{\frac{\a}{2}}+2\,C_0C_2\int_0^t(t-s)^{\frac{\a}{2}-1}\cdot\big[2\mC_{n-1}(s,\eps)
+\hat \mC_{n-1}(s,\eps)\big]\dif s\no\\
&\qquad\qquad\quad+2\,C_0C_2\int_0^t(t-s)^{\frac{\a}{2}-1}\cdot\big[2\mC_{n-1}(s,\eps)
+\kappa\,\hat \mC_{n-1}(s,\eps)\big]\dif s,\label{cn1}\\
&\hat\mC_{n}(t,\eps):=C_3\,e^{-\frac{\gamma t}{\eps}}\no\\
&\qquad\qquad\quad+\frac{2\,C_0C_3}{\eps}\int_0^t\,\Big(\frac{t-s}{\eps}\Big)^{\frac{\b}{2}-1}\cdot
e^{-\frac{\gamma(t-s)}{\eps}}\cdot\big[2\mC_{n-1}(s,\eps)+\kappa\,\hat \mC_{n-1}(s,\eps)\big]\dif s.\label{cn2}
\end{align}
Then we can deduce that
\begin{align*}
&|b_{n+1,\eps}(t,x,y)-b_{n+1}(t,x,y)|\\
&=\Big|b\big(x,\cL_{X_t^{n,\eps}},y,\cL_{Y_t^{n,\eps}}\big)-b\Big(x,\cL_{\bar X_t^{n}},y,\int_{\mR^{d_1}}\zeta_{n}^{t,\tilde x}\cL_{\bar X_t^{n}}(\dif\tilde x)\Big)\Big|\\
&\leq C_0(1+|y|^p)\big[2\mC_{n}(t,\eps)+\hat \mC_{n}(t,\eps)\big],
\end{align*}
and
\begin{align*}
&|b_{n+1}(t_1,x,y)-b_{n+1}(t_2,x,y)|\\
&=\Big|b\Big(x,\cL_{\bar X_{t_1}^{n}},y,\int_{\mR^{d_1}}\zeta_{n}^{t_1,\tilde x}\cL_{\bar X_{t_1}^{n}}(\dif\tilde x)\Big)-b\Big(x,\cL_{\bar X_{t_2}^{n}},y,\int_{\mR^{d_1}}\zeta_{n}^{t_2,\tilde x}\cL_{\bar X_{t_2}^{n}}(\dif\tilde x)\Big)\Big|\\
&\leq C_0(1+|y|^p)|t_1-t_2|^\frac{\a}{2}.
\end{align*}
The same conclusions hold for the other coefficients. Using the similar argument as the  proof for $n=2$ and by Theorem \ref{non-aut}, we can get that
\begin{align*}
&\big|\varphi(\cL_{X_t^{n+1,\eps}})-\varphi(\cL_{\bar X_t^{n+1}})\big|
\leq C_2\,\eps^{\frac{\a}{2}}\\
&\quad+2\,C_0C_2\int_0^t(t-s)^{\frac{\a}{2}-1}\cdot\big[2\mC_{n}(s,\eps)+\hat \mC_{n}(s,\eps)\big]\dif s\\
&\quad+2\,C_0C_2\int_0^t(t-s)^{\frac{\a}{2}-1}\cdot\big[2\mC_{n}(s,\eps)+\kappa\,\hat \mC_{n}(s,\eps)\big]\dif s=\mC_{n+1}(t,\eps),
\end{align*}
and
\begin{align*}
&\Big|\psi(\cL_{Y_t^{n+1,\eps}})-\psi\Big(\int_{\mR^{d_1}}\zeta_{n+1}^{t,x}\cL_{\bar X_t^{n+1}}(\dif x)\Big)\Big|\\
&\leq \mC_{n+1}(t,\eps)+C_3\,e^{-\frac{\gamma t}{\eps}}+\frac{2\,C_0C_3}{\eps}\int_0^t\,\Big(\frac{t-s}{\eps}\Big)^{\frac{\b}{2}-1}\\
&\qquad\qquad\qquad\times e^{-\frac{\gamma(t-s)}{\eps}}\cdot\big[2\mC_{n}(s,\eps)+\kappa\,\hat \mC_{n}(s,\eps)\big]\dif s=\mC_{n+1}(t,\eps)+\hat \mC_{n+1}(t,\eps),
\end{align*}
which in turn imply that the estimates (\ref{psi1}) and (\ref{psi2}) hold for $n+1$.

\vspace{1mm}
\noindent{\it iv)} It remains to provide the uniform control of the constants $\mC_{n}(t,\eps)$ and $\hat \mC_{n}(t,\eps)$ in estimates (\ref{psi1}) and (\ref{psi2}) with respect to $n\geq1$.
We have  the following claim:

\begin{description}
  \item[{\bf Claim}] There exist a sufficiently small $t_0>0$ and  constants $\Sigma_0, \tilde\Sigma_0>0$ (depending on $t_0$) and $\Sigma_1>0$, $0<\gamma_0<\gamma$ (independent of $t_0$) such that for every $t\in[0,t_0]$ and $n\geq1$,
$$
\mC_{n}(t,\eps)\leq\Sigma_0\,\eps^{\frac{\a}{2}}
$$
and
$$
\hat \mC_{n}(t,\eps)\leq\tilde\Sigma_0\,\eps^{\frac{\a}{2}}+\Sigma_1\,e^{-\frac{\gamma_0t}{\eps}}.
$$
\end{description}
Once these two estimates are proved, we can get that the estimates (\ref{weakxn}) and (\ref{weakyn}) hold for  $t\in[0,t_0]$. Since the system (\ref{sde40}) are classical SDE and $t_0$ does not depend on the initial condition, we can repeat the above argument on the interval $[t_0,2t_0]$ and iterate up to any finite time interval $[0,T]$ as in \cite{V0}.  The proof can be finished.

\vspace{1mm}
{\bf Proof of the Claim:}
Let us define
$$
h_n(t,\eps):=2\mC_n(t,\eps)+\kappa\,\hat \mC_n(t,\eps).
$$
Pulsing (\ref{cn1}) and (\ref{cn2}),  we have
\begin{align*}
h_{n}(t,\eps)&=2\,C_2\,\eps^{\frac{\a}{2}}+\kappa\,C_3\,e^{-\frac{\gamma t}{\eps}}+\Big[\frac{4\,C_0C_2}{\kappa}+4\,C_0C_2\Big]\cdot\int_0^t(t-s)^{\frac{\a}{2}-1}\cdot h_{n-1}(s,\eps)\dif s\\
&\quad+\frac{2\,C_0\kappa\,C_3}{\eps}\int_0^t\,\Big(\frac{t-s}{\eps}\Big)^{\frac{\b}{2}-1}\cdot
e^{-\frac{\gamma(t-s)}{\eps}}\cdot h_{n-1}(s,\eps)\dif s.
\end{align*}
Taking $\kappa$ small enough such that $\kappa\,C_3<\gamma$ and by induction, we can deduce that there exist small $t_0>0$ and $0<\l_0<\gamma$ such that for every $n\geq1$,
\begin{align*}
h_{n}(t,\eps)\leq \hat\Sigma_0\,\eps^{\frac{\a}{2}}+\hat\Sigma_1\,e^{-\frac{\gamma_0t}{\eps}}.
\end{align*}
where $\hat\Sigma_1>0$ is a constant independent of $t_0$. Taking this back into the definitions (\ref{cn1}) and (\ref{cn2}), we deduce that
\begin{align*}
\mC_{n}(t,\eps)&\leq C_2\,\eps^{\frac{\a}{2}}+\Big[\frac{2\,C_0C_2}{\kappa}+2\,C_0C_2\Big]\cdot\int_0^t(t-s)^{\frac{\a}{2}-1}
\cdot\big[\hat\Sigma_0\,\eps^{\frac{\a}{2}}+\hat\Sigma_1\,e^{-\frac{\gamma_0s}{\eps}}\big]\dif s\\
&\leq C_2\,\eps^{\frac{\a}{2}}+\Big[\frac{2\,C_0C_2}{\kappa}+2\,C_0C_2\Big]\cdot\frac{2}{\a}\cdot\hat\Sigma_0\,t^{\frac{\a}{2}}\eps^{\frac{\a}{2}}
+\Big[\frac{2\,C_0C_2}{\kappa}+2\,C_0C_2\Big]\cdot\hat\Sigma_1\tilde C\eps^{\frac{\a}{2}},
\end{align*}
and
\begin{align*}
\hat \mC_{n}(t,\eps)&\leq C_3\,e^{-\frac{\gamma t}{\eps}}+\frac{2\,C_0C_3}{\eps}\int_0^t\,\Big(\frac{t-s}{\eps}\Big)^{\frac{\b}{2}-1}
e^{-\frac{\gamma(t-s)}{\eps}}
\big[\hat\Sigma_0\,\eps^{\frac{\a}{2}\wedge1}+\hat\Sigma_1\,e^{-\frac{\gamma_0t}{\eps}}\big]\dif s\\
&\leq C_3\,e^{-\frac{\gamma t}{\eps}}+2\,C_0C_3\hat\Sigma_0\gamma^{-\frac{\b}{2}}\Gamma\Big(\frac{\b}{2}\Big)\eps^{\frac{\a}{2}}
+2\,C_0C_3\hat\Sigma_1(\gamma-\gamma_0)^{-\frac{\b}{2}}\Gamma\Big(\frac{\b}{2}\Big)\eps^{-\frac{\gamma_0t}{\eps}},
\end{align*}
which in turn imply the desired results.
\end{proof}

It is not easy to seek the limit of the non-linear system (\ref{sde1}) directly from the averaged equation (\ref{aven}) and the frozen equation (\ref{Yn}).  Below, we provide an alternative form of formulation for the  averaged approximation systems, which shall be used  to derive the limit as $n\to\infty$.  For simplicity, given a sequence of $\{\mu_n\}_{n\geq 1}$, we define
$$
\zeta_0^{x,\mu_{-1}}=\cL_\eta\quad\text{and}\quad \zeta_1^{x,\mu_0}=\zeta_1^x,
$$
where $\zeta^{x}_1$ is the invariant measure of the system (\ref{Yn1}). Let
$\zeta_n^{x,\mu_{n-1}}(\dif y)$ be the unique invariant measure for the following frozen equation:
\begin{align}\label{Ynn}
&\dif Y_s^{n,(x,\mu_{n-1})}=F\Big(x,\mu_{n-1},Y_s^{n,(x,\mu_{n-1})},\int_{\mR^{d_1}}\zeta_{n-1}^{\tilde x,\mu_{n-2}}\mu_{n-1}(\dif\tilde x)\Big)\dif s\no\\
&\qquad\qquad\qquad+G\Big(x,\mu_{n-1},Y_s^{n,(x,\mu_{n-1})},\int_{\mR^{d_1}}\zeta_{n-1}^{\tilde x,\mu_{n-2}}\mu_{n-1}(\dif\tilde x)\Big)\dif \tilde W_s,
\end{align}
where $\tilde W_s$ is a standard Brownian motion, and for $x\in\mR^{d_1}$, define
\begin{align*}
&\bar b_n(x,\mu_{n-1}):=\int_{\mR^{d_2}}b\Big(x,\mu_{n-1},y,\int_{\mR^{d_1}}\zeta_{n-1}^{\tilde x,\mu_{n-2}}\mu_{n-1}(\dif\tilde x)\Big)\zeta_n^{x,\mu_{n-1}}(\dif y),\\
&\bar\sigma_n(x,\mu_{n-1}):=\sigma\Big(x,\mu_{n-1},\int_{\mR^{d_1}}\zeta_{n-1}^{\tilde x,\mu_{n-2}}\mu_{n-1}(\dif\tilde x)\Big).
\end{align*}
We have the following result.

\bc\label{weaknn}
For the  averaged approximation systems (\ref{aven}), we have for every $n\geq1$,
\begin{align}\label{avenn}
\dif \bar X_t^{n}=\bar b_n(\bar X_t^{n},\cL_{\bar X_t^{n-1}})\dif t+\bar\sigma_n(\bar X_t^{n},\cL_{\bar X_t^{n-1}})\dif W_t^1,\quad \bar X_0^n=\xi,
\end{align}
and for every $t\in[0,T]$ and $\psi\in C_p^{(2,\beta)}(\sP_2(\mR^{d_2}))$,
\begin{align}\label{weakynn}
\Big|\psi(\cL_{Y_t^{n,\eps}})-\psi\Big(\int_{\mR^{d_1}}\zeta_{n}^{x, \cL_{\bar X_t^{n-1}}}\cL_{\bar X_t^{n}}(\dif x)\Big)\Big|\leq \tilde C_T\,\eps^{\frac{\a}{2}}+C_0\,e^{-\frac{\gamma_0 t}{\eps}},
\end{align}
where $\zeta_n^{x,\mu_{n-1}}$ is the unique invariant measure of (\ref{Ynn})  with $\mu_{n-1}:=\cL_{\bar X_t^{n-1}}$.
\ec
\begin{proof}
For $n=1$, it is easy to see that the conclusions are true by (\ref{aven1}) and (\ref{Yn1}). Applying (\ref{bn}) and (\ref{Fn}) in Theorem \ref{weakn}, we have
\begin{align*}
&b_2(t,x,y)=b\Big(x,\cL_{\bar X_t^1},y,\int_{\mR^{d_1}}\zeta_{1}^{\tilde x}\cL_{\bar X_t^1}(\dif\tilde x)\Big)=b\Big(x,\mu_1,y,\int_{\mR^{d_1}}\zeta_{1}^{\tilde x,\mu_0}\mu_1(\dif\tilde x)\Big),\\
&\sigma_2(t,x)=\sigma\Big(x,\cL_{\bar X_t^1},\int_{\mR^{d_1}}\zeta_{1}^{\tilde x}\cL_{\bar X_t^1}(\dif\tilde x)\Big)=\sigma\Big(x,\mu_1,\int_{\mR^{d_1}}\zeta_{1}^{\tilde x,\mu_0}\mu_1(\dif\tilde x)\Big),\\
&F_2(t,x,y)=F\Big(x,\cL_{\bar X_t^1},y,\int_{\mR^{d_1}}\zeta_{1}^{\tilde x}\cL_{\bar X_t^1}(\dif\tilde x)\Big)=F\Big(x,\mu_1,y,\int_{\mR^{d_1}}\zeta_{1}^{\tilde x,\mu_0}\mu_1(\dif\tilde x)\Big),
\end{align*}
and
\begin{align*}
G_2(t,x,y)=G\Big(x,\cL_{\bar X_t^1},y,\int_{\mR^{d_1}}\zeta_{1}^{\tilde x}\cL_{\bar X_t^1}(\dif\tilde x)\Big)=G\Big(x,\mu_1,y,\int_{\mR^{d_1}}\zeta_{1}^{\tilde x,\mu_0}\mu_1(\dif\tilde x)\Big).
\end{align*}
Consequently, the frozen equation (\ref{Yn}) with $n=2$ can be rewritten as
\begin{align*}
&\dif Y_s^{2,(x,\mu_1)}=F\Big(x,\mu_1,Y_s^{2,(x,\mu_1)},\int_{\mR^{d_1}}\zeta_1^{\tilde x,\mu_0}\mu_1(\dif\tilde x)\Big)\dif s\\
&\qquad\qquad\quad+G\Big(x,\mu_1,Y_s^{2,(x,\mu_1)},\int_{\mR^{d_1}}\zeta_1^{\tilde x,\mu_0}\mu_1(\dif\tilde x)\Big)\dif \tilde W_s,
\end{align*}
and the corresponding invariant measure $\zeta_2^{t,x}$ equals to $\zeta_2^{x,\mu_1}$ with $\mu_1=\cL_{\bar X_t^1}$. These together with (\ref{barbn}) yield that
\begin{align*}
\bar b_2(t,x)=\int_{\mR^{d_2}}b\Big(x,\mu_1,y,\int_{\mR^{d_1}}\zeta_{1}^{\tilde x,\mu_0}\mu_1(\dif\tilde x)\Big)\zeta_2^{x,\mu_1}(\dif y)=\bar b_2(x,\mu_1),
\end{align*}
and
\begin{align*}
\sigma_2(t,x)=\sigma\Big(x,\mu_1,\int_{\mR^{d_1}}\zeta_{1}^{\tilde x,\mu_0}\mu_1(\dif\tilde x)\Big)=\bar\sigma_2(x,\mu_1),
\end{align*}
which in turn implies that the averaged equation (\ref{avenn}) with $n=2$ holds. In addition, by the estimate (\ref{weakyn}) we arrive at
\begin{align*}
\Big|\psi(\cL_{Y_t^{2,\eps}})-\psi\Big(\int_{\mR^{d_1}}\zeta_{2}^{x,\cL_{\bar X_t^1}}\cL_{\bar X_t^{2}}(\dif x)\Big)\Big|\leq\tilde C_T\,\eps^{\frac{\a}{2}}+C_0\,e^{-\frac{\l_0 t}{\eps}}.
\end{align*}
Assume that the averaged equation (\ref{avenn}), the estimate (\ref{weakynn}) and the frozen equation (\ref{Ynn}) hold for $n-1\geq2$, then we have
\begin{align*}
&b_n(t,x,y)=b\Big(x,\mu_{n-1},y,\int_{\mR^{d_1}}\zeta_{n-1}^{\tilde x,\mu_{n-2}}\mu_{n-1}(\dif\tilde x)\Big),\\
&\sigma_n(t,x)=\sigma\Big(x,\mu_{n-1},\int_{\mR^{d_1}}\zeta_{n-1}^{\tilde x,\mu_{n-2}}\mu_{n-1}(\dif\tilde x)\Big),\\
&F_n(t,x,y)=F\Big(x,\mu_{n-1},y,\int_{\mR^{d_1}}\zeta_{n-1}^{\tilde x,\mu_{n-2}}\mu_{n-1}(\dif\tilde x)\Big),\\
&G_n(t,x,y)=G\Big(x,\mu_{n-1},y,\int_{\mR^{d_1}}\zeta_{n-1}^{\tilde x,\mu_{n-2}}\mu_{n-1}(\dif\tilde x)\Big).
\end{align*}
Similarly, the frozen equation (\ref{Yn}) can be rewritten as
\begin{align*}
&\dif Y_s^{n,(x,\mu_{n-1})}=F\Big(x,\mu_{n-1},Y_s^{n,(x,\mu_{n-1})},\int_{\mR^{d_1}}\zeta_{n-1}^{\tilde x,\mu_{n-2}}\mu_{n-1}(\dif\tilde x)\Big)\dif s\\
&\qquad\qquad\qquad+G\Big(x,\mu_{n-1},Y_s^{n,(x,\mu_{n-1})},\int_{\mR^{d_1}}\zeta_{n-1}^{\tilde x,\mu_{n-2}}\mu_{n-1}(\dif\tilde x)\Big)\dif \tilde W_s,
\end{align*}
that is, (\ref{Ynn}) holds for given $n$. Moreover, the corresponding invariant measure $\zeta_n^{t,x}$ equals to $\zeta_n^{x,\mu_{n-1}}$, and
\begin{align*}
&\bar b_n(t,x)=\int_{\mR^{d_2}}b\Big(x,\mu_{n-1},y,\int_{\mR^{d_1}}\zeta_{n-1}^{\tilde x,\mu_{n-2}}\mu_{n-1}(\dif\tilde x)\Big)\zeta_n^{x,\mu_{n-1}}(\dif y)=\bar b_n(x,\mu_{n-1}),\\
&\sigma_n(t,x)=\sigma\Big(x,\mu_{n-1},\int_{\mR^{d_1}}\zeta_{n-1}^{\tilde x,\mu_{n-2}}\mu_{n-1}(\dif\tilde x)\Big)=\bar\sigma_n(x,\mu_{n-1}).
\end{align*}
As a result, the averaged equation (\ref{avenn}) and the estimate (\ref{weakynn}) hold for given $n$. Thus the proof is finished.
\end{proof}

\subsection{Characterization of the limits}

In this subsection, we  give the proof of the convergence of the distributions of the slow and fast process by using the approximation systems (\ref{sde40}).


We first provide the following result.

\bl
Assume that (\ref{dissf}) holds. Then we have that for any $q\geq 2$,
\begin{align}\label{tight}
\sup_{n\geq 1}\sup_{x\in\mR^{d_1}}\int_{\mR^{d_2}}|y|^{q}\zeta_n^{x,\mu_{n-1}}(\dif y)\leq C_0<\infty.
\end{align}
\el

\begin{proof}
For simplicity, we only prove the estimate (\ref{tight}) for $q=2$, the general case follows by the same argument. For any $n\geq1$, using It\^o's formula and (\ref{dissf}) we have
	\begin{equation*}
\begin{split}
		\dif \left|Y_s^{n,(x,\mu_{n-1})}\right|^2=& 2\left\langle Y_s^{n,(x,\mu_{n-1})},F\Big(x,\mu_{n-1},Y_s^{n,(x,\mu_{n-1})},\int_{\mR^{d_1}}\zeta_{n-1}^{\tilde x,\mu_{n-2}}\mu_{n-1}(\dif\tilde x)\Big)\right\rangle\dif s\\
&+\left\|G\Big(x,\mu_{n-1},Y_s^{n,(x,\mu_{n-1})},\int_{\mR^{d_1}}\zeta_{n-1}^{\tilde x,\mu_{n-2}}\mu_{n-1}(\dif\tilde x)\Big)\right\|^2\dif s\\
		&+2\left\langle Y_s^{n,(x,\mu_{n-1})},G\Big(x,\mu_{n-1},Y_s^{n,(x,\mu_{n-1})},\int_{\mR^{d_1}}\zeta_{n-1}^{\tilde x,\mu_{n-2}}\mu_{n-1}(\dif\tilde x)\Big)\dif W_s\right\rangle\\
		\leq&\left(-C_1\left|Y_s^{n,(x,\mu_{n-1})}\right|^{2}+C_2\left\|\int_{\mR^{d_1}}\zeta_{n-1}^{\tilde x,\mu_{n-2}}\mu_{n-1}(\dif\tilde x)\right\|_{2}^{2}+C_3\right)\dif s\\
&+2\left\langle Y_s^{n,(x,\mu_{n-1})},G\Big(x,\mu_{n-1},Y_s^{n,(x,\mu_{n-1})},\int_{\mR^{d_1}}\zeta_{n-1}^{\tilde x,\mu_{n-2}}\mu_{n-1}(\dif\tilde x)\Big)\dif W_s\right\rangle.
\end{split}
\end{equation*}
	Writing the above inequality in integral form and taking expectation, we get that
	\begin{align}\label{Ito}
      \begin{split}
			&\mathbb{E}\left[\left|Y_s^{n,(x,\mu_{n-1})}\right|^2\right]-\mathbb{E}|\eta|^2\\
&\leq -C_1\int_0^s\mathbb{E}\left[\left|Y_r^{n,(x,\mu_{n-1})}\right|^2\right]\dif r+C_2\left\|\int_{\mR^{d_1}}\zeta_{n-1}^{\tilde x,\mu_{n-2}}\mu_{n-1}(\dif\tilde x)\right\|_2^2\cdot s+C_3\cdot s.
      \end{split}
	\end{align}
Dividing both sides of \eqref{Ito} by $s$, letting $s\rightarrow\infty$ and using the ergodic theorem, we have that
 \begin{align*}
    	C_1\int_{\mR^{d_2}}|y|^{2}\zeta_n^{x,\mu_{n-1}}(\dif y)&\leq C_2\int_{\mR^{d_1}}\int_{\mR^{d_2}}|y|^{2}\zeta_{n-1}^{\tilde x,\mu_{n-2}}(\dif y)\mu_{n-1}(\dif \tilde x)+C_3\\
    &\leq C_2\sup_{x\in\mR^{d_1}}\int_{\mR^{d_2}}|y|^{2}\zeta_{n-1}^{x,\mu_{n-2}}(\dif y)+C_3.
 \end{align*}
As a result, it holds that
   \begin{align*}
    \begin{split}
    	\sup_{x\in\mR^{d_1}}\int_{\mR^{d_2}}|y|^{2}\zeta_n^{x,\mu_{n-1}}(\dif y)\leq \frac{C_2}{C_1}\left(\sup_{x\in\mR^{d_1}}\int_{\mR^{d_2}}|y|^{2}\zeta_{n-1}^{x,\mu_{n-2}}(\dif y)\right)+\frac{C_3}{C_1}.
    \end{split}
 \end{align*}
By induction, we deduce  that
\begin{align*}
		\sup_{x\in\mR^{d_1}}\int_{\mR^{d_2}}|y|^{2}\zeta_n^{x,\mu_{n-1}}(\dif y) &\leq \left(\frac{C_2}{C_1}\right)^n\left(\sup_{x\in\mR^{d_1}}\int_{\mR^{d_2}}|y|^{2}\zeta_0^{x,\mu_{-1}}(\dif y)\right) +\frac{C_3}{C_1}\left[\sum_{i=0}^{n-1}\left(\frac{C_2}{C_1}\right)^i\right]\\
&=\mE|\eta|^2\cdot\left(\frac{C_2}{C_1}\right)^n +\frac{C_3}{C_1}\left[\sum_{i=0}^{n-1}\left(\frac{C_2}{C_1}\right)^i\right].
\end{align*}
Therefore, under the assumption that $C_2<C_1$, we obtain
$$\sup_{n\geq 1}\sup_{x\in\mR^{d_1}}\int_{\mR^{d_2}}|y|^{2}\zeta_n^{x,\mu_{n-1}}(\dif y)<\infty.$$
The proof is finished.
\end{proof}

The above result implies the tightness of $\{\zeta_n^{x,\mu_{n-1}}\}_{n\geq 1}$. Below, we proceed to identify the limit of $\zeta_n^{x,\mu_{n-1}}$. For $M>0$ and two probability measures $\mu_1,\mu_2\in\sP_2(\mR^{d_1})$, we define
\begin{align*}
\rho_{\a,M}(\mu_1,\mu_2):=\sup_{\|\varphi\|_{C_b^{(2,\a)}}\leq M}\big|\varphi(\mu_1)-\varphi(\mu_2)\big|,
\end{align*}
and
\begin{align}\label{dis}
\tilde\rho_V(\zeta_n^{\cdot,\mu_{n-1}},\zeta^{\cdot,\mu}):=\sup_{x\in\mR^{d_1}}\rho_V(\zeta_n^{x,\mu_{n-1}},\zeta^{x,\mu}).
\end{align}
We have the following result.

\bl\label{inv}
Assume that $\mu_n$ converge weakly to $\mu$. Then there exists a $\zeta^{x,\mu}$ such that for every $\vartheta\geq 1$,
\begin{align}\label{con}
\lim_{n\to\infty}\sup_{x\in\mR^{d_1}}W_\vartheta (\zeta_n^{x,\mu_{n-1}},\zeta^{x,\mu})=0.
\end{align}
Moreover, $\zeta^{x,\mu}$ is the invariant measure of the following system:
\begin{align}\label{Yxmu}
&\dif Y_s^{x,\mu}=F\Big(x,\mu,Y_s^{x,\mu},\int_{\mR^{d_1}}\zeta^{\tilde x,\mu}\mu(\dif\tilde x)\Big)\dif s+G\Big(x,\mu,Y_s^{x,\mu},\int_{\mR^{d_1}}\zeta^{\tilde x,\mu}\mu(\dif\tilde x)\Big)\dif \tilde W_s.
\end{align}
By the uniqueness of the solution, $\zeta^{x,\mu}$ is also an invariant measure of the McKean-Vlasov equation (\ref{frozen0}),
and there exist constants $C_0, M>0$ such that
\begin{align}\label{ax}
\tilde \rho_V(\zeta_n^{\cdot,\mu_{n-1}},\zeta^{\cdot,\mu})\leq C_0\Big(\rho_{\a,M}(\mu_{n-1},\mu)+\kappa\,\tilde\rho_V(\zeta_{n-1}^{\cdot,\mu_{n-2}},\zeta^{\cdot,\mu})\Big),
\end{align}
where $\kappa$ is the constant in (\ref{dissf2}).
\el
\begin{proof}
The existence of $\zeta^{x,\mu}$ and the convergence in (\ref{con}) follows by the estimate (\ref{tight}).
It remains to show that  $\zeta^{x,\mu}$ is an invariant measure of the equation (\ref{Yxmu}). We deduce that for any $s>0$ and $g\in C_0^\infty(\mR^{d_2})$,
\begin{align*}
|\mE g(Y_s^{x,\mu})-\<g,\zeta^{x,\mu}\>|&\leq |\mE g(Y_s^{x,\mu})-\mE g(Y_s^{n,(x,\mu_{n-1})})|\\
&\quad+|\mE g(Y_s^{n,(x,\mu_{n-1})})-\<g,\zeta_n^{x,\mu_{n-1}}\>|+|\<g,\zeta_n^{x,\mu_{n-1}}\>-\<g,\zeta^{x,\mu}\>|.
\end{align*}
By (\ref{con}), we have that
$$
\lim_{n\to\infty}|\<g,\zeta_n^{x,\mu_{n-1}}\>-\<g,\zeta^{x,\mu}\>|=0.
$$
For the second term, by the assumption (\ref{dissf}) and the estimate (\ref{tight}), we have that
\begin{align*}
&2\left\langle y,F\Big(x,\mu_{n-1},y,\int_{\mR^{d_1}}\zeta_{n-1}^{\tilde x,\mu_{n-2}}\mu_{n-1}(\dif\tilde x)\Big)\right\rangle+\left\|G\Big(x,\mu_{n-1},y,\int_{\mR^{d_1}}\zeta_{n-1}^{\tilde x,\mu_{n-2}}\mu_{n-1}(\dif\tilde x)\Big)\right\|^2\\
&\leq -C_1|y|^2+C_2\left\|\int_{\mR^{d_1}}\zeta_{n-1}^{\tilde x,\mu_{n-2}}\mu_{n-1}(\dif\tilde x)\right\|_2^2+C_3\\
&\leq -C_1|y|^2+C_2\int_{\mR^{d_1}}\left\|\zeta_{n-1}^{\tilde x,\mu_{n-2}}\right\|_2^2\mu_{n-1}(\dif\tilde x)+C_3\leq-C_1|y|^2+\tilde C_3,
\end{align*}
where $\tilde C_3$ is independent of $n$.
Since $\zeta_n^{x,\mu_{n-1}}$ is the invariant measure for $Y_s^{n,(x,\mu_{n-1})}$, and the dissipative condition holds uniformly with respect to $n$, we get that
$$
|\mE g(Y_s^{n,(x,\mu_{n-1})})-\<g,\zeta_n^{x,\mu_{n-1}}\>|\leq C_0\,e^{-\gamma s},
$$
where $C_0, \gamma$ are constants independent of $n$. To control the first term, we use Lemma \ref{diff2} to deduce that
\begin{align*}
&\big|\mE g(Y_s^{x,\mu})-\mE g(Y_s^{n,(x,\mu_{n-1})})\big|\leq C_4\bigg(\Big\|F\Big(x,\mu,\cdot,\int_{\mR^{d_1}}\zeta^{\tilde x,\mu}\mu(\dif\tilde x)\Big)\\
&\quad-F\Big(x,\mu_{n-1},\cdot,\int_{\mR^{d_1}}\zeta_{n-1}^{\tilde x,\mu_{n-2}}\mu_{n-1}(\dif\tilde x)\Big)\Big\|_{L_p^\infty}+\Big\|G\Big(x,\mu,\cdot,\int_{\mR^{d_1}}\zeta^{\tilde x,\mu}\mu(\dif\tilde x)\Big)\\
&\quad-G\Big(x,\mu_{n-1},\cdot,\int_{\mR^{d_1}}\zeta_{n-1}^{\tilde x,\mu_{n-2}}\mu_{n-1}(\dif\tilde x)\Big)\Big\|_{L_p^\infty}\bigg).
\end{align*}
By the assumptions on the coefficients, we have
\begin{align*}
&\Big\|F\Big(x,\mu,\cdot,\int_{\mR^{d_1}}\zeta^{\tilde x,\mu}\mu(\dif\tilde x)\Big)-F\Big(x,\mu_{n-1},\cdot,\int_{\mR^{d_1}}\zeta_{n-1}^{\tilde x,\mu_{n-2}}\mu_{n-1}(\dif\tilde x)\Big)\Big\|_{L_p^\infty}\\
&\leq C_4\bigg(\rho_{\a,M}(\mu_{n-1},\mu)+\kappa\,\rho_V\bigg(\int_{\mR^{d_1}}\zeta^{\tilde x,\mu}\mu(\dif\tilde x),\int_{\mR^{d_1}}\zeta_{n-1}^{\tilde x,\mu_{n-2}}\mu_{n-1}(\dif\tilde x)\bigg)\bigg)\\
&\leq C_4\bigg(\rho_{\a,M}(\mu_{n-1},\mu)+\kappa\,\rho_V\bigg(\int_{\mR^{d_1}}\zeta^{\tilde x,\mu}\mu(\dif\tilde x),\int_{\mR^{d_1}}\zeta^{\tilde x,\mu}\mu_{n-1}(\dif\tilde x)\bigg)\\
&\qquad+\kappa\,\rho_V\bigg(\int_{\mR^{d_1}}\zeta^{\tilde x,\mu}\mu_{n-1}(\dif\tilde x),\int_{\mR^{d_1}}\zeta_{n-1}^{\tilde x,\mu_{n-2}}\mu_{n-1}(\dif\tilde x)\bigg)\bigg),
\end{align*}
where in the first inequality we used the assumption (\ref{dissf2}). On the one hand, by the definition (\ref{tv}) we have
\begin{align*}
&\rho_\cV\bigg(\int_{\mR^{d_1}}\zeta^{\tilde x,\mu}\mu(\dif\tilde x),\int_{\mR^{d_1}}\zeta^{\tilde x,\mu}\mu_{n-1}(\dif\tilde x)\bigg)\\
&=\sup_{\|f\|_{1+\cV}\leq 1}\int_{\mR^{d_1}}\int_{\mR^{d_2}}f(y)\zeta^{\tilde x,\mu}(\dif y)\Big(\mu(\dif\tilde x)-\mu_{n-1}(\dif\tilde x)\Big),
\end{align*}
and we have by Lemma \ref{lem3} below that
$$
\bar f(\tilde x)=\int_{\mR^{d_2}}f(y)\zeta^{\tilde x,\mu}(\dif y)\in C_b^{\a}.
$$
Thus we can choose $M$ large enough such that
\begin{align*}
\rho_\cV\bigg(\int_{\mR^{d_1}}\zeta^{\tilde x,\mu}\mu(\dif\tilde x),\int_{\mR^{d_1}}\zeta^{\tilde x,\mu}\mu_{n-1}(\dif\tilde x)\bigg)\leq \rho_{\a,M}(\mu_{n-1},\mu).
\end{align*}
On the other hand,  by the Minkowski inequality, we  have that
\begin{align*}
&\rho_V\bigg(\int_{\mR^{d_1}}\zeta^{\tilde x,\mu}\mu_{n-1}(\dif\tilde x),\int_{\mR^{d_1}}\zeta_{n-1}^{\tilde x,\mu_{n-2}}\mu_{n-1}(\dif\tilde x)\bigg)\\
&\leq \int_{\mR^{d_1}}\rho_V(\zeta^{\tilde x,\mu},\zeta_{n-1}^{\tilde x,\mu_{n-2}})\mu_{n-1}(\dif\tilde x)\leq \tilde\rho_V(\zeta^{\cdot,\mu},\zeta_{n-1}^{\cdot,\mu_{n-2}}).
\end{align*}
As a result, we get
\begin{align}\label{666}
&\Big\|F\Big(x,\mu,\cdot,\int_{\mR^{d_1}}\zeta^{\tilde x,\mu}\mu(\dif\tilde x)\Big)-F\Big(x,\mu_{n-1},\cdot,\int_{\mR^{d_1}}\zeta_{n-1}^{\tilde x,\mu_{n-2}}\mu_{n-1}(\dif\tilde x)\Big)\Big\|_{L_p^\infty}\no\\
&\leq C_4\Big(\rho_{\a,M}(\mu_{n-1},\mu)+\kappa\tilde\rho_V(\zeta^{\cdot,\mu},\zeta_{n-1}^{\cdot,\mu_{n-2}})\Big).
\end{align}
The same estimate holds for the coefficient $G$. We obtain that
\begin{align}\label{66}
\big|\mE g(Y_s^{x,\mu})-\mE g(Y_s^{n,(x,\mu_{n-1})})\big|
\leq C_4\Big(\rho_{\a,M}(\mu_{n-1},\mu)+\kappa\tilde\rho_V(\zeta^{\cdot,\mu},\zeta_{n-1}^{\cdot,\mu_{n-2}})\Big),
\end{align}
Letting $n\to\infty$, we obtain that
$$
|\mE g(Y_s^{x,\mu})-\<g,\zeta^{x,\mu}\>|\leq C_0\,e^{-\gamma s},
$$
which implies $\zeta^{x,\mu}$ is the invariant measure of $Y_s^{x,\mu}$. The estimate (\ref{ax}) can be proved similarly as (\ref{66}), we omit the details.
\end{proof}

We  need the following regularity result for the averaged coefficients.
\bl\label{lem3}
Assume that $f\in C_p^{\a,(2,\a),\beta,(2,\beta)}$ with $0<\alpha,\beta\leq 2$, and define
\begin{align*}
\bar f(x,\mu):=\int_{\mR^{d_2}} f\Big(x,\mu,y,\int_{\mR^{d_1}}\zeta^{\tilde x,\mu}\mu(\dif \tilde x)\Big)\zeta^{x,\mu}(\dif y),
\end{align*}
where $\zeta^{x,\mu}$ is the unique invariant measure of system (\ref{Yxmu}). Then we have $\bar f\in C_b^{\a, (2,\a)}$.
\el
\begin{proof}
We only prove the result for $1<\a,\beta\leq 1$, the case that $1<\a,\beta\leq 2$ can be proved similarly.	Let
	$$
	\tilde\zeta^\mu:=\int_{\mR^{d_1}}\zeta^{\tilde x,\mu}\mu(\dif\tilde x),
	$$
	and define
	$$
	\tilde F(x,\mu,y):=F(x,\mu,y,\tilde\zeta^\mu),\quad \tilde G(x,\mu,y):=G(x,\mu,y,\tilde\zeta^\mu).
	$$
	Then the system (\ref{Yxmu}) can be written as
	\begin{align*}
	&\dif Y_s^{x,\mu}=\tilde F(x,\mu,Y_s^{x,\mu})\dif s+\tilde G(x,\mu,Y_s^{x,\mu})\dif \tilde W_s,
	\end{align*}
	where $(x,\mu)$ are parameters. For every fixed $\mu$, it is easy to see that
	$\tilde F(\cdot,\mu,y), \tilde G(\cdot,\mu,y)\in C_b^\a$ and $\tilde F(x,\mu,\cdot), \tilde G(x,\mu,\cdot)\in C_b^\beta$. Thus, by the regularity of  the averaged functions in the classical multi-scale SDEs (see e.g. Corollary \ref{cor} or \cite[Lemma 3.2]{RX2}), we have
	\begin{align*}
	\bar f(\cdot,\mu)=\int_{\mR^{d_2}} f(\cdot,\mu,y,\tilde\zeta^\mu)\zeta^{\cdot,\mu}(\dif y)\in C_b^\a.
	\end{align*}
	The above argument is not suitable to study the regularity of $\bar f$ with respect to $\mu$  as the regularity of the coefficients $\tilde F(x,\cdot,y), \tilde G(x,\cdot,y)$ are unknown. To prove  $\bar f(x,\cdot)\in C_b^{(2,\a)}$, it is enough to show  that for every $x\in\mR^{d_1}$, $\psi_1\in C_p^\beta(\mR^{d_2})$ and $\psi_2\in C_p^{(2,\beta)}(\sP_2(\mR^{d_2}))$, we have
	$$
	\hat\psi_1(\mu):=\int_{\mR^{d_2}}\psi_1(y)\zeta^{x,\mu}(\dif y)\in C_b^{(2,\a)}\quad \text{and}\quad\hat\psi_2(\mu):=\psi_2(\tilde\zeta^{\mu})\in C_b^{(2,\a)}.
	$$
	For these, we use the approximation argument.  Instead of (\ref{Ynn}), we consider that for fixed $\mu$,
	\begin{align*}
	&\dif \hat Y_s^{n,(x,\mu)}=F\Big(x,\mu,\hat Y_s^{n,(x,\mu)},\int_{\mR^{d_1}}\hat \zeta_{n-1}^{\tilde x,\mu}\mu(\dif\tilde x)\Big)\dif s\no\\
	&\qquad\qquad\qquad+G\Big(x,\mu,\hat Y_s^{n,(x,\mu)},\int_{\mR^{d_1}}\hat \zeta_{n-1}^{\tilde x,\mu}\mu(\dif\tilde x)\Big)\dif \tilde W_s,
	\end{align*}
	where $\hat \zeta_{n}^{x,\mu}$ is the unique invariant measure of $\hat Y_s^{n,(x,\mu)}$. Define
	$$
	\hat\psi_{1,n}(\mu):=\int_{\mR^{d_2}}\psi_1(y)\hat\zeta_n^{x,\mu}(\dif y)\quad \text{and}\quad\hat\psi_{2,n}(\mu):=\psi_2\left(\int_{\mR^{d_1}}\hat\zeta_n^{x,\mu}\mu(\dif x)\right).
	$$
	Then, by induction and corollary \ref{cor}  we have that
	$$
	\sup_{n\geq 1}\big(\|\hat\psi_{1,n}\|_{C_b^{(2,\a)}}+\|\hat\psi_{2,n}\|_{C_b^{(2,\a)}}\big)<\infty.
	$$
	Moreover, by exactly the same procedure as in Lemma \ref{inv}, we have
	\begin{align*}
	\lim_{n\to\infty}W_\vartheta (\hat\zeta_n^{x,\mu},\zeta^{x,\mu})=0,
	\end{align*}
	which in turn implies that
	$$
	\lim_{n\to\infty}\hat\psi_{i,n}(\mu)=\hat\psi_i(\mu),~i=1,2.
	$$
	As a result, we have  $\hat\psi_i(\cdot)\in C_b^{(2,\a)}.$ The proof is finished.
\end{proof}

We have the following result.
\bl
There exists $M>0$ such that for every $n\geq 1$, we have
\begin{align}\label{ay}
\sup_{s\in[0,t]}\rho_{\a,M}(\cL_{\bar X_s^n},\cL_{\bar X_s})\leq C_t\sup_{s\in[0,t]}\Big(\rho_{\a,M}(\cL_{\bar X_s^{n-1}},\cL_{\bar X_s})+ +\tilde\rho_V(\zeta_{n-1}^{\cdot,\cL_{\bar X_s^{n-2}}},\zeta^{\cdot,\cL_{\bar X_s}})\Big).
\end{align}
where $\tilde\rho_V$ is defined by (\ref{dis}), and $C_t>0$ is a constant with $\lim_{t\to0}C_t=0$.
\el
\begin{proof}
Using Lemma \ref{diff}, we have for every $\varphi\in C_b^{(2,\a)}$,
\begin{align*}
\big|\varphi(\cL_{\bar X_t^n})-\varphi(\cL_{\bar X_t})\big|&\leq C_t\sup_{s\in[0,t]}\Big(\|\bar b_n(\cdot,\cL_{\bar X_s^{n-1}})-\bar b(\cdot,\cL_{\bar X_s})\|_\infty\\
&\qquad\qquad\qquad+\|\bar \sigma_n(\cdot,\cL_{\bar X_s^{n-1}})-\bar \sigma(\cdot,\cL_{\bar X_s})\|_\infty\Big),
\end{align*}
where $C_t>0$ is a constant with $\lim_{t\to0}C_t=0$.
By definition and as in the proof of (\ref{666}), we have
\begin{align*}
|\bar b_n(x,\mu_{n-1})-\bar b(x,\mu)|&=\bigg|\int_{\mR^{d_2}}b\Big(x,\mu_{n-1},y,\int_{\mR^{d_1}}\zeta_{n-1}^{\tilde x,\mu_{n-2}}\mu_{n-1}(\dif\tilde x)\Big)\zeta_n^{x,\mu_{n-1}}(\dif y)\\
&\quad-\int_{\mR^{d_2}}b\Big(x,\mu,y,\int_{\mR^{d_1}}\zeta^{\tilde x,\mu}\mu(\dif\tilde x)\Big)\zeta^{x,\mu}(\dif y)\bigg|\\
&\leq C_0\Big(\rho_{\a,M}(\mu_{n-1},\mu)+\rho_\cV(\zeta_n^{x,\mu_{n-1}},\zeta^{x,\mu})+\tilde\rho_\cV(\zeta_{n-1}^{\cdot,\mu_{n-2}},\zeta^{\cdot,\mu})\Big),
\end{align*}
and similarly,
\begin{align*}
|\bar \sigma_n(x,\mu_{n-1})-\bar \sigma(x,\mu)|\leq C_0\Big(\rho_{\a,M}(\mu_{n-1},\mu)+\rho_\cV(\zeta_n^{x,\mu_{n-1}},\zeta^{x,\mu})+\tilde\rho_\cV(\zeta_{n-1}^{\cdot,\mu_{n-2}},\zeta^{\cdot,\mu})\Big),
\end{align*}
which in turn imply that
\begin{align}\label{55}
\big|\varphi(\cL_{\bar X_t^n})-\varphi(\cL_{\bar X_t})\big|&\leq C_t\sup_{s\in[0,t]}\Big(\rho_{\a,M}(\cL_{\bar X_s^{n-1}},\cL_{\bar X_s})+\tilde\rho_\cV(\zeta_n^{\cdot,\cL_{\bar X_s^{n-1}}},\zeta^{\cdot,\cL_{\bar X_s}})\no\\
&\qquad+\tilde\rho_\cV(\zeta_{n-1}^{\cdot,\cL_{\bar X_s^{n-2}}},\zeta^{\cdot,\cL_{\bar X_s}})\Big).
\end{align}
Furthermore, by estimate (\ref{ax}) we have that
$$
\tilde\rho_\cV(\zeta_n^{\cdot,\cL_{\bar X_t^{n-1}}},\zeta^{\cdot,\cL_{\bar X_t}})\leq C_0\Big(\rho_{\a,M}(\cL_{\bar X_t^{n-1}},\cL_{\bar X_t})+\kappa\,\tilde\rho_V(\zeta_{n-1}^{\cdot,\cL_{\bar X_t^{n-2}}},\zeta^{\cdot,\cL_{\bar X_t}})\Big).
$$
Taking this back into (\ref{55}) yields the desired result.
\end{proof}

Now, we proceed to give:

\begin{proof}[{\bf Proof of Theorem \ref{main1} (i) and (ii)}]
Let $(X_t^\eps,Y_t^\eps)$ and $\bar X_t$ satisfy the system (\ref{sde1}) and the averaged system (\ref{ave}), respectively,   and for every $n\geq 1$, $(X^{n,\eps}_t,Y^{n,\eps}_t)$ be the solution of the approximation system (\ref{sde40}).
Then by Theorem \ref{weakn} and Corollary \ref{weaknn}, we deduce that for every $t\in[0,T]$,
\begin{align}\label{3}
\big|\varphi(\cL_{X_t^{\eps}})-\varphi(\cL_{\bar X_t})\big|&\leq \big|\varphi(\cL_{X_t^{\eps}})-\varphi(\cL_{ X^{n,\eps}_t})\big|\no\\
&\quad+\big|\varphi(\cL_{X_t^{n,\eps}})-\varphi(\cL_{\bar X_t^{n}})\big|+\big|\varphi(\cL_{\bar X_t^{n}})-\varphi(\cL_{\bar X_t})\big|\no\\
&\leq C_T\,\eps^\frac{\alpha}{2}+\big|\varphi(\cL_{X_t^{\eps}})-\varphi(\cL_{ X^{n,\eps}_t})\big|+\big|\varphi(\cL_{\bar X_t^{n}})-\varphi(\cL_{\bar X_t})\big|,
\end{align}
and
\begin{align*}
\big|\psi(\cL_{Y_t^{\eps}})-\psi(\tilde\zeta^{\cL_{\bar X_t}}) \big|&\leq \big|\psi(\cL_{Y_t^{\eps}})-\psi(\cL_{Y_t^{n,\eps}})\big|+\Big|\psi(\cL_{Y_t^{n,\eps}})
-\psi\Big(\int_{\mR^{d_1}}\zeta_n^{x,\cL_{\bar X_t^{n-1}}}\cL_{\bar X_t^n}(\dif x)\Big) \Big|\\
&\quad+\Big|\psi\Big(\int_{\mR^{d_1}}\zeta_n^{x,\cL_{\bar X_t^{n-1}}}\cL_{\bar X_t^n}(\dif x)\Big)-\psi\Big(\int_{\mR^{d_1}}\zeta^{x,\cL_{\bar X_t}}\cL_{\bar X_t}(\dif x)\Big) \Big|\\
&\leq\tilde C_T\,\eps^{\frac{\alpha}{2}}+C_0\,\e^{-\frac{\gamma_0 t}{\eps}}+\big|\psi(\cL_{Y_t^{\eps}})-\psi(\cL_{Y_t^{n,\eps}})\big|\\
&\quad+\Big|\psi\Big(\int_{\mR^{d_1}}\zeta_n^{x,\cL_{\bar X_t^{n-1}}}\cL_{\bar X_t^n}(\dif x)\Big)-\psi\Big(\int_{\mR^{d_1}}\zeta^{x,\cL_{\bar X_t}}\cL_{\bar X_t}(\dif x)\Big) \Big|,
\end{align*}
where  $C_T, \tilde C_T, C_0, \gamma_0>0$ are constants independent of $n$, and $\tilde\zeta^{\cL_{\bar X_t}}$ is defined by
$$\tilde\zeta^{\cL_{\bar X_t}}:=\int_{\mR^{d_1}}\zeta^{x,\cL_{\bar X_t}}\cL_{\bar X_t}(\dif x).$$
By the convergence of the heat kernel (see \cite[Section 5]{CF2}) of $(X_t^{n,\eps}, Y_t^{n,\eps})$ to $(X_t^{\eps},Y_t^{\eps})$, we have that for every fixed $\eps>0$,
\begin{align*}
&\lim_{n\to\infty}\Big(\big|\varphi(\cL_{ X^{n,\eps}_t})-\varphi(\cL_{X_t^{\eps}})\big|+
\big|\psi(\cL_{Y_t^{n,\eps}})-\psi(\cL_{Y_t^{\eps}})\big|\Big)=0.
\end{align*}
Below, we show  that there exists a small $t_0>0$ such that for every $t\in[0,t_0]$ and $\varphi\in C_b^{(2,\a)}(\mR^{d_1})$,
\begin{align}\label{limit}
\lim_{n\to\infty}\big|\varphi(\cL_{\bar X_t^{n}})-\varphi(\cL_{\bar X_t})\big|=0.
\end{align}
For this, let us define
$$
\sC_n(t):=\sup_{s\in[0,t]}\rho_{\a,M}(\cL_{\bar X_s^n},\cL_{\bar X_s})\quad\text{and}\quad\sD_n(t):=\sup_{s\in[0,t]}\tilde\rho_V(\zeta_n^{\cdot,\cL_{\bar X_s^{n-1}}},\zeta^{\cdot,\cL_{\bar X_s}}).
$$
Combing (\ref{ax}) and (\ref{ay}), we deduce that
\begin{align*}
\sC_n(t)\leq C_t\Big(\sC_{n-1}(t)+\sD_{n-1}(t)\Big).
\end{align*}
and
\begin{align*}
\sD_n(t)\leq C_0\Big(\sC_{n-1}(t)+\kappa\sD_{n-1}(t)\Big).
\end{align*}
Thus we have
\begin{align*}
\sC_n(t)+\kappa\sD_n(t)\leq (C_t/\kappa+C_0\kappa)\Big(\sC_{n-1}(t)+\kappa\sD_{n-1}(t)\Big).
\end{align*}
Taking $t_0$ small enough such that $C_{t_0}/\kappa+C_0\kappa<1$, then we have that for every $t\in[0,t_0]$,
$$
\lim_{n\to\infty}\(\sC_n(t)+\kappa\sD_n(t)\)=0.
$$
This together with (\ref{limit}) also implies that for $t\in[0,t_0]$,
\begin{align*}
\lim_{n\to\infty}\Big|\psi\Big(\int_{\mR^{d_1}}\zeta_n^{x,\cL_{\bar X_t^{n-1}}}\cL_{\bar X_t^n}(\dif x)\Big)-\psi\Big(\int_{\mR^{d_1}}\zeta^{x,\cL_{\bar X_t}}\cL_{\bar X_t}(\dif x)\Big) \Big|=0.
\end{align*}
For general $t\in[0,T]$, the above convergence follows by the semigroup property. The proof is finished.
\end{proof}

\subsection{Weak convergence implies the strong convergence}

Let us first point out that in Theorem \ref{weakn}, we can  obtain simultaneously the strong convergence of $X_t^{n,\eps}$ to $\bar X_t^n$ as direct results of {\it(i)} in Theorem \ref{non-aut}. Namely, we can get that for every $n\geq 1$,
\begin{align*}
\mE|X_t^{n,\eps}-\bar X_t^n|^2\leq C_t\,\eps^{\a\wedge1},
\end{align*}
where $C_t>0$ is independent of $n$, and $\bar X_t^n$ satisfy equation (\ref{avenn}).
Then in order to prove the strong convergence (\ref{es3}) of $X_t^\eps$ to $\bar X_t$ in Theorem \ref{main1} (iii), we can deduce as in (\ref{3}) that
\begin{align*}
\mE|X_t^\eps-\bar X_t|^2&\leq \mE|X_t^\eps-X_t^{n,\eps}|^2+\mE|X_t^{n,\eps}-\bar X_t^n|^2+\mE|\bar X_t^{n}-\bar X_t|^2\\
&\leq C_t\,\eps^{\a\wedge1}+\mE|X_t^\eps-X_t^{n,\eps}|^2+\mE|\bar X_t^{n}-\bar X_t|^2.
\end{align*}
Thus for every $\eps>0$, we need to show the {\bf strong convergence} of $X_t^{n,\eps}$ and $\bar X_t^{n}$ to $X_t^\eps$ and $\bar X_t$ as $n\to\infty$, respectively. But due to the low regularity assumptions on the coefficients of the systems (only H\"older continuous), even for the proof of the strong convergence of $X_t^{n,\eps}$ to $X_t^\eps$ for every fixed $\eps>0$, we shall need to use Zvonkin's argument to transform the original systems into a new one with better coefficients (see e.g. \cite{C}), which is quite complicated.
To avoid this, we make use of the idea that to transform the non-linear system into a non-autonomous linear system again. It turns out to be quite easy, as we shall see,  that the strong convergence of $X_t^\eps$ in the averaging principle of the non-linear stochastic system (\ref{sde1}) follows directly from the convergence of the distribution of $X_t^\eps$ and $Y_t^\eps$.

For every $\eps>0$, let us  define
\begin{eqnarray}
\begin{split}\label{bfeps}
b_\eps(t,x,y)&:=b(x,\cL_{X_t^\eps},y,\cL_{Y_t^\eps}),\\ \sigma_\eps(t,x)&:=\sigma(x,\cL_{X_t^\eps},\cL_{Y_t^\eps}),\\ F_\eps(t,x,y)&:=F(x,\cL_{X_t^\eps},y,\cL_{Y_t^\eps}),\\
 G_\eps(t,x,y)&:=G(x,\cL_{X_t^\eps},y,\cL_{Y_t^\eps}).
\end{split}
\end{eqnarray}
Then, the system (\ref{sde1}) can be rewritten as
\begin{equation} \label{sde3}
\left\{ \begin{aligned}
&\dif X^{\eps}_t =b_\eps(t,X^{\eps}_t,Y^{\eps}_t)\dif t
+\sigma_\eps(t,X^{\eps}_t)\dif W^1_t,\qquad\qquad\quad\, X^{\eps}_0=\xi,\\
&\dif Y^{\eps}_t =\frac{1}{\eps}F_\eps(t,X^{\eps}_t,Y^{\eps}_t)\dif t+\frac{1}{\sqrt{\eps}}G_\eps(t,X^{\eps}_t,Y^{\eps}_t)\dif W_t^2,\quad   Y^{\eps}_0=\eta.
\end{aligned} \right.
\end{equation}
This system is exactly the form of (\ref{sde22}), and the convergence of distributions of $X_t^\eps$ and $Y_t^\eps$ obtained in the previous subsection imply the convergence of the coefficients $b_\eps,\sigma_\eps, F_\eps$ and $G_\eps$. Thus, the strong convergence of $X_t^\eps$ can be obtained by the strong convergence of the non-autonomous system (\ref{sde22}) obtained in Theorem \ref{non-aut}  {\it (i)} directly.

\vspace{2mm}
We proceed to give:

\begin{proof}[{\bf Proof of Theorem \ref{main1} (iii)}]
For every $x\in\mR^{d_1}$, $y\in\mR^{d_2}$ and $t>0$, let us define
\begin{eqnarray}
\begin{split}\label{bfhat}
\hat b(t,x,y)&:=b\big(x,\cL_{\bar X_t},y,\tilde\zeta^{\cL_{\bar X_t}}\big),\\
\hat \sigma(t,x)&:=\sigma\big(x,\cL_{\bar X_t},\tilde\zeta^{\cL_{\bar X_t}}\big),\\
\hat F(t,x,y)&:=F\big(x,\cL_{\bar X_t},y,\tilde\zeta^{\cL_{\bar X_t}}\big),\\
\hat G(t,x,y)&:=G\big(x,\cL_{\bar X_t},y,\tilde\zeta^{\cL_{\bar X_t}}\big),
\end{split}
\end{eqnarray}
where $\bar X_t$ satisfies the averaged equation (\ref{ave}), and $\tilde \zeta^\mu$ is defined by (\ref{zmu}).
By the assumptions that $b\in C_p^{\a,(2,\a),\b,(2,\b)}$, the definitions (\ref{bfeps}) and (\ref{bfhat}), and using the  convergence of the distributions of $X_t^\eps$ and $Y_t^\eps$ obtained in estimates (\ref{es1}) and (\ref{es2}), we have that there exists a constant $C_t>0$ independent of $\eps$ such that
\begin{align*}
\big|b_{\eps}(t,x,y)-\hat b(t,x,y)\big|&=\Big|b\big(x,\cL_{X_t^\eps},y,\cL_{Y_t^\eps}\big)-b\big(x,\cL_{\bar X_t},y,\tilde\zeta^{\cL_{\bar X_t}}\big)\Big|\\
&\leq C_t(1+|y|^p)\big(\eps^{\frac{\a}{2}}+e^{-\frac{\gamma_0 t}{\eps}}\big),
\end{align*}
which in turn implies that
\begin{align*}
&\|b_\eps(t,\cdot,\cdot)-\hat b(t,\cdot,\cdot)\|_{L_p^\infty}\leq C_t\big(\eps^{\frac{\a}{2}}+e^{-\frac{\gamma_0 t}{\eps}}\big).
\end{align*}
Meanwhile, we have that $\hat b(t,\cdot,\cdot)\in C_p^{\a,\b}$, and by Lemma \ref{lem3} we deduce that  for every $t_1,t_2\in\mR_+$,
\begin{align*}
\big|\hat b(t_1,x,y)-\hat b(t_2,x,y)\big|&=\Big|b\big(x,\cL_{\bar X_{t_1}},y,\tilde\zeta^{\cL_{\bar X_{t_1}}}\big)-b\big(x,\cL_{\bar X_{t_2}},y,\tilde\zeta^{\cL_{\bar X_{t_2}}}\big)\Big|\\
&\leq C_0(1+|y|^p)\cW_2(\cL_{\bar X_{t_1}},\cL_{\bar X_{t_2}})^\a\leq C_0(1+|y|^p)|t_1-t_2|^\frac{\a}{2},
\end{align*}
As a result, we get $\hat b\in C_p^{\a/2,\a,\b}$. Similarly, we have that
\begin{align*}
\|\sigma_\eps(t,\cdot)-\hat \sigma(t,\cdot)\|_{L^\infty}&+\|F_\eps(t,\cdot,\cdot)-\hat F(t,\cdot,\cdot)\|_{L_p^\infty}\\
&+\|G_\eps(t,\cdot,\cdot)-\hat G(t,\cdot,\cdot)\|_{L_p^\infty}\leq C_t\big(\eps^{\frac{\a}{2}}+e^{-\frac{\gamma_0 t}{\eps}}\big),
\end{align*}
and $\hat F,\hat G\in C_p^{\a/2,\a,\b}$ and $\hat \sigma\in C_b^{\a/2,1}$.
As a result of estimate (\ref{st-non}) in Theorem \ref{non-aut}, we obtain
\begin{align*}
\sup_{t\in[0,T]}\mE\big|X_t^\eps-\bar{\hat X}_t\big|^2&\leq C_T\bigg(\eps^{\a\wedge1}+\int_0^T\Big[\|b_\eps(s,\cdot,\cdot)-\hat b(s,\cdot,\cdot)\|^2_{L_p^\infty}+\|\sigma_\eps(s,\cdot)-\hat \sigma(s,\cdot)\|^2_{L^\infty}\no\\
&\qquad\quad+\|F_\eps(s,\cdot,\cdot)-\hat F(s,\cdot,\cdot)\|^2_{L_p^\infty}+\|G_\eps(s,\cdot,\cdot)-\hat G(s,\cdot,\cdot)\|^2_{L_p^\infty}\Big]\dif s\bigg)\\
&\leq C_T\,\eps^{\a\wedge1}+C_T\int_0^T\big(\eps^\a+e^{-\frac{2\gamma_0 s}{\eps}}\big)\dif s\leq C_T\,\eps^{\a\wedge1},
\end{align*}
where $\bar{\hat X}_t$ satisfy the equations (\ref{ave1}) with the coefficients $\hat b, \hat \sigma, \hat F$ and $\hat G$ given by (\ref{bfhat}). Hence, we need only prove that $\bar X_t=\bar{\hat X}_t$, where $\bar X_t$ satisfies (\ref{ave}). Applying (\ref{frozen1}) and (\ref{bfhat}) we have
\begin{align*}
\dif\hat Y_s^{t,x}&=\hat F(t,x,\hat Y_s^{t,x})\dif s+\hat G(t,x,\hat Y_s^{t,x})\dif \hat W_s^2\\
&=F(x,\cL_{\bar X_t},\hat Y_s^{t,x},\tilde \zeta^{\cL_{\bar X_t}})\dif s+G(x,\cL_{\bar X_t},\hat Y_s^{t,x},\tilde \zeta^{\cL_{\bar X_t}})\dif \hat W_s^2.
\end{align*}
Note that in the above equation, the parameter $t$ is fixed. This together with the frozen equation (\ref{frozen}) yields that $\zeta^{t,x}(\dif y)=\zeta^{x,\cL_{\bar X_t}}(\dif y)$. Consequently, we deduce that
\begin{align*}
\bar{\hat b}(t,x)&=\int_{\mR^{d_2}}\hat b(t,x,y)\zeta^{t,x}(\dif y)\\
&=\int_{\mR^{d_2}} b(x,\cL_{\bar X_t},y,\tilde\zeta^{\cL_{\bar X_t}})\zeta^{x,\cL_{\bar X_t}}(\dif y)=\bar b(x,\cL_{\bar X_t})
\end{align*}
and
\begin{align*}
\hat \sigma(t,x)=\sigma(x,\cL_{\bar X_t},\tilde\zeta^{\cL_{\bar X_t}})=\bar \sigma(x,\cL_{\bar X_t}),
\end{align*}
which further imply that
\begin{align*}
\dif \bar X_t&=\bar b(\bar X_t,\cL_{\bar X_t})\dif t+\bar \sigma(\bar X_t,\cL_{\bar X_t})\dif W^1_t\\
&=\bar{\hat b}(t,\bar X_t)\dif t+\hat\sigma(t,\bar X_t)\dif W^1_t.
\end{align*}
In view of the strong uniqueness of the solution to (\ref{ave1}), we have $\bar X_t=\bar{\hat X}_t$. Thus, the proof is finished.
\end{proof}

\end{document}